\documentclass[11pt,a4paper]{article}

\usepackage[english,francais]{babel} 
\usepackage{ucs}
\usepackage[utf8]{inputenc}

\RequirePackage{color} 

\RequirePackage{amsmath,amssymb,amsfonts,bbm,latexsym,mathrsfs}

\usepackage{graphicx}

\newcommand{\bu}{{\bf u}}

\newcommand{\MM}{{\cal M}}
\newcommand{\mm}{\mathtt{m}} 
\newcommand{\dis}{{\rm dis\, }}
\newcommand{\dw}{\mathtt{d}_{\rm GHP}} 
\newcommand{\bA}{{\bf A}}
\newcommand{\bN}{{\bf N}}
\newcommand{\mw}{\mathbb{M}_{\rm w}} 
\newcommand{\X}{{\rm X}}
\newcommand{\pmw}{{\rm P}\mw} 
\newcommand{\K}{\mathbb{K}}
\newcommand{\CC}{{\cal C}}

\newcommand{\ii}{{\rm i}} 
\newcommand{\dgh}{\mathtt{d}_{{\rm GH}}}

\newcommand{\be}{{\bf e}}
\newcommand{\bq}{{\bf q}}
\newcommand{\bv}{{\bf v}} 
\newcommand{\BWN}{{\bf WN}} \newcommand{\bW}{{\bf W}}
\newcommand{\bm}{{\bf m}} 
\newcommand{\br}{{\bf r}}
\newcommand{\bl}{{\bf l}} 
\newcommand{\bz}{{\bf z}} 
\newcommand{\bx}{{\bf x}} 
\newcommand{\by}{{\bf y}}
\newcommand{\bC}{{\bf C}}
\newcommand{\bS}{{\bf S}}
\newcommand{\bP}{{\bf P}} 
\newcommand{\bQ}{{\bf Q}}
\newcommand{\bE}{{\bf E}} 
\newcommand{\RR}{{\cal R}} 
\newcommand{\mg}{\mathfrak{m}} 
\newcommand{\lgo}{\mathfrak{l}} 
\newcommand{\eg}{\mathfrak{e}} 
\newcommand{\fgo}{\mathfrak{f}}

\newcommand{\R}{\mathbb{R}}
\newcommand{\N}{\mathbb{N}} 
\newcommand{\M}{\mathbb{M}}

\newcommand{\pmc}{{\rm P}\mathbb{M}} 
\newcommand{\Z}{\mathbb{Z}}

\newcommand{\C}{\mathbb{C}} 
\newcommand{\FF} { {\cal F }} 
\def\build#1_#2^#3{\mathrel{ \mathop{\kern 0pt#1}\limits_{#2}^{#3}}}

\def\cq{$\hfill \square$}
\def\un{\underline}
\def\d{{\rm d}}
\def\n{{\cal N}}
\def\W{{\cal W}}

\def\eps{\varepsilon}
\def\QQ{{\sf Q}}
\def\DSN{{\sf DS}}
\def\SN{{\sf S}}
\def\WN{{\sf WN}}
\def\PN{{\sf PN}}
\def\bM{{\bf M}}
\def\LM{{\sf LM}}
\def\CLM{{\sf CLM}}
\def\BLM{{\bf LM}}
\def\ov{\overline}
\def\wh{\widehat}
\def\wt{\widetilde}
\def\rem{\noindent{\bf Remark. }}
\def\rems{\noindent{\bf Remarks. }}
\def\proof{\noindent{\bf Proof. }}
\def\diam{{\rm diam\,}}
\def\supp{{\rm supp\,}}
\def\cov{{\rm Cov\,}}
\def\aut{{\rm Aut\,}}

\newcommand{\ind}{\mathbbm{1}}

\newtheorem{thm}{Theorem}
\newtheorem{lmm}{Lemma}
\newtheorem{prp}{Proposition}
\newtheorem{defn}{Definition}
\newtheorem{crl}{Corollary}

\usepackage{fancyhdr}

\usepackage[outerbars]{changebar} 

\title{Tessellations of random maps of arbitrary genus
\\(Mosaïques sur des cartes aléatoires en genre arbitraire)}

\author{Gr\'egory Miermont\thanks{This research was supported by
    the Fondation des Sciences Mathématiques de Paris and by the ANR grant ANR-08-BLAN-0190.} \\
  CNRS \& DMA, \'Ecole Normale Sup\'erieure}

\begin{document}

\selectlanguage{english}

\maketitle

\begin{abstract}
We investigate Voronoi-like tessellations of bipartite
quadrangulations on surfaces of arbitrary genus, by using a natural
generalization of a bijection of Marcus and Schaeffer allowing one to
encode such structures by labeled maps with a fixed number of
faces. We investigate the scaling limits of the latter. Applications
include asymptotic enumeration results for quadrangulations, and
typical metric properties of randomly sampled quadrangulations. In
particular, we show that scaling limits of these random
quadrangulations are such that almost every pair of points are linked
by a unique geodesic.

\bigskip

\selectlanguage{francais}

Nous examinons les propriétés de mosaïques de type Voronoi sur des
quadrangulations biparties de genre arbitraire. Ceci est rendu
possible par une généralisation naturelle d'une bijection de Marcus et
Schaeffer, permettant de décrire ces mosaïques par des cartes
étiquetées avec un nombre fixé de faces, dont nous déterminons les
limites d'échelle. Parmi les applications de ces résultats, figurent
le comptage asymptotique des quadrangulations, ainsi que
des propriétés métriques typiques de quadrangulations
choisies au hasard. En particulier, nous montrons que les limites
d'échelles de quadrangulations aléatoires sont telles que presque
toute paire de points est liée par un unique chemin géodésique.
\end{abstract}

\selectlanguage{english}

\noindent{\bf Keywords.} Random maps, scaling limits, random snakes,
asymptotic enumeration, geodesics. 

\bigskip

\noindent{\bf MSC.} 60C05 ; 05C30 ; 60F05

\section{Introduction}\label{sec:introduction}

\subsection{Motivation}\label{sec:maps-surf-theor}

A map is a graph embedded into a compact orientable surface without
boundary, yielding a cell decomposition of the surface, and considered
up to orientation-preserving homeomorphisms. Random maps are
considered in the physics literature on quantum gravity as discretized
versions of random surfaces \cite{ADJ}. This approach allows one to
perform computations of certain integrals with respect to an
ill-defined measure on surfaces, by approximating them by finite sums
over maps. From a mathematical perspective, this leads to the
stimulating problem of existence of a measure on compact surfaces
arising as the scaling limit of, say, uniform random triangulations of
the sphere (or the compact orientable surface of genus $g$) with a
large number of faces.

Important progress has been made in this direction in the recent
years, thanks to bijective approaches initiated in Schaeffer's thesis
\cite{schaeffer98}. They allowed Chassaing \& Schaeffer \cite{CSise}
to obtain the scaling limits for the radius and profile of a uniform
rooted planar quadrangulation with $n$ faces, considered as a metric
space by endowing its vertices with the usual graph distance. In
particular, they showed that typical graph distances between vertices
of such random quadrangulations are of order $n^{1/4}$.
Generalizations of this result have been obtained for much more
general families of random maps by Marckert, Weill and the author in
\cite{jfmgm05,mierinv,weill06,mierweill}, relying on generalizations
of Schaeffer's bijection by Bouttier, Di Francesco \& Guitter
\cite{BdFGmobiles}.

An important step has been accomplished by Le Gall \cite{legall06},
who showed that scaling limits of these random quadrangulations,
considered in the Gromov-Hausdorff sense \cite{gromov99}, must be
homeomorphic to a metric quotient of the so-called Brownian Continuum
Random Tree of Aldous \cite{aldouscrt91}. This random quotient was
introduced under the name of {\em Brownian map} by Marckert \&
Mokkadem \cite{MM05}, who proved convergence of random
quadrangulations towards this limit in a sense yet different from
convergence in distribution in the Gromov-Hausdorff
topology. Subsequently, Le Gall \& Paulin \cite{lgp} inferred that the
Brownian map is homeomorphic to the two-dimensional sphere. At the
present stage, it is however not known whether the scaling limit is
uniquely determined as the Brownian map, which would lead to a
satisfactory answer to the above mentioned problem.

A natural idea is to characterize the limit through its
``finite-dimensional marginal distributions''. From the point of view
of metric measure spaces, this would mean to characterize the laws of
mutual distances between an arbitrary number of randomly sampled
points. This, however, seems hard to obtain. The approach of the
present work, which deals with combinatorial, probabilistic and
geometric aspects of maps, was initially motivated by another notion
of finite-dimensional marginals, which takes into account the planar
structure of the graphs that are considered, and roughly consists in
studying the Voronoi tessellation on the map with sources taken at
random.

Our central combinatorial result is a natural generalization of
Schaeffer's bijective construction, to be introduced in Section
\ref{sec:coding-delay-voron}. We do not restrict ourselves to the
planar case, as the most natural framework is to consider maps of
arbitrary (but fixed) genus $g$, so our bijection is really a
generalized version of the Marcus-Schaeffer bijection in arbitrary
genus \cite{MaSc01,ChMaSc}.
Let us explain briefly the idea of our construction. Recall that the
Marcus-Schaeffer bijection encodes a genus-$g$ bipartite
quadrangulation with a distinguished vertex by a labeled map of genus
$g$ with one face, in such a way that the labels keep track of graph
distances to the distinguished vertex in the quadrangulation. Our
construction encodes a genus-$g$ bipartite quadrangulation with $k$
distinguished vertices by a labeled map of genus $g$ with $k$ faces
(Theorem \ref{BIJ} and Corollary \ref{sec:rooting-1}). In some sense,
the faces of the encoding map correspond to the Voronoi tiles in the
quadrangulation, with sources at the $k$ distinguished vertices, and
the labels allow to keep track of graph distance of a vertex to the
source of the Voronoi tile it belongs to. An important fact is that
there is some flexibility in the construction, allowing one to
consider a generalized Voronoi tessellation in which the distances to
the $k$ sources are measured after addition of a {\em delay} depending
on the source.

The resulting set of labeled maps with a fixed number of faces is
combinatorially much simpler than that of bipartite quadrangulations,
and we are able to study their scaling limits. This is done in Section
\ref{sec:discr-cont-path} (Theorem \ref{sec:cont-meas-label-1} and
Propositions \ref{sec:exceptional-cases-1},
\ref{sec:exceptional-cases-2}), after showing how to encode labeled
maps by appropriate processes in Section
\ref{sec:struct-label-maps}. Sections \ref{sec:appl-asympt-enum} and
\ref{sec:proof-main-results} give three applications of these scaling
limit results for labeled maps, to combinatorial and geometric
properties of bipartite quadrangulations

The first application is a new derivation of known asymptotic
enumeration results (Theorem \ref{sec:main-results}), to which the
short Section \ref{sec:appl-asympt-enum} is devoted. These are
initially due to Bender \& Canfield \cite{BenCan86}, who obtained the
asymptotic number of rooted maps of genus $g$ with $n$
edges\footnote{This is equivalent to our result, since such maps are
  in one-to-one correspondence with bipartite quadrangulations of same
  genus with $n$ faces \cite[Proposition 1]{MaSc01}.} by recursive
decomposition methods, of a very different nature from our bijective
study. These results have also been obtained in the recent work of
Chapuy, Marcus \& Schaeffer \cite{ChMaSc}, who completed the exact
enumeration of maps of genus $g$ that was initiated in
\cite{MaSc01}. The starting point of our study (using the
Marcus-Schaeffer bijection, which corresponds to the case $k=1$ of our
bijection) is the same as in \cite{ChMaSc}, but the rest of our
approach is different, the arguments of
Sect.\ \ref{sec:discr-cont-path} being of a more probabilistic nature.

Our other results concern the metric structure of randomly sampled
bipartite quadrangulations. Let us indicate briefly what we mean by
``random'' in this paper. Most the articles on the topic have focused
on scaling limits as $n\to\infty$ of random families of maps
conditioned to have $n$ faces or $n$ vertices. We prefer to randomize
the number of faces, by using natural $\sigma$-finite measures on
bipartite quadrangulations of fixed genus $g$, called the Boltzmann
measures $\QQ_g$ (the term being inspired from \cite{angelpeeling},
see also \cite{jfmgm05,mierinv}), and which are obtained by assigning
an appropriate weight to the faces. In the Physics terminology, these
measures correspond to the so-called {\em grand-canonical} measures,
while the measures with fixed number of faces are the {\em
  microcanonical} measures. These measures are used to define natural
``Boltzmann-Gibbs'' distributions on quadrangulations, depending on an
inverse temperature parameter $\beta$ that allows one to tune the
average size of the quadrangulation and take scaling limits.

The existence of these scaling limits, considered with respect to the
Gromov-Hausdorff-Prokhorov topology on metric measure spaces, and
stated in Theorem \ref{sec:mains-results-1}, is then obtained from the
study of Section \ref{sec:discr-cont-path}. This generalizes to
arbitrary genera the fact that planar quadrangulations admit scaling
limits for the Gromov-Hausdorff topology, as shown in Le Gall
\cite{legall06}.  The theory of metrics on weighted metric spaces that
is needed here is developed in Section \ref{sec:grom-hausd-weight}.

The main result of the present paper, Theorem
\ref{sec:mains-results-2}, gives qualitative information on the metric
structure of the scaling limits of random bipartite
quadrangulations. We show that these scaling limits are geodesic
weighted metric spaces, in which two typical points are linked by a
unique geodesic. The idea is to use our bijection for quadrangulations
with $k=2$ marked vertices, in order to show that, if $x,y$ are
typical points in the limiting random metric space $(X,d)$ arising as
a scaling limit of random bipartite quadrangulations, the intersection
of a geodesic between $x$ and $y$ and a geometric locus of the form
$$\{z\in X:d(x,z)-d(y,z)=D\}\, ,\qquad D\in [-d(x,y),d(x,y)]\, ,$$ is
a.s.\ reduced to a single point.

Theorem \ref{sec:mains-results-2} is strongly reminiscent of the
combinatorial results of the recent paper by Bouttier \& Guitter
\cite{BoGu08a}. They show that in the planar case $g=0$, two typical
vertices in a random quadrangulation with $n$ faces cannot be linked
by two geodesics which are ``very far'' in the scale $n^{1/4}$, as
$n\to\infty$. Our result confirms that geodesics between typical
points become unique in the scaling limit, answering one of the
questions raised in Section 5 of \cite{BoGu08a}, though in our
slightly different setting of Boltzmann-distributed maps. We also
mention the recent work by Le Gall \cite{legall08} on similar topics
as the present paper and \cite{BoGu08a}, and providing results on the
exceptional geodesics as well.

In this work, the cardinality of a set $A$ is denoted by $|A|$. 

\subsection{Embedded graphs and maps}\label{sec:embedded-graphs-maps-1}

Let us introduce some formalism for embedded graphs and maps, which we
partly borrow from \cite{therryYM}.

\medskip

\noindent{\bf Graphs on surfaces. }Let $S$ be a compact connected
orientable two-dimensional surface without boundary. It is well-known
that such surfaces are characterized, up to homeomorphism, by an
integer $g\geq 0$, called the genus of $S$, and that the topology of
$S$ is that of the connected sum $\mathbb{S}_g$ of $g$ tori. The
surface of genus $0$ is the two-dimensional sphere.

A {\em half-edge}, or oriented edge $e$ in $S$ is a continuous path
$c:[0,1]\to S$, considered up to continuous increasing
reparametrization, which is either injective on $[0,1]$, or is
injective on $[0,1)$ and satisfies $c(0)=c(1)$. In the latter case,
$e$ is called a loop. The order on $[0,1]$ induces a natural
orientation of the half-edge $e$, and its start and end points are
defined as $e^-=c(0)$ and $e^+=c(1)$, a definition that depends only
on the half-edge $e$ and not on the function $c$. Similarly, the image
of $e$ is defined as ${\rm Im}(e)=c([0,1])$, and the interior of $e$
is $c((0,1))$.  The reversal $\ov{e}$ of the half-edge $e$ is defined
as the reparametrization class of $c(1-\cdot)$, so that
$\ov{e}^+=e^-,\ov{e}^-=e^+$. An {\em edge} is a pair of the form
$\be=\{e,\ov{e}\}$ where $e$ is a half-edge.

An embedded graph in $S$, or simply a graph on $S$, is a pair
$(V,\bE)$, where
\begin{itemize}
\item $\bE$ is a non-empty, finite set of edges such that  two
  distinct edges intersect, if at all, only at some of their endpoints
\item $V$ is the set of vertices, i.e.\ of end-points of elements of
  $\bE$.
\end{itemize}

We let $E$ be the set of half-edges corresponding to $\bE$. The set $E$
has even cardinality and comes with the involution $e\mapsto
\ov{e}$. An orientation of the edges is a choice of one half-edge
inside each edge to form a set $E_{1/2}\subset E$ with same
cardinality as $\bE$.

A graph on $S$ determines faces, which are the connected components of
the complement of $V\cup \bigcup_{e\in E}{\rm Im}(e)$ in $S$. Let $F$
be the set of faces of the graph.

\medskip

\noindent{\bf Maps.} 
We say that the embedded graph $(V,\bE)$ is a {\em map} if all faces
are simply connected (hence determining a cell complex structure on
$S$). A map is necessarily a connected graph by \cite[Lemma
1.5]{therryYM}.

The degree of a face in a map is the number of edges it is incident
to, where it is understood that edges that are incident twice to the
same face should be counted twice. By Euler's formula,
\begin{equation}\label{euler}
|V|-|\bE|+|F|=\chi(g)\, ,
\end{equation} 
where $\chi(g)=2-2g$ is the Euler characteristic of $S$.

A map is rooted if it comes with a distinguished half-edge $e_*\in E$.

We declare two (rooted) maps $\bm,\bm'$ on surfaces $S,S'$ to be
isomorphic if there exists a homeomorphism $S\to S'$ that preserves
the orientation, and sends the vertices and edges of $\bm$ to those of
$\bm'$ (and maps the root of $\bm$ to that of $\bm'$). Note that the number
and degrees of vertices, faces, edges and genus are preserved under
such transformations. We will systematically identify maps that are
isomorphic, making the set of maps a countable set.

We will usually denote (equivalence classes of) maps by bold lowercase
letters $\bm,\bq,\ldots$. The standard notation for the root of a
rooted map will be $e_*$, and it should be clear according to the
context to which map $e_*$ refers. Also, when dealing with maps, it
will often be understood that we are working with a particular graph
belonging to the class of that map. This allows to consider the set of
(half-)edges, vertices and faces of a (isomorphism class of a) map
$\bm$, denoted by $E(\bm),\bE(\bm),V(\bm),F(\bm)$, or simply
$E,\bE,V,F$ when there is no ambiguity.

A {\em quadrangulation} is a map whose faces all have degree $4$.  A
map is called {\em bipartite} if its vertices can be colored in black
or white, in such a way that two neighboring vertices are assigned
different colors. We let $\bQ_g$ be the set of rooted bipartite
quadrangulation of genus $g$ (planar quadrangulations are all
bipartite, but it is not the case in higher genus).
For every $\bq\in \bQ_g$ we
have $4|F(\bq)|=2|\bE(\bq)|=|E(\bq)|$ and $|V(\bq)|=|F(\bq)|+\chi(g)$
by (\ref{euler}).

\medskip

\noindent{\bf Distances in maps. }For a map $\bm$, and two vertices
$x,y$, we say that the half-edges $e_1\ldots,e_n$ form a {\em chain}
from $x$ to $y$ if $e_1^-=x,e_n^+=y$ and $e_i^+=e_{i+1}^-$ for $1\leq
i\leq n-1$.  We define the graph distance $d_\bm(x,y)$ as the minimal
$n$ such that there exists a chain of $n$ half-edges from $x$ to $y$,
and a chain from $x$ to $y$ with length $d_\bm(x,y)$ is called {\em
  geodesic}.

\subsection{Spaces of (weighted) metric spaces}\label{sec:spac-weight-metr} 

The main motivation of this paper is to study the metric structure of
large random maps, endowed with their graph distance. This needs to
introduce some concepts of metric geometry, essentially originating
from Gromov's work \cite{gromov99}, which have developed a growing
interest in the probabilist community, starting from works of Evans
and coauthors. See
\cite{evpiwin,evanswinter,GPW06,duqlegprep,legall06}. All results that
are needed here will be presented in Section
\ref{sec:grom-hausd-weight}.

A weighted metric space is a triple $(X,d,\mu)$ where $(X,d)$ is
metric space and $\mu$ a Borel probability distribution on
$(X,d)$. Two metric spaces are isometric if there exists a bijective
isometry from $X$ onto $X'$. The isometry class of the space $(X,d)$
is denoted by $[X,d]$. Two weighted metric spaces $(X,d,\mu)$ and
$(X',d',\mu')$ are said to be isometry-equivalent if there exists a
bijective isometry $\phi:X\to X'$ such that $\phi_*\mu=\mu'$. The
isometry-equivalence class of $(X,d,\mu)$ is denoted by $[X,d,\mu]$.

We let $\M$ be the set of isometry classes of compact metric spaces,
and $\mw$ the set of isometry-equivalence classes of compact weighted
metric spaces. As we will see in Section \ref{sec:grom-hausd-weight},
it is possible to endow $\mw$ with a distance $\dw$, the
Gromov-Hausdorff-Prokhorov distance, that makes $(\mw,\dw)$ a Polish
space, see Theorem \ref{sec:weight-grom-hausd} and Proposition
\ref{sec:relat-with-prev-1}. The topology on $\mw$ will be implicitly
that inherited from $\dw$. If $a>0$, and $\X=[X,d]\in \M$
(resp.\ $\X=[X,d,\mu]\in \mw$), we let $a\X=[X,ad]$
(resp.\ $a\X=[X,ad,\mu]$).

Also, a metric space $(X,d)$ is called a {\em geodesic metric space}
if for which every pair of points $x,y\in X$, there exists a geodesic
path between $x$ and $y$, i.e.\ an isometry $\gamma_{xy}:[0,d(x,y)]\to
X$, with $\gamma_{xy}(0)=x,\gamma_{xy}(d(x,y))=y$. In particular,
geodesic metric spaces are arcwise connected. 
For compact metric spaces, being a geodesic space is equivalent to
being a {\em path metric space}, i.e.\ an arcwise connected metric
space $(X,d)$ such that the distance $d(x,y)$ is achieved as the
minimum over all continuous paths $c:[0,1]\to X$ with $c(0)=x,c(1)=y$,
of the length of $c$, defined as 
$${\rm length}(c)=\sup \left\{
  \sum_{i=0}^{n-1}d(c(t_i),c(t_{i+1})): n\geq 1, 0=t_0<t_1<\ldots
  <t_n=1\right\} \, .$$ We let $\pmc$ be the set of isometry classes
of compact path (or geodesic) metric spaces.

\subsection{Main results}\label{sec:mains-results}

Our core combinatorial result (Theorem \ref{BIJ}) is a generalization
of the Marcus-Schaeffer bijection, encoding $k$-pointed bipartite
quadrangulations (with delays between sources) with the help of
labeled maps with $k$ faces. The central statement in this paper is
Theorem \ref{sec:cont-meas-label-1} (and its exceptions, Propositions
\ref{sec:exceptional-cases-1} and \ref{sec:exceptional-cases-2}),
which allows us to take scaling limits of this simpler class of maps.
We do not reproduce the exact statements of these theorems in this
paragraph, as important notational background is needed.

We give several applications of these results to various
(combinatorial and probabilistic) aspects of maps. The first is a
counting result, giving the following asymptotic behavior for the
cardinality of the set $\bQ_g^n$ of bipartite quadrangulations of
genus $g$ with $n$ faces, consistently with \cite{BenCan86,ChMaSc}:
\begin{thm}\label{sec:main-results}
For any $g\geq 0$, it holds that
$$\left|\bQ_g^n\right|\build\sim_{n\to\infty}^{} C_g \, 12^n\,
n^{-5\chi(g)/4}\, ,$$ for some constant $C_g\in(0,\infty)$
defined at (\ref{eq:25}).
\end{thm}
In the case $g=0$, this result is a trivial consequence of the exact
formula (see \cite{CSise})
$$|\bQ_0^n|=\frac{2}{n+2}\, 3^n\, {\rm Cat}_n\, ,$$ where ${\rm
  Cat}_n=\frac{1}{n+1}\binom{2n}{n}$ is the $n$-th Catalan number,
yielding the value $C_0=2\pi^{-1/2}$. 

From this statement, we see that the generating function for bipartite
quadrangulations counted with respect to the number of faces has
radius of convergence $12^{-1}$, and it is natural to introduce the
{\em critical Boltzmann measure} $\QQ_g$ on $\bQ_g$ assigning mass
$12^{-|F(\bq)|}$ to the element $\bq$. Since
$|V(\bq)|=|F(\bq)|+\chi(g)$, the function $\beta\mapsto
\QQ_g(e^{-\beta V_\bq})$ is analytic in $\beta\in(0,\infty)$, where
$V_\bq=|V(\bq)|$ is the ``volume'' of $\bq$.  We define the {\em
  Boltzmann-Gibbs probability distribution} on $\bQ_g$ with inverse
temperature $\beta$ by
$$\QQ_g^{(\beta)}=\frac{\QQ_g(e^{-\beta V_\bq}\,
  \d\bq)}{\QQ_g(e^{-\beta V_\bq})}\, ,\qquad g\geq 1\, ,$$
and 
$$\QQ_0^{(\beta)}=\frac{\QQ_g(V_\bq^2e^{-\beta V_\bq}\,
  \d\bq)}{\QQ_g(V_\bq^2e^{-\beta V_\bq})}\, .$$ The reason for the
difference between the cases $g=0$ and $g\geq 1$ is technical and will
be explained in due course.  

With every $\bq\in \bQ_g$ we associate the metric space
$\X_\bq=[V(\bq),d_\bq]$ and the weighted space $\X_\bq^{\rm
  w}=[V(\bq),d_\bq,\mu_\bq]$, where
$$\mu_\bq=\frac{1}{V_\bq}\sum_{v\in V(\bq)}\delta_v$$ is the uniform
distribution on the vertices of $\bq$.

\begin{thm}\label{sec:mains-results-1}
For every $g\geq 0$ and $\beta>0$, the probability measures
\begin{equation}\label{eq:10}
\QQ^{(\beta/a)}_g\Big(\big\{\big(a^{-1/4}\X^{\rm w}_\bq,a^{-1}V_\bq
\big)\in \cdot\big\}\Big)\, , \qquad a>1
\end{equation} 
form a relatively compact family of probability distributions on
$\mw\times (0,\infty)$, endowed with the topology of weak convergence.
\end{thm}

\begin{thm}\label{sec:mains-results-2}
A limiting point $\mathscr{S}_g^{(\beta)}$ of (\ref{eq:10}), which
is a probability measure on $\mw\times (0,\infty)$, is supported on spaces
$([X,d,\mu],\mathcal{V})$ such that
\begin{enumerate}
\item $[X,d]$ is an element of $\pmc$,
\item
$\mu$ is a diffuse measure of full support,
\item 
 for $\mu\otimes\mu$ almost every $x,y\in X$, there exists a unique
 geodesic joining $x$ and $y$. 
\end{enumerate}
\end{thm}

\noindent{\bf Comments.}



\noindent{\bf 1.} Introducing the extra parameter $\beta$ and dividing
it by a parameter $a\to \infty$ is referred to in Physics as
``approaching the critical point''. One possible interpretation is to
see it as a device to compute various Laplace transforms of
``continuum'' limits for discrete models (see also Remark {\bf
  3.}). It is also reminiscent of the basic method of singularity
analysis consisting in approaching a singular point, and this is
exactly what we do in Sect.\ \ref{sec:appl-asympt-enum} as a special
case. On a very basic level, one should just note that letting
$a\to\infty$ has the effect of making the expectation under
$\QQ_g^{(\beta/a)}$ of quantities such as $V_\bq$ go to infinity, so
that we are intuitively taking a limit as the size of the map gets
large. Alternatively, we could also have considered conditioned laws
of the form $\QQ_g(\cdot\,|\, \{K^{-1}\leq a^{-1}V_\bq\leq K\})$ for
some fixed $K>1$, and let $a\to\infty$.

\medskip

\noindent{\bf 2.}  Considering uniform measures on bipartite
quadrangulations with $n$ faces amounts to conditioning the measure
$\QQ_g^{(\beta)}$ on $\{|F(\bq)|=n\}$, regardless of $\beta$, and a
similar result as Theorems \ref{sec:mains-results-1} and
\ref{sec:mains-results-2} could thus be considered as a conditioned
version of these statements. We believe that such generalizations can
be obtained at the price of technical complications. Note that in the
case $g=0$, a result similar to Theorem \ref{sec:mains-results-1}
appears in Le Gall \cite{legall06} in such a conditioned setting,
though the convergence does not really take $\mu_\bq$ into account
(see however \cite[Remark a., p.646]{legall06}).

%

\medskip

\noindent{\bf 3.} As mentioned in the beginning of the article, a
satisfactory statement would be to get rid of the extraction in
Theorem \ref{sec:mains-results-2}, by showing that there is actually a
unique limiting point to the set of measures (\ref{eq:10}). A further
guess would be that this unique measure could be written in the form
$\mathscr{S}_g(e^{-\beta\mathcal{V}}\d(\X,{\cal
  V}))/\mathscr{S}_g(e^{-\beta\mathcal{V}})$ for some $\sigma$-finite
measure $\mathscr{S}_g$ on $\mw\times (0,\infty)$.

\medskip

\noindent{\bf 4.} Last, a natural question is to ask whether our
results are universal. That is, is it still true for, say,
Boltzmann-Gibbs distributions on maps of genus $g$ with faces of
arbitrary valences and weight $w_k$ on a face of degree $k$? The
natural generalization of $\QQ_g$ would be the measure assigning
weight $\prod_{f\in F(\bm)} w_{\deg_\bm f}$ on the rooted map $\bm$,
for a particular ``critical'' tuning of the weights $w$. In this vein,
the results of \cite{jfmgm05,mierinv,weill06,mierweill} should be
generalizable to our context.

\bigskip

\noindent{\bf Acknowledgments. }This work was initiated while visiting
the Centre for Mathematical Studies in Cambridge (UK), supported by a
London Mathematical Society funding, where the question about
uniqueness of geodesics arose in motivating discussions with James
Norris. Thanks also to Guillaume Chapuy and Gilles Schaeffer for
interesting discussions around \cite{ChMaSc}. We finally thank the two
anonymous referees and the Editorial Board for very useful comments,
which were greatly appreciated.

\section[Delayed Voronoi tessellations]{Delayed Voronoi
  tessellations of bipartite
  quadrangulations}\label{sec:coding-delay-voron}

\subsection{Fatgraphs}\label{sec:fat-graphs}

As mentioned before, maps are of combinatorial nature, and there
exists a purely combinatorial description of maps which will be useful
to our purpose. A nice introduction to the different definitions of
maps can be found in \cite[Chapter 1]{LaZv04}.

A map $\bm$ can be encoded in a triple of permutations\footnote{Note
  that in \cite{LaZv04}, the permutations act to the right, while we
  prefer to make them act to the left, which results in differences
  like the formula $\varphi\alpha\sigma=1$ below instead of
  $\sigma\alpha\varphi=1$, referenced in Remark 1.3.19 in this
  book.}. Start from a particular representation of the map $\bm$ on
the surface $S$. Around each vertex, the outgoing half-edges are
cyclically ordered in a non-ambiguous way, by considering the first
intersection points of the half-edges with a small circle around the
vertex, oriented counterclockwise.  These cyclic orders associated
with distinct vertices involve distinct elements of $E$, and thus
determine a permutation $\sigma_\bm$ of $E$, whose cycles are
naturally associated with the vertex set $V$. We also let
$\alpha_\bm(e)=\ov{e}$, and
$\varphi_\bm=\sigma_\bm^{-1}\alpha_\bm^{-1}$. The triple
$(\sigma_\bm,\alpha_\bm,\varphi_\bm)$ is a {\em fatgraph} structure on
$E$.  We will usually drop the subscripts $\bm$ when the situation is
unambiguous.

\begin{figure}
\begin{center}
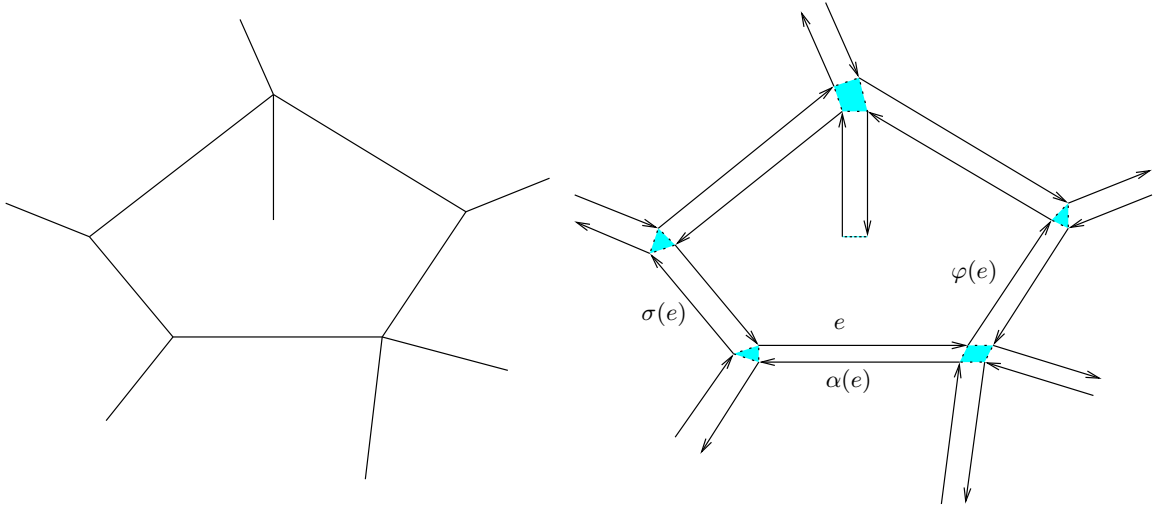
\caption{A portion of a map, emphasizing a face of degree $7$, and an
  illustration of the fatgraph structure. The action of the
  permutations $\sigma,\alpha,\varphi$ are shown on a particular
  half-edge $e$.}
\label{fig:fatgraph}
\end{center}
\end{figure}

To understand what $\varphi_\bm$ is, the reader might check on
examples that it corresponds to the following. Visiting the boundary
of each face of the graph in counterclockwise order, determines
unambiguously a cyclic ordering of the half-edges belonging to this
oriented boundary. The cycles of $\varphi_\bm$ are exactly the
so-defined cyclic orders, which are in a one-to-one natural
correspondence with faces of the graph. See Figure \ref{fig:fatgraph}
for an illustration. One deduces that the group generated by
$(\sigma_\bm,\alpha_\bm,\varphi_\bm)$ acts transitively on $E$,
amounting to the fact that maps are connected.

Conversely, every triple $(\sigma,\alpha,\varphi)$ of permutations on
a finite set $E$ such that $\varphi\alpha\sigma=1$, $\alpha$ is an
involution without fixed points, and the permutation group generated
by $(\sigma,\alpha,\varphi)$ acts transitively on $E$, can be
represented in the above way as the triple associated to a uniquely
determined map of some genus $g$, determined by Euler's formula.

A fancy intuitive image of a fatgraph is to consider it as a road
network with cars driving on the left lane. An edge and its reversal
constitute the two lanes of the road, while the vertices are
roundabouts connecting the roads. The cycles of $\sigma$ correspond to
the outgoing roads from the roundabouts, the cycles of $\varphi$
correspond to the roads taken by a driver always taking first left at
each roundabout, and the cycles of $\alpha$ correspond to a U-turn.

In a fatgraph, each half-edge $e$ belongs to exactly one cycle of
$\sigma$ and $\varphi$, and thus determines unambiguously a vertex and
face of the map, called the vertex and face that are incident to $e$:
in terms of graphs, they correspond to the origin vertex of $e$ and
the face lying to the left of $e$, respectively. By slight abuse of
notation, for $f$ a face of a map, we will write $e\in f$ if the
half-edge $e$ is incident to $f$, i.e.\ belongs to the cycle
corresponding to $f$ in the fatgraph structure.

If $e_1,e_2$ are two distinct half-edges incident to the same face, we
can break the cycle of $\varphi$ containing $e_1,e_2$ into two
linearly ordered sequences of edges which are denoted by $[e_1,e_2]$
and $[e_2,e_1]$. Formally:
$$[e_1,e_2]=\left\{\varphi^i(e_1):i\geq
0,e_2\notin\{\varphi^j(e_1),0\leq j\leq i-1\}\right\}\, ,$$ so that
$[e,e]=\{e\}$ with this definition. We also define
$(e_1,e_2]=[e_1,e_2]\setminus\{e_1\}$, and so on.

Another notion that will be useful is that of {\em corner}. In a map,
the corner incident to the half-edge $e$ is the intersection of a
small neighborhood around the vertex incident to $e$, and the angular
region comprised between $e$ and $\sigma(e)$ (to be perfectly
accurate, we should define the corner as the germ of such sets). In
the sequel, we will often speak indistinctly of a half-edge or its
incident corner.

\subsection{Quadrangulations with sources and delays}\label{sec:bijection}

We first consider unrooted maps, and will re-inject the root
afterwards.  For $g\geq 0,k\geq 1$, we let $\bQ^\circ_{g,k}$ be the
set of $3$-tuples of the form $(\bq,\bx,D)$ where
\begin{description}
\item[(Q1)] $\bq$ is a bipartite (unrooted)
  quadrangulation 
of genus $g$,
\item[(Q2)] $\bx=(x_1,\ldots,x_k)$ are elements of $V(\bq)$,
\item[(Q3)] $D=[d_1,\ldots,d_k]$ is an element of the set ${\cal
  D}(\bq,\bx)$ of {\em delay vectors}, i.e.\ elements of $\Z^k$
  defined modulo addition of a common integer, such that 
\begin{equation}\label{eq:11}
|d_i-d_j|<d_{\bq}(x_i,x_j)\qquad \mbox{ for every }\qquad 1\leq i\neq
j\leq k\, ,
\end{equation} and such that
\begin{equation}\label{eq:12}
d_\bq(x_i,x_j)+d_i-d_j\in 2\N\, ,\qquad 1\leq i, j\leq k\, .
\end{equation}
\end{description}
When $k=1$ this definition is void, and $D=[0]$ will usually be
omitted from the notation.  Note that a delay vector such as in (Q3)
exists if and only if $\min_{i\neq j}d_\bq(x_i,x_j)\geq 2$. Indeed, in
that case we can take a bicoloration of $V(\bq)$ in black and white,
and notice that two vertices are at even distance if and only if they
are of the same color. Setting $d_i=0$ if $x_i$ is white, and $d_i=1$
if $x_i$ is black, yields a vector $(d_1,\ldots,d_k)$ with the wanted
property. Conversely, if $d(x_i,x_j)=1$ for some $i\neq j$, then
obviously we should have $d_i\neq d_j$ and $|d_i-d_j|<1$ which is
impossible. This implies some restrictions on the quadrangulations
involved in $\bQ^\circ_{g,k}$ --- for instance, they should have at
least $k$ distinct vertices.

With any element of $\bQ^\circ_{g,k}$, we can associate a partition of
the set of edges of $\bq$ that is very much reminiscent of Voronoi
tessellations in a Euclidean space. We fix arbitrarily a
representative of a delay vector $(d_1,\ldots,d_k)$, and define
$\bl_i(x)=d_\bq(x,x_i)+d_i$ for $x\in V(\bq)$.  Then for $1\leq i\leq
k$, and for $x$ and $y$ adjacent in $\bq$,
\begin{equation}\label{eq:5}
|\bl_i(x)-\bl_i(y)|=1\, ,
\end{equation}
since this quantity can be at most one and vertices are at even
distance of $x_i$ if and only if they have the same color in a
bipartite coloration of the vertices. 
Thus, defining
$$\bl(x)=\min_i(d_i+d_{\bq}(x,x_i))\, ,$$ we deduce that
$|\bl(x)-\bl(y)|\leq 1$ for every adjacent $x,y$. 

On the other hand, by
(Q3) we have $d(x_1,x_i)+d_1-d_i\in 2\Z$, which shows that
$\bl_i(x_1)-\bl_j(x_1)\in 2\Z$ for every $i,j$. By (\ref{eq:5}), this
property propagates from $x_1$ to all other vertices, and since the
graph is connected it holds that $\bl_i-\bl_j$ takes its values in $2\Z$
for every $i,j$. Consequently, we cannot have $\bl_i(x)=\bl_j(y)$ with
adjacent $x,y$ for any $i,j$, and this implies that $\bl(x)\neq \bl(y)$
for adjacent $x,y$. Therefore
$$|\bl(x)-\bl(y)|=1\, ,\qquad \mbox{ for all adjacent }x,y\, .$$ 

This allows us to fix an orientation $E_{1/2}\subset E(\bq)$ by
enforcing $\bl(e^+)=\bl(e^-)-1$ for all $e\in E_{1/2}$. The function
$\bl$ decreases along chains of edges of $E_{1/2}$, and plainly, the
only vertices without edges pointing outwards are $x_1,\ldots,x_k$.
Thus, any chain of edges of $E_{1/2}$ can be extended into a maximal
chain ending at some element of $\{x_1,\ldots,x_k\}$.  Also, any such
chain with initial vertex $x$ is geodesic, and the endpoint of a
maximal extension must be a $x_i$ with $i$ such that
$\bl(x)=\bl_i(x)$. To see this, recall that each of the functions
$\bl_i$ makes a step in $\{-1,1\}$ when passing from a vertex to one
of its neighbors, by (\ref{eq:5}). Moreover, we just saw that the
function $\bl$ decreases by $1$ when going along one step in the
chain. So if the chain terminates at $x_i$, it must hold that
$\bl(y)=\bl_i(y)$ for every vertex $y$ along the chain. Since
$\bl(x)=d_\bq(x_i,x)+d_i$ and $\bl(x_i)=d_i$, the chain has length
$d_\bq(x_i,x)$ and is a geodesic.

We color each edge of $E_{1/2}$ with one in $k$ colors in the
following way. For $e\in E_{1/2}$, consider the maximal path starting
with $e$ and turning as much to the left as possible, that is, after
arriving at a vertex, it takes the first outgoing edge it encounters
in clockwise order around the vertex. If this path ends at vertex
$x_i$, assign the color $i$ to the edge $e$. 

\begin{defn}
The set $E_i(\bq,\bx,D)$ of edges colored $i$ is called the
$D$-delayed Voronoi tile with sources $\bx$ and center $x_i$.
\end{defn}

The $D$-delayed Voronoi tiles form a partition of $E_{1/2}$. The
reason for the name is that when $D=0$, the set $E_i(\bq,\bx,0)$
contains all the edges pointing from a vertex that is strictly closer
to $x_i$ than any other source, which is the usual property of Voronoi
tiles in Euclidean space. A general delay vector $D$ corresponds to
the fact that distances from different sources are measured with a
relative advantage of $d_i-d_j$ of source $x_j$ on source $x_i$.  Note
that it might be that no edge pointing from a vertex $x$ for which
$d_\bq(x,x_i)=\min_j d_\bq(x,x_j)$ lies in $E_i(\bq,\bx,0)$. In fact,
it holds that if
$$V_i(\bq,\bx,D)=\{x_i\}\cup\{e^-:e\in E_i(\bq,\bx,D)\}\, ,\qquad
1\leq i\leq k\, ,$$ then these sets cover $V(\bq)$ and we have
$$\left\{x\in V(\bq): \bl_i(x)=\bl(x)<\min_{j\neq i}\bl_j(x)\right\}
\subsetneq V_i(\bq,\bx,D)\subseteq \{x\in V(\bq):\bl_i(x)=\bl(x)\}\,
.$$ It is not difficult to see that the first inclusion is always
strict as soon as $k\geq 2$, because one can find some $x\in
V_i(\bq,\bx,D)$ and some $j\neq i$ such that $x\in V_j(\bq,\bx,D)$,
and therefore it holds that $\bl_i(x)=\bl_j(x)=\bl(x)$. The second
inclusion can be strict, although this is harder to see at first: on
Figure \ref{fig:exemple}, the only vertex $x$ with label $3$ is at
distance $3$ of $x_1$ and $2$ of $x_2$, and $d_1=0,d_2=1$, so that
$\bl_1(x)=\bl_2(x)=\bl(x)$, however, $x$ does not belong to
$V_1(\bq,\bx,D)$, because the (only) leftmost geodesic started from
$x$ ends at $x_2$.

\subsection{The bijection}\label{sec:bijection-1}

An embedded graph is said to be labeled if its vertices are assigned
integer values $\bl=(\bl(x),x\in V)$, in such a way that
$|\bl(x)-\bl(y)|\leq 1$ for every adjacent $x,y$. Labeled maps are
defined accordingly.  We let $\BLM^\circ_{g,k}$ be the set of pairs of
the form $(\bm, [\bl])$ such that
\begin{description}
\item[(LM1)] $\bm$ is a map of genus $g$ with $k$ faces indexed as
  $f_1(\bm),\ldots,f_k(\bm)$
\item[(LM2)] $\bl=(\bl(x),x\in V(\bm))$ is a labeling of $\bm$, and
  $[\bl]$ is its class up to an additive constant in $\Z$.
\end{description}

Our goal in the remaining of this section is to construct and study a
bijective mapping $\Psi_{g,k}^\circ$ between $\bQ^\circ_{g,k}$ and
$\BLM^\circ_{g,k}$, and its inverse mapping $\Phi_{g,k}^\circ$. We
first make some preliminary remarks.

\medskip

\noindent{\bf 1.} The construction starts from a particular embedding
of an element of $\bQ_{g,k}^{\circ}$, and builds an embedding of the
image labeled map, by deletion of some of the edges and vertices of
the initial quadrangulation, and addition of new edges.

In particular, if $(\bq,\bx,D)$ and $(\bm,[\bl])$ correspond by the
construction, the set of vertices of $\bm$ will be identified with the
set $V(\bq)\setminus\{x_1,\ldots,x_k\}$ of vertices of $\bq$ distinct
from the sources. This will be important in the sequel. This is also
the reason why we will give the same name $\bl$ to labeling functions
on a quadrangulation $\bq$ and its image map $\bm$, as the labeling
function on $\bm$ is indeed the restriction to the vertices of $\bm$
of the labeling function $\bl$ on $\bq$ introduced above.

\medskip

\noindent{\bf 2.}  The mapping $\Psi^\circ_{g,k}$ is a natural analog
of a bijection of Schaeffer between pointed planar maps and
well-labeled trees \cite{CSise,BdFGmobiles}. More precisely, in the
case $(g,k)=(0,1)$, the set $\BLM^\circ_{0,1}$ is constituted of
labeled (unrooted) planar trees, which themselves are in one-to-one
correspondence with {\em well-labeled} trees, for which all labels are
positive and at least one vertex has label $1$, by taking any labeling
$\bl$ and adding $-\min \bl+1$. 
A generalization of Schaeffer's bijection to genus $g\geq 1$ is
discussed in Marcus \& Schaeffer \cite{MaSc01}. 
It is recovered as the
special case $k=1$ of our construction, where the set
$\BLM^\circ_{g,1}$ is constituted of labeled maps of genus $g$ with
only one face, called {\em $g$-trees} in \cite{MaSc01}.

\medskip

\noindent{\bf 3.} We can understand our generalized bijection as
running simultaneously the Marcus-Schaeffer construction at $k$
distinct competing vertices. Informally, let water flow at unit speed
from the sources $x_1,\ldots,x_k$, in such a way that the water starts
diffusing from $x_i$ at time $d_i$, and takes unit time to go through
an edge. When water currents emanating from different edges meet at a
vertex (whenever the water initially comes from the same source of
from different sources), they can go on flowing into unvisited edges
only by respecting the rules of a roundabout, i.e.\ edges that can be
attained by turning around the vertex counterclockwise and not
crossing any other current. The process ends when the water cannot
flow anymore, and the tile $E_i(\bq,\bx,D)$ is the set of reversed
half-edges that are visited by the water flowing from source
$x_i$. The parity condition on $D$ implies that the water flows
emanating from different sources never meet in the middle of an edge,
but always at vertices.

\medskip

\begin{figure}
\begin{center}
\includegraphics[scale=.8]{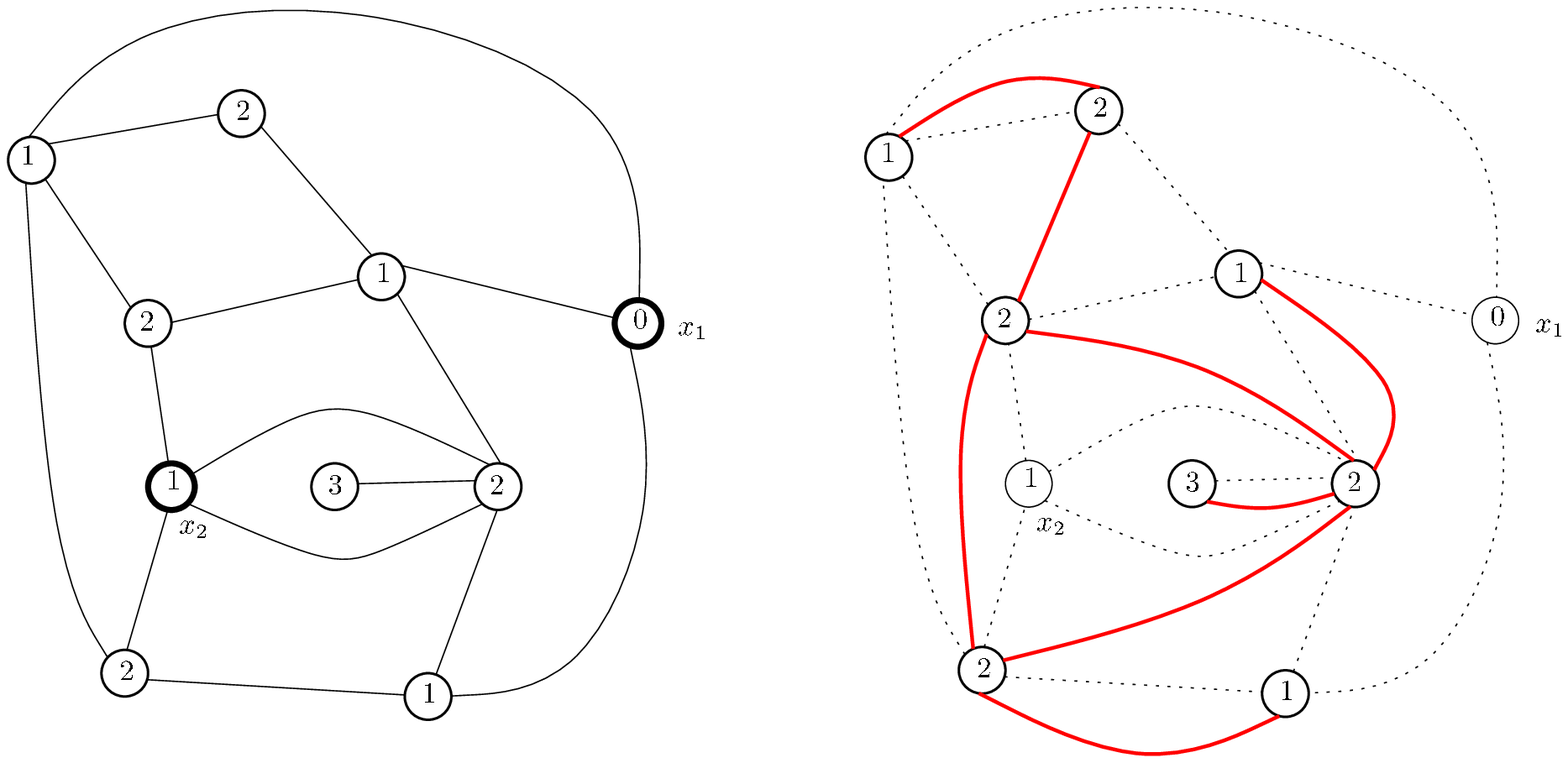}
\end{center}
\caption{An example of a $2$-pointed quadrangulation of genus $0$,
  with delays $d_1=0,d_2=1$ (left) and the associated labeled map with
  $2$ faces by $\Psi^\circ_{0,2}$}
\label{fig:exemple}
\end{figure}

We now give the rigorous construction, and refer to Figure
\ref{fig:exemple} for an example.

\subsection{Construction of the bijection}\label{sec:proof-theor-refs}

\begin{figure}
\begin{center}
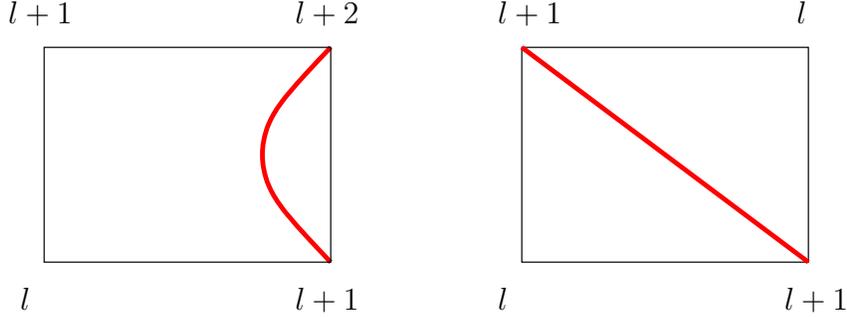
\caption{Adding new edges (thick lines) inside simple and confluent
  faces of a graphical representation of an element of
  $\bQ_{g,k}^\circ$}
\label{fig:simconf}
\end{center}
\end{figure}

We start from an element $(\bq,\bx,D)\in\bQ^\circ_{g,k}$, and take
some particular embedding of $\bq$. 
Choosing arbitrarily $(d_1,\ldots,d_k)\in D$, we label the
vertices by the function $\bl:x\in V(\bq)\mapsto
\min_i(d_i+d_\bq(x,x_i))$.  By (\ref{eq:5}), the labels of vertices
encountered when exploring a face of $\bq$, starting from the vertex
with lowest label and turning counterclockwise, must be either of the
form $l,l+1,l+2,l+1$ or $l,l+1,l,l+1$. Following \cite{MaSc01}, we
call faces {\em simple} and {\em confluent} accordingly.  We then
perform the construction of \cite{MaSc01} on this labeled object, that
is, we add extra edges to the map $\bm$, one inside each face:
\begin{itemize}
\item if the face is simple, then the added edge links the vertex with
  highest label to its successor in clockwise order around the face,
  hence splitting the face into a face with degree $4$, and a face of
  degree $2$
\item if the face is confluent, then the added edge links the two
  vertices with highest label, hence splitting the face into two
  triangular faces. 
\end{itemize}
In this way, we have defined a new map (graph) $\bq'$ whose faces have
degree $2,3$ or $4$, and with labeled vertices, whose set is still
$V(\bq)$. See Figure \ref{fig:simconf}.

The newly added edges, excluding the ones of the original map, define
a graph $\bm$ on the surface on which $\bq$ is drawn. It is easily
obtained from the construction that the vertex set of this new graph
is exactly $V(\bq)\setminus\{x_1,\ldots,x_k\}$. Indeed, the vertices
$x_1,\ldots,x_k$ have neighbors with strictly larger labels, and thus
have the least labels in any face they are incident to, and cannot be
incident to a newly added edge. Conversely, if $x$ is a vertex not in
$\{x_1,\ldots,x_k\}$, with label $\bl(x)=l$, we let $e_1,\ldots,e_r$
be the outgoing edge cycle associated with $x$ in the fatgraph
structure of $\bq$. Then there exists $i$ such that
$\bl(e_i^+)=l-1$. It is then easy that the new edge added in the face
of $\bq$ incident to $\ov{e}_i$ (i.e.\ located to the right of $e_i$)
must be incident to $x$, regardless of the face being simple or
confluent.

\begin{lmm}\label{sec:proof-theorem-refbij}
  The embedded graph $\bm$ is a map. Each face of $\bm$ contains
  exactly one element of $\{x_1,\ldots,x_k\}$, and is indexed
  accordingly as $f_1(\bm),\ldots,f_k(\bm)$.
\end{lmm}

Once this is proved, and recalling that vertices of $\bq$ are vertices
of $\bm$ as well, we use $\bl=(\bl(x),x\in V(\bq)\setminus
\{x_1,\ldots,x_k\})$, as a labeling function on $V(\bm)$.  It is now
straightforward to check that $\Psi^\circ_{g,k}(\bq,\bx,D):=(\bm,[\bl])$
satisfies properties (LM1,LM2).

\bigskip

\begin{figure}
\begin{center}
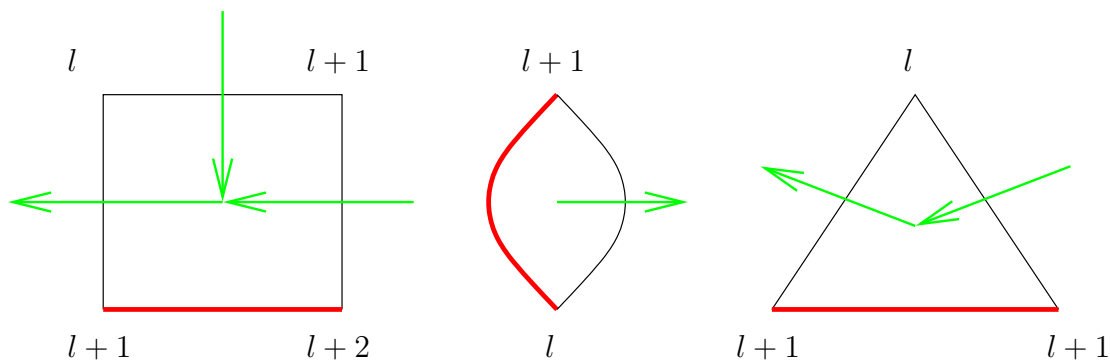
\caption{Drawing and orienting the dual edges not crossing the newly
  added edges inside each of the possible types of faces of $\bq'$}
\label{fig:schema}
\end{center}
\end{figure}

\noindent{\bf Proof of Lemma \ref{sec:proof-theorem-refbij}. } We
proceed in a similar way as in \cite{MaSc01}. Namely, we consider the
dual edges of $\bq'$ that do not cross the newly added edges, and give
them an orientation according to Figure \ref{fig:schema}, i.e.\ in
such a way that the vertex located to the right of the oriented dual
edge has strictly smaller label than the one located to the left. In
doing so, we create an oriented spanning subgraph of the dual graph of
$\bq'$, call $E_{1/2}^\circ$ the set of its edges.

By Figure \ref{fig:schema}, for every face $f$ of $\bq'$, there is a
single element of $E^\circ_{1/2}$ that goes out of it. Therefore,
there is a unique, infinite chain of $E^\circ_{1/2}$, say
$(e^\circ_m,m\geq 1)$, starting from $f$. This infinite chain
eventually cycles. By inspection of Figure \ref{fig:schema}, if $l_i$
is the label of the vertex of $\bq$ located to the right of
$e^\circ_m$, then it must hold that $l_m\geq l_{m+1}$ for every $m\geq
1$. But since $(e^\circ_m,m\geq 1)$ eventually cycles, the sequence
$(l_m,m\geq 1)$ is eventually constant, equal to $l$ say. Another look
at the first and third situations of Figure \ref{fig:schema} shows
that the only possibility is that the cycle (i.e.\ the sequence
$e^\circ_m,m\geq N$ for $n$ large enough) surrounds a single vertex of
$\bq$. This vertex must be an element of the set of sources
$\{x_1,\ldots,x_k\}$, since these are the only vertices of $\bq'$ that
are incident only to edges of $\bq$, and hence that can be surrounded
by a cycle of $E^\circ_{1/2}$.

This implies that the (oriented) graph $E^\circ_{1/2}$ has $k$
distinct connected components, and hence $\bm$ has $k$ faces. Now
starting from the union of faces and edges of $\bq'$ that are incident
to $x_i$, which is a topological disk because $x_i$ is incident
exactly once to each face, we can add the faces of $\bq'$ leading to
$x_i$ in the dual path of edges of $E^\circ_{1/2}$ without ever making
a cycle, so that the faces of $\bm$ are simply connected, hence $\bm$
is a map. \cq

\subsection{Converse construction}\label{sec:converse-contruction}

In the sequel, the label $\bl(e)$ of a half-edge (or the associated
corner) $e$ will be a notational shorthand for the label of its origin
$\bl(e^-)$.

The idea of the converse construction is simply to perform the
Marcus-Schaeffer inverse construction separately inside each face of
$\bm$.  Specifically, given $(\bm,[\bl])$, we add an extra vertex
$x_i$ inside the $i$-th face of $\bm$. We choose an arbitrary labeling
$\bl\in[\bl]$ and extend it to $V(\bm)\cup\{x_1,\ldots,x_k\}$ by
$$\bl(x_i)=\min_{e\in f_i(\bm):e^-=x} \bl(e)-1\, .$$ We let
$(\sigma,\alpha,\varphi)$ be the permutation $3$-tuple associated with
the fatgraph structure on $\bm$.

We define the successor $s(e)$ of an edge $e\in f_i(\bm)$ by
$s(e)=x_i$ if $\bl(e)=\bl(x_i)+1$, and otherwise,
$$s(e)=\varphi^{m(e)}(e)\, ,$$
where 
$$m(e)=\min\{m\geq 1:\bl(\varphi^m(e))=\bl(e)-1\}\, .$$ We then draw
an edge between the corners incident to each edges $e$ and $s(e)$ for
each $e\in E(\bm)$, in such a way that these edges do not cross each
other. This is possible, because by the constraint on the labels, we
have $|\bl(\varphi^m(e))-\bl(\varphi^{m+1}(e))|\leq 1$, which entails
that if $e$ has successor $s(e)\notin\{x_1,\ldots,x_k\}$ and $e'\in
[e,s(e)]$, then we have $\bl(e')\geq \bl(e)$ and thus $s(e')$ must be
in $(e',s(e)]$. If $s(e)=x_i$ then $e$ has minimal label inside the
  face it is incident to, and the labeling constraint prevents $e$
  from being in any set of the form $[e',s(e')]$.

When the new edges have been added in this proper non-crossing way, we
delete the edges of $\bm$ and keep only the newly added ones. 

\begin{lmm}\label{sec:converse-contruction-1}
The following properties hold.
\begin{enumerate}
\item
The graph thus obtained is a map $\bq$, which is a bipartite
  quadrangulation. 
\item
The labels $(\bl(x),x\in V(\bq))$ are equal to
$\bl(x)=\min_{1\leq i\leq k}(d_\bq(x_i,x)+d_i)$, where $d_i=\bl(x_i)$.
\item
Setting $D=[d_1,\ldots,d_k]$, the triple
$(\bq,\bx,D)$ is an element of $\bQ^\circ_{g,k}$. 
\item
The set of oriented edges of $\bq$ from $e$ to
$s(e)$, for $e\in f_i(\bm)$, is $E_i(\bq,\bx,D)$.
\end{enumerate}
\end{lmm}

\noindent{\bf Proof. }  We first check property 1. Let $e$ be an edge
of $\bm$ with label $\bl(e)=l$. Its successor $s(e)$ has label
$l-1$. If $\varphi(e)$ has also label $l$, then its successor
$s(\varphi(e))=s(e)$, so that the edges between $e,\varphi(e)$ and
their common successor form a triangle with the base edge $e$. This is
the same in the face incident to the reversal $\ov{e}$, so that the
deletion of $\{e,\ov{e}\}$ after adding the new edges yields a
$4$-valent face.

If $\varphi(e)\neq l$, then either $\varphi(e)$ has label $l+1$, or
$\varphi(e)=s(e)$ by definition of the successor. In the first case,
notice that we also have that $s(s(\varphi(e)))=s(e)$, because
$s(\varphi(e))$ must have label $l$ and lie in $[e,s(e)]$.  Thus, the
newly added edges between (the corners associated with) $e$ and
$s(e)$, $\varphi(e)$ and $s(\varphi(e))$, and $s(\varphi(e))$ and
$s(e)$, must form a $4$-valent face with the edge $e$. Now, in the
face $\ov{e}$ is incident to, observe that the successor of $\ov{e}$
must be $\sigma^{-1}(e)=\varphi(\ov{e})$, as they are successive
vertices with labels $l+1$ and $l$. Thus, the newly added edges around
the edges $e,\ov{e}$ form a $4$-valent face. It remains to treat the
case where $s(e)=\varphi(e)$, but then we are in the situation just
discussed in the face incident to $\ov{e}$.

Therefore, any (non-oriented) edge of $\bm$ is surrounded by a
$4$-valent face in the new graph, and these are the only possible
faces, as any new edge is drawn from some corner of $\bm$. Thus, the
new graph is a map on the same surface as $\bm$. It has only square
faces, and coloring vertices according to the parity of their labels
yields a bipartite coloration.

To check property 2., let $x$ be a vertex of $\bm$ incident to the
face $f_i(\bm)$, with label $l$, and let $e\in f_i(\bm)$ be an edge
pointing out from $x$. Then the (half)-edges drawn between the corners
$e,s(e),s^2(e),\ldots,x_i$ form a chain with length $\bl(x)-\bl(x_i)$,
and by definition, it turns to the left as much as possible among the
half-edges of $\bq$ along which $\bl$ decreases. Moreover, any path
using only the edges in $\bq$ must be at least as long, because these
edges always link vertices with labels that differ by $1$
exactly. Hence, this path is a geodesic and one deduces
$d_\bq(x_i,v)+d_i=l$. For the same reason, we must have
$d_\bq(x_j,v)+d_j\geq l$ for every $j\in\{1,\ldots,k\}$, so that
$l=\min_j(d_\bq(v,x_j)+d_j)$. This gives, in passing, that the edge of
$\bq$ from $e$ to $s(e)$ is an element of $E_i(\bq,\bx,D)$, yielding
the property 4.

We now check that $D=[d_1,\ldots,d_k]$ is a delay vector, which will
entail property 3. First, we have seen that for any vertex $x$, the
length of a shortest path from $x_i$ to $x_j$ passing through $x$ must
have length at least $2\bl(x)-d_i-d_j$. Hence $d_\bq(x_i,x_j)$ is at
least $2l_*-d_i-d_j$, where $l_*$ is the maximal label on a geodesic
path from $x_i$ to $x_j$. Thus $d_\bq(x_i,x_j)$ is strictly larger
than $|d_i-d_j|$ if $l_*>d_i\vee d_j$. Obviously this is the case
since the $x_i$'s are pairwise distinct and have neighboring vertices
with strictly larger labels, by construction.

Finally, we infer that $d_\bq(x_i,x_j)+d_i-d_j\in 2\N$ from the fact
that $x_i$ can be colored according to the parity of $d_i-d_j$. If
$x_i$ and $x_j$ have same color, then they are at even distance and
$d_i-d_j$ is even, and otherwise both $d_\bq(x_i,x_j)$ and $d_i-d_j$
are odd.  \cq

\bigskip

Note that, due to property 4.\ of Lemma
\ref{sec:converse-contruction-1}, the set of half-edges of $\bm$ is
identified with the set $E_{1/2}$ of edges of $\bq$ with the special
orientation of Sect.\ \ref{sec:bijection}, in such a way that
half-edges of $\bm$ incident to $f_i(\bm)$ correspond to elements of
$E_i(\bq,\bx,D)$. This will be useful in the sequel.

We now state the main result of this section.

\begin{thm}\label{BIJ}
  Let $k\geq 1$ and $g\geq 0$. The mappings $\Psi^\circ_{g,k}$ and
  $\Phi^\circ_{g,k}$ are inverse of each other, and therefore are
  bijections between the sets $\bQ^\circ_{g,k}$ and
  $\BLM^\circ_{g,k}$. Moreover, if
  $(\bm,[\bl])=\Psi^\circ_{g,k}(\bq,\bx ,D)$,
\begin{enumerate}
\item $|F(\bq)|=|\bE(\bq)|/2=|\bE(\bm)|$, more precisely
  $|E_i(\bq,\bx,D)|=\deg_\bm (f_i(\bm)),1\leq i\leq k$.
\item $|V(\bq)|=|V(\bm)|+k$, 
\item For every $1\leq i\leq k$ and $d\geq 1$,
$$\left|\left\{e\in E_i(\bq,\bx,D):d_\bq(e^-,x_i)=d\right\}\right|
  =\left|\left\{e\in f_i(\bm):\bl(e)-\!\!\! \min_{e' \in
    f_i(\bm)}\!\!\! \bl(e')+1=d\right\}\right|\, $$ (note that the
  last appearing quantity is independent of the choice of
  $\bl\in[\bl]$).
\end{enumerate}
\end{thm}

Properties 1,2,3 are easy consequences of the construction and Lemma
\ref{sec:converse-contruction-1}. They are far from being the only
important properties of $\Psi^\circ_{g,k}$. Other consequences will be
derived in Section \ref{sec:bounds-distances}.

\medskip

\proof It is easy to see why
$\Psi^\circ_{g,k}(\Phi^\circ_{g,k}(\bm,[\bl]))=(\bm,[\bl])$ from the
construction of the $4$-valent faces of $\Phi^\circ_{g,k}(\bm,[\bl])$
depicted in the previous proof. We just notice that in each of the
three cases that are detailed ($e\in E(\bm)$ is such that
$\bl(e^-)=\bl(e^+)$, $\bl(e^-)=\bl(e^+)+1$ or $\bl(e^-)=\bl(e^+)-1$),
the edge $e$ that is deleted in the end is precisely re-added by the
construction of $\Psi^\circ_{g,k}$ (the first case is that of a
confluent face, and the two others correspond to simple faces). The
labels remain unchanged in this step thanks to the last sentence of
Lemma \ref{sec:converse-contruction-1}.

We now prove that
$\Phi^\circ_{g,k}(\Psi^\circ_{g,k}(\bq,\bx,D))=(\bq,\bx,D)$.  To see
this we must check that the edges of $\bq$ can all be obtained as
edges between a corner and its successor in
$\Psi^\circ_{g,k}(\bq,\bx,D)=(\bm,[\bl])$.  Necessarily, edges must
join corners with labels differing by $1$, and cannot cross the edges
of $\bm$, so that they all lie inside some of its faces. Moreover,
inside the face $f_i(\bm)$, all corners with label $d_i+1$ must be
re-linked to the vertex $x_i$. So let $e,e'\in f_i(\bm)$ be linked by
an edge of $\bq$, for which $\bl(e)=\bl(e')+1=l+1$, but such that
$s(e)\neq e'$. We choose moreover $e,e'$ in such a way that $l$ is
minimal, that is, all the edges between a corner with label $\leq l$
and its successor are edges of $\bq$. In this case, the paths
$s(e),s^2(e),\ldots,x_i$ and $e',s(e'),\ldots,x_i$ are paths in $\bq$,
and cannot cross the edge between $e$ and $e'$. On the other hand,
$s(e),s^2(e),\ldots$ must lie in $[e,e']$, while
$s(e'),s^2(e'),\ldots$ must all lie in $[e',e]$. In particular, adding
the edge between $e$ and $e'$ to $f_i(\bm)$ breaks the latter in two
connected regions, both of which have to contain $x_i$, and this is
impossible.  \cq

\subsection{Rooting}\label{sec:rooting}

We let $\bQ_{g,k}$ be the set of $3$-tuples $(\bq,\bx,D)$, where $\bq$
is an element of $\bQ_g$, and if $\bq^\circ$ denotes the associated
unrooted map, $(\bq^\circ,\bx,D)\in \bQ^\circ_{g,k}$. Similarly, we
let $\BLM_{g,k}$ be the set of pairs $(\bm,[\bl])$ where
$(\bm^\circ,[\bl])\in \BLM^\circ_{g,k}$.

Let $(\bq,\bx,D)\in \bQ_{g,k}$, and
$(\bm^\circ,[\bl])=\Psi^\circ_{g,k}(\bq^\circ,\bx,D)$. With root edge
$e_*$ of $\bq$, corresponds naturally one of the corners of
$(\bm,[\bl])$, i.e.\ it is the unique corner such that the edge drawn
between $e$ and $s(e)$ is either equal to $e_*$ or its reversal
$\ov{e}_*$ when performing the converse construction. We root
$\bm^\circ$ at this corner, obtaining a rooted map $\bm$, and let
$\Psi_{g,k}(\bq,\bx,D)=(\bm,[\bl])$. Notice that since any edge of
$\bq$ joins two vertices with labels that differ by $1$, the two
possible orientations of a base edge of $\bq$ yield two distinct
rooted and $k$-pointed quadrangulations $(\bq,\bx)$. The following
corollary of Theorem \ref{BIJ} follows.

\begin{crl}\label{sec:rooting-1}
  The mapping $\Psi_{g,k}:\bQ_{g,k}\to \BLM_{g,k}$ is two-to-one. 
\end{crl}

Besides getting rid of symmetries, the nice feature of rooting is that
it will allow us to choose a particular labeling function inside its
class, as we will see in Sect.\ \ref{sec:struct-label-maps}. 

\subsection{Bounds on distances}\label{sec:bounds-distances}

The construction of $\Phi^\circ_{g,k}$ allows us to get useful bounds
on distances between vertices in the quadrangulation, from the label
function on the associated map.

\begin{lmm}\label{sec:bounds-distances-1}
Let $(\bq,\bx,D)\in \bQ_{g,k}$ and
$(\bm,[\bl])=\Psi_{g,k}(\bq,\bx,D)$. Take two half-edges $e_1,e_2\in
E_i(\bq,\bx,D)$ for some $i$. Still denoting by $e_1,e_2$ the
half-edges of $f_i(\bm)$ associated with $e_1,e_2$ by $\Psi_{g,k}$, we
have
$$|\bl(e_1)-\bl(e_2)|\leq d_\bq(e_1^-,e_2^-)\leq
\bl(e_1)+\bl(e_2)-2\min_{e\in [e_1,e_2]}\bl(e) +2$$ (here the interval
$[e_1,e_2]$ is the one in $\bm$).
\end{lmm}

Note that we could have interchanged the roles of $e_1$ and $e_2$, so
that we have in fact
$$d_\bq(e_1^-,e_2^-)\leq \bl(e_1)+\bl(e_2)-2\max\left(\min_{e\in
  [e_1,e_2]}\bl(e),\min_{e\in [e_2,e_1]}\bl(e) \right)+2$$

\proof The lower bound is obvious from the fact that for $e\in
E_i(\bq,\bx,D)$, $\bl(e)-\bl(x_i)=d_\bq(e^-,x_i)$ and by the triangular
inequality.  For the upper bound, consider the last half-edge $e$,
among the sequence $s^m(e_1),m\geq 0$, that lies in
$[e_1,e_2]$. Necessarily, $\bl(e)\leq \bl(e_2)$ as otherwise, the
successor of $e$ would be in $[e,e_2]$. From this, it follows that
$s(e)=s^{\bl(e_2)-\bl(e)+1}(e_2)$. Therefore, there exists a chain in
$\bq$ from $e_1^-$ to $e_2^-$ with length
$(\bl(e_1)-\bl(e))+1+(\bl(e_2)-\bl(e)+1)$, which yields the result. \cq

\section{Decomposition of labeled
  maps}\label{sec:struct-label-maps}

Our next task is to make one more step in our bijective study, by
pulling apart elementary pieces that compose an element of
$\BLM_{g,k}$. Our presentation resembles that of \cite{MR1802530}, the
two differences being that our maps are rooted, which allows us to get
rid of the symmetry factors, and, more importantly, the vertices are
labeled.  In the following, it will often be implicit that
$(g,k)\notin\{ (0,1),(0,2)\}$, which need special treatment.

\subsection{Map reductions}\label{sec:map-reductions}

We let $\bM_{g,k}$ be the set of rooted maps of genus $g$ with $k$
indexed faces. Elements of $\bM_{g,k}$ for fixed $k$ but with a large
number of vertices, show a dendritic pattern branching on a graph with
a simple structure.  In order to study this in a precise way, we are
going to perform several reductions on the maps.  We let
$\bM_{g,k}^{\geq 2}$ and $\bM_{g,k}^{\geq 3}$ be the subsets of
$\bM_{g,k}$ consisting of the maps in which all vertices have degree
at least $2$, resp.\ at least $3$. It is easy to see that
$\bM_{g,k}^{\geq 2}$ is an infinite set (take a simple loop with an
arbitrary number of edges). On the other hand, $\bM_{g,k}^{\geq 3}$ is
finite according to the combination of Euler's formula
\begin{equation}\label{eq:17}
|V|-|\bE|+k=2-2g
\end{equation}
and the fact that
\begin{equation}\label{eq:4}
2|\bE|=\sum_{v\in V}\deg v\geq 3|V|\, ,
\end{equation}
which entails $|\bE|\leq 3k+6g-6$. Note that equality holds if and
only if every vertex has degree $3$ exactly. We let $\bM_{g,k}^3$ be
the subset of $\bM_{g,k}^{\geq 3}$ for which this extra constraint is
satisfied. According to these definitions, the set $\bM^{\geq
  3}_{0,2}$ should be empty. We make an exception by letting
$\bM^3_{0,2}=\bM^{\geq 3}_{0,2}$ consist in the two {\em loop maps}
made of a rooted loop with an origin, bounding two faces with degree
$1$ (there are two such maps according to whether the face indexed $1$
stands to the left or to the right of the root edge). Elements of
$\bM^{\geq 3}_{g,k}$ will usually be denoted by the letter $\mg$, and
their root by $\eg_*$.

There are natural projection maps
$$\bM_{g,k}\longrightarrow \bM_{g,k}^{\geq 2}\quad ,\qquad
\bM_{g,k}\longrightarrow \bM_{g,k}^{\geq 3}\, ,$$ which roughly
consist in taking the largest subgraph belonging to $\bM^{\geq
  2}_{g,k}$ (resp.\ $\bM^{\geq 3}_{g,k}$), for some appropriate
definition of subgraphs, where merging edges are allowed.  More
precisely, starting from (a graphical realization of) $\bm\in
\bM_{g,k}$, we delete one by one the edges incident to vertices with
degree $1$ until no such vertex remains. This operation does not
disconnect the graph, nor does it modify its genus or the number of
its faces. The outcome is thus a map $\bm^{\geq 2}$ with $k$ faces and
without vertices of degree $1$, and which does not depend on the order
of the edges that were removed\footnote{Alternatively, elements of
  $E(\bm^{\geq 2})$ are those that can start an infinite
  non-backtracking chain $e_1,e_2,\ldots$, i.e.\ such that $e_i\neq
  \ov{e}_{i+1}$ for $i\geq 1$.}.  If $\bm$ has root $e_*$, we
naturally root $\bm^{\geq 2}$ at the first element of
$\varphi_\bm^i(e),i\geq 0$ that belongs to $E(\bm^{\geq 2})$. This
yields an element $\bm^{\geq 2}\in\bM^{\geq 2}_{g,k}$.

\begin{figure}[!h]
\begin{center}
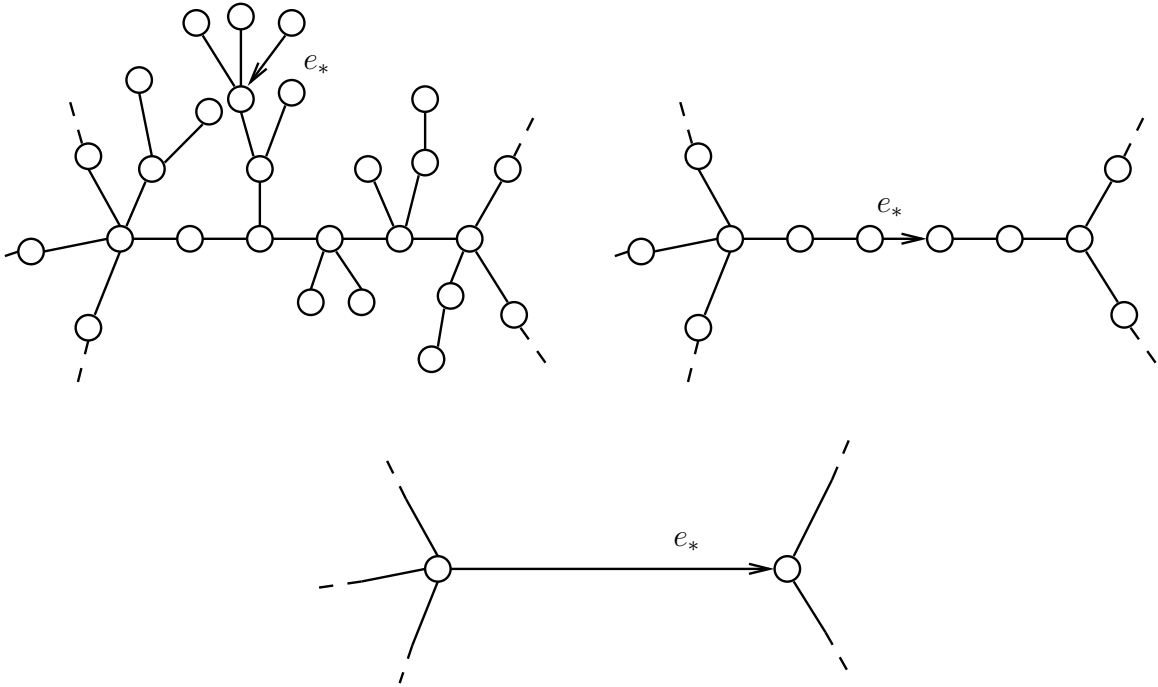
\caption{Top-left: a portion of a map $\bm$, emphasizing the tree
  components branching on $\bm^{\geq 2}$ depicted on top-right.
  Bottom: the projection $\bm^{\geq 3}$. The picture shows how the
  various projections are re-rooted depending on the location of the
  root of $\bm$.}  
\label{fig:proj}
\end{center}
\end{figure}

Likewise, starting from $\bm\in\bM_{g,k}$, we can remove all the
vertices of degree $2$ in $\bm^{\geq 2}$: we concatenate the
half-edges $e_1,\ldots,e_r$ of a simple chain linking two vertices of
$\bm^{\geq 2}$ of degrees $\geq 3$, and traversing only degree-$2$
vertices, into a single edge $\eg=e_1\ldots e_r$. We root the
resulting graph at the one edge $\eg_*=e_1\ldots e_r$ in which the
root of $\bm^{\geq 2}$ is involved. This does not modify the genus nor
the number of faces of the map, and these merging operations commute,
and thus yield a well-defined element $\bm^{\geq 3}$ of $\bM^{\geq
  3}_{g,k}$.  It is also obvious by definition that $(\bm^{\geq
  2})^{\geq 3}=\bm^{\geq 3}$. There are natural inclusions
$V(\bm^{\geq 3})\subseteq V(\bm^{\geq 2})\subseteq V(\bm)$ and
$E(\bm^{\geq 2})\subseteq E(\bm)$.  The map $\bm^{\geq 3}$ is
sometimes called the {\em homotopy type} of $\bm$. See Figure
\ref{fig:proj}.

\subsection{Some notations for continuous paths}\label{sec:some-notations}

Certain continuous functions will be crucial to encode maps, so let us
introduce some notations.  A {\em path} is a continuous function
$(w(s),0\leq s\leq \tau(w))$ taking its values in $\R$, such that
$w(0)=0$. The number $\tau(w)\geq 0$ is called its duration, and
$\wh{w}=w(\tau(w))$ its terminal value. We also let
$\ov{w}=w(\tau(w)-\cdot)-\wh{w}$ be the reversed path (beware that
this does not mean the same as the reversal of an edge in a graph). We
let $\W$ be the set of paths, and make it a Polish space by endowing
it with the distance
$${\rm dist}(w,w')=\left(\sup_{s\geq 0}|w(s\wedge
  \tau(w))-w'(s\wedge \tau(w'))|\right)\vee |\tau(w)-\tau(w')|\, .$$
If $w,w'\in \W$, we let $ww'$ be the concatenation of $w$ and $w'$
defined by
$$ww'(s)=\left\{\begin{array}{cl} w(s) & \mbox{if }0\leq s\leq
\tau(w)\\ \wh{w}+w'(s-\tau(w)) & \mbox{if }\tau(w)\leq s\leq
\tau(w)+\tau(w')=\tau(ww')
\end{array}\right.$$ 

For $(g,k)\notin\{(0,1),(0,2)\}$ introduce the sets
$$\bC_{g,k}=\left\{\left(\mg,\lgo,t^*\right):\mg\in \bM_{g,k}^{\geq
  3}, \lgo=(\lgo^\eg,\eg\in E(\mg))\in (\W^3)^{E(\mg)}, 0\leq t^*\leq
\tau^{\eg_*}\right\}$$ endowed with the natural product topology.  We
make the exceptions $\bC_{0,1}=\{(c,z)\in \W^2\}$ and
$$\bC_{0,2}=\{(\mg,\lgo,t^*):\mg\in \bM_{0,2}^{\geq
  3},\lgo=((c,z),\lgo^{\eg_*},\lgo^{\ov{\eg}_*}) \in \W^2\times
\W^3\times \W^3,0\leq t^*\leq\tau(c)\}\, .$$
For $(g,k)\notin\{(0,1),(0,2)\}$ the generic notations are
\begin{equation}\label{eq:14}
\lgo^\eg=(w^\eg,c^\eg,z^\eg)\, ,\qquad r_\eg=\tau(w^\eg)\, ,\qquad
\tau^\eg=\tau(c^\eg)\, ,\qquad
\tau_\lgo=\sum_{\eg\in E(\mg)}\tau^\eg\, .
\end{equation} 
For $g=0$ and $k=1, 2$ we let $\tau=\tau(c)$, for $k=2$ we also let
$\tau_\lgo=\tau+\tau^{\eg_*}+\tau^{\ov{\eg}_*}$.  Our next task is to
bijectively encode elements $\bm$ of $\BLM_{g,k}$ with certain
elements of $\bC_{g,k}$, whose first component is the homotopy type of
$\bm$.

\subsection{Level sets of $\bM_{g,k}\to\bM^{\geq
    3}_{g,k}$}\label{sec:map-reductions-1}

Let us see how to get back from the projection $\mg=\bm^{\geq 3}$ to
the initial map. First of all, every $\eg\in E(\mg)$ must be split
into a certain number of edges $e_1\ldots e_{r_\eg}$, with $r_\eg\geq
1$, to yield the element $\bm^{\geq 2}$, and obviously
$r_\eg=r_{\ov{\eg}}:=r_\be$ (where $\be=\{\eg,\ov{\eg}\}$) so what is
needed to do this are positive integers $(r_\be,\be\in \bE(\mg))$.

Then, consider elements $e\in E(\bm)$ such that the first element of
$\varphi_\bm^i(e),i\geq 0$ belonging to $E(\bm^{\geq 2})$ is one of
the edges $e_1,\ldots,e_{r_\be}$. These edges form an interval, and we
let $\eg(0)$ be the first element in this interval, and label the
others as $\eg(i)=\varphi_\bm^i(\eg(0)),0\leq i\leq \tau^\eg-1$.  With
this labeling, we see that every half-edge of $\bm$ is of the form
$\eg(i)$ for a unique $\eg\in \mg$. Moreover, the half edges of the
intervals $(e_{j-1},e_j),1\leq j\leq r_\be$ must form an acyclic
graph, i.e.\ a tree, possibly reduced to a single root vertex (in
particular, their set has even cardinality and is stable under
edge-reversal). Here, we take the convention that $(e_0,e_1)$ are the
elements of the list $\eg(0),\eg(1),\ldots$ that precede $e_1$
strictly. Otherwise said, we can interpret
$\eg(0),\ldots,\eg(\tau^\eg-1)$ as a forest of $r_\be$ trees, with
roots bound together by the floor half-edges $e_1,\ldots,e_{r_\be}$.
If $\eg_*$ is the root edge of $\bm^{\geq 3}$, then it deserves
special attention. Indeed, in this case, one of the edges of the
interval $[\eg_*(0),\eg_*(\tau^{\eg_*}-1)]$ must be distinguished as
the initial root edge of $\bm$, which amounts to distinguishing an
integer $t^*\in [0,\tau^\eg)$.

The contour function associated with this forest is defined by
$$c^\eg(i)=d_\bm(\eg(i)^-,V(\bm^{\geq 2})) - \un{c}^\eg(i)\, , \qquad
0\leq i\leq \tau^\eg\, ,$$ where
$$\un{c}^\eg(i)=|\{j:e_j\in [\eg(0),\eg(i)]\}|=-\inf_{0\leq j\leq
  i}c^\eg(j)\, ,\qquad 0\leq i\leq \tau^\eg-1\, ,$$ and
$\un{c}^\eg(\tau^\eg)=r_\be$, is the number of trees in the forest
already explored before time $i$.  We extend $c^\eg,\un{c}^\eg$ into
continuous functions (still denoted by $c^\eg,\un{c}^\eg$) by linear
interpolation between values taken at integers. The result is a
concatenation of $r_\be$ {\em Dyck paths}\footnote{Here we will call
  Dyck path a piecewise linear element $c\in \W$ starting at $0$, with
  slope $\pm1$ between integer times, and which attains the value $-1$
  for the first time at $\tau(c)$.} encoding the trees of the forest,
and achieving a new negative infimum when a floor edge of the forest
is traversed (note that the terminal value $-r_\be$ of $c^\eg$ is
attained for the first time at the final time $\tau^\eg$).  It is an
easy and well-known fact that the contour function entirely determines
the forest. In this way we have shown

\begin{lmm}\label{sec:level-sets-bm_g}
Assume $(g,k)\notin\{(0,1),(0,2)\}$. An element $\bm\in\bM_{g,k}$ is
unambiguously described by its projection $\mg=\bm^{\geq 3}$, a family
$(c^\eg,\eg\in E(\mg))\in \W^{E(\mg)}$ such that $c^\eg$ is the
concatenation of $r_\be$ Dyck paths, and an integer $t^*\in
[0,\tau^{\eg_*})$.
\end{lmm}

\begin{figure}
\begin{center}
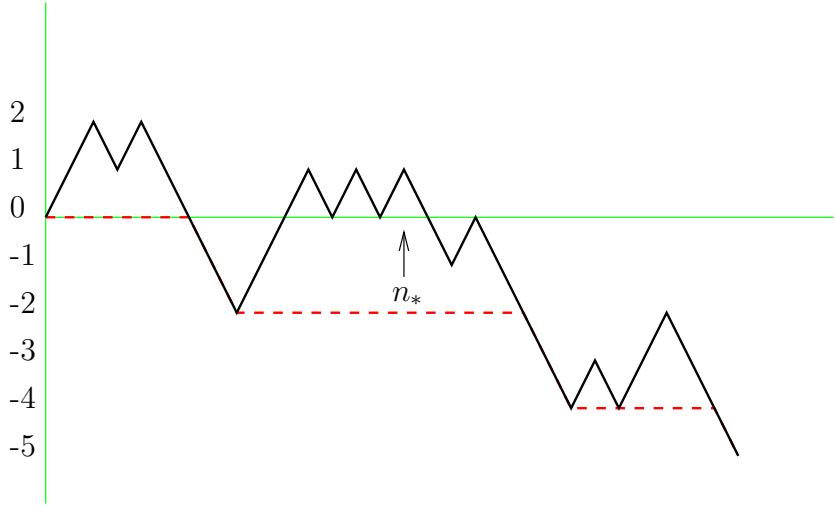
\caption{The contour process $c^{\eg_*}$ associated with the edge
  $\eg_*$ in the reconstruction of the map $\bm$ from $\bm^{\geq 3}$
  as depicted in Figure \ref{fig:proj}. The process $-\un{c}^{\eg_*}$
  is drawn in thick dashed lines.  The location of the root of $\bm$
  is distinguished by the time $t^*$.}
\label{fig:contour}
\end{center}
\end{figure}

See Figure \ref{fig:contour} for an illustration. The cases
$(g,k)=(0,1)$ or $(0,2)$ are special and requires a little
discussion. An element of $\bM_{0,1}$ is simply a rooted planar tree
with at least one edge, and is described by a single Dyck path $c$
with duration $\tau(c)>1$.  When $(g,k)=(0,2)$, the construction of
$\bm^{\geq 3}$ from an element $\bm\in \bM_{0,2}$ shows that the root
vertex of $\bm^{\geq 3}$ can be any of the vertices of the loop
constituting $\bm^{\geq 2}$. This difficulty does not appear for other
values of $(g,k)$, since the root edge of $\bm^{\geq 3}$ has an
unambiguous origin in the original map $\bm$.  While this does not
affect the final outcome considered as a map, there are several
non-equivalent ways to obtain the same result. To lift the ambiguity,
we use the extra convention that the root vertex in $\bm^{\geq 3}$ is
$e^-$, where $e$ is the first element of $E(\bm^{\geq 2})$ in the
sequence $(\varphi_\bm^i(e_*),i\geq 0)$.  Otherwise said, we single
out the tree component that contains the root and branching on the
unique loop of $\bm$, by saying that its root vertex must be the root
vertex of $\bm^{\geq 3}$. The other tree components are obtained by
taking concatenations of Dyck paths so that an element $\bm\in
\bM_{0,2}$ is unambiguously described by the element $\bm^{\geq 3}$, a
Dyck path $c$, an integer time in $[0,\tau(c))$ and two paths
  $c^{\eg_*},c^{\ov{\eg}_*}$, the first being a concatenation of
  $r_{\eg_*}-1$ Dyck paths, and the second the concatenation of
  $r_{\eg_*}$ Dyck paths.

\subsection{Labeling the maps}\label{sec:labeling-maps}

In order to complete the description of elements of $\BLM_{g,k}$, we
need to incorporate the labeling function $[\bl]$. We can use the
origin of the root edge $\eg_*$ of $\mg=\bm^{\geq 3}$ as a reference
for $[\bl]$, i.e.\ we can choose the element $\bl\in[\bl]$ such that
$\bl(\eg_*^-)=0$. This convention will be in force from this section
onwards. It is then natural to introduce a labeling function $l^\eg$
for every $\eg\in E(\mg)$, defined by
\begin{equation}\label{eq:7}
l^\eg_i=\bl(\eg(i))-\bl(\eg(0))\, ,\qquad0\leq i< \tau^\eg\, ,
\end{equation}
and $l^\eg_{\tau^\eg}=\bl(\eg(\tau^{\eg}-1)^+)-\bl(\eg(0))$, where the
notation $\eg(0),\eg(1),\ldots$ is as above the list of edges of the
trees branching on $\eg$. We interpolate linearly between integer
times, in order to define a continuous piecewise linear function with
slope in $\{-1,0,1\}$, starting at $0$ and with same duration as
$c^\eg$.  Another process that will be of crucial importance is the
process
\begin{equation}\label{eq:20}
l_i=\bl(\varphi_\bm^i(\eg_*(0)))\, \qquad 0\leq i\leq \deg_{\bm}(f_*)-1\, ,
\end{equation}
where $f_*$ is the face of $\bm$ incident to $e_*$, and
$l_{\deg_\bm(f_*)}=0$.  As usual, it is turned into an element of $\W$
by linear interpolation. In the case where $\bm$ has only one face, we
see that $l$ provides a natural exploration of the labels of all
vertices of $\bm$ in the fatgraph face order. It is easily checked
that $l$ is the concatenation
\begin{equation}\label{eq:21}
l=l^{\eg_*}l^{\varphi_\mg(\eg_*)}\ldots
l^{\varphi_\mg^{\deg_\mg(\fgo_*)-1}(\eg_*)}\, ,
\end{equation}
where $\fgo_*$ is the face of $\mg$ incident to $\eg_*$.  In order to
study scaling limits, we will still need to decompose $l^\eg$ into two
more elementary pieces, which is the purpose of the two next sections.

\subsubsection{Walk networks}\label{sec:walk-networks}

A {\em walk} is an element $w\in \W$ taking integer values on $\Z_+$,
interpolating linearly between these values, and satisfying
$|w(i)-w(i-1)|\leq 1$ for $1\leq i\leq r:=\tau(w)$ (where it is
assumed that $r\geq 1$). Let $\bW$ be the set of walks. A walk from
$0$ to $x\in \Z$ is a walk with final value $\wh{w}=x$, their set is
written $\bW_{0\to x}$.

Let $\mg\in\bM^{\geq 3}_{g,k}$ for some $g,k\geq 0$. We let
$\ell(V)$ and $\ell(E)$ be the spaces of functions $V(\mg)\to\R$
and $E(\mg)\to\R$ respectively.  We also let $\ell_0(V)$ be those
functions $\bl$ of $\ell(V)$ with $\bl(\eg_*^-)=0$. A function
$\theta\in\ell(E)$ is {\em antisymmetric} if
$\theta(\eg)=-\theta(\ov{\eg})$ for every $\eg$, and we let
$\ell_-(E)$ be the space of antisymmetric functions. For $\bl\in
\ell(V)$, we let $\nabla \bl\in \ell_-(E)$ be defined by $\nabla
\bl(\eg)=\bl(\eg^+)-\bl(\eg^-)$. As $\mg$ is connected, the kernel of
the linear map $\nabla$ is the space of constant functions, so the
restriction of $\nabla$ to $\ell_0(V)$ is an isomorphism of vector
spaces.  A function $\theta$ in the image $\nabla_V$ of $\nabla$ is
called a {\em potential}.

A {\em walk network} on $\mg$ is a collection $(w^\eg,\eg\in
E(\bm))\in \bW^{E(\bm)}$ such that
$$w^{\ov{\eg}}(t)=w^\eg(r_\be-t)-\wh{w}^\eg\, ,\qquad 0\leq t\leq
r_\be=\tau(w^\eg)=\tau(w^{\ov{\eg}})\, ,$$ and such that the
antisymmetric function $\theta=\wh{w}=(\wh{w}^\eg,\eg\in E(\mg))$ is a
potential\footnote{Alternatively, walk networks on $\mg$ are in
  correspondence with labeled maps $\bm\in \bM_{g,k}^{\geq 2}$ such
  that $\bm^{\geq 3}=\mg$.}.  We let $\BWN_{\mg}$ be the set of walk
networks on $\mg$. If we are given an antisymmetric function
$\theta\in \ell_-(E)$ assuming integers values, we let
$\BWN_\mg(\theta)$ be the set of walk networks $(w^\eg,\eg\in
E(\mg))\in \BWN_{\mg}$ such that $\wh{w}=\theta$.

If $(\bm,[\bl])\in \BLM_{g,k}$, and using the notations of Section
\ref{sec:map-reductions-1}, we set, for every $\eg\in E(\mg)$ of the
form $\eg=e_1\ldots e_{r_\be}$,
$$w^\eg(i)=\bl(e_{i+1}^-)-\bl(e_1^-)\, , \qquad 0\leq i\leq
r_\be-1\, ,$$ and $w^\eg(r_\be)=\bl(e_{r_\be}^+)-\bl(e_1^-)$. This
defines a walk with duration $r_\be$. Moreover, the family
$(w^\eg,\eg\in E(\mg))$ is a walk network. The reversal property
$w^{\ov{\eg}}(\cdot)=w^\eg(r_\be-\cdot)-\wh{w}^\eg$ is clear since
$\ov{\eg}=\ov{e}_{r_\be}\ldots\ov{e}_1$, and by definition the
function $(\wh{w}^\eg,\eg\in E(\mg))$ is $\nabla\bl|_{V(\mg)}$.

\subsubsection{Discrete snakes}\label{sec:discrete-snakes}

Finally, we say that the pair $(c,z)\in \W^2$ is a discrete snake on
$r\geq 1$ tree components if
\begin{itemize}
\item $c$ is a concatenation of $r$ Dyck paths,
\item $z$ is taking integer values on $\Z_+$, interpolating linearly
  between these values, with integer steps in $\{-1,0,1\}$, satisfies
  $\tau(z)=\tau(c)$, and $z(i)=z(j)$ whenever $c(i)=c(j)=\min_{i\wedge
    j\leq k\leq i\vee j}c(k)$ (which means that $i$ and $j$ are times
  where the same vertex of the forest encoded by $c$ is explored),
\item moreover, $z(i)=0$ whenever $i$ is an integer time at which $c$
  achieves a new minimum (this means that $i$ is a time at which a new
  tree component is explored, and the label of the root is fixed to be
  $0$). 
\end{itemize}
In words, a discrete snake is the process of heights and labels when
exploring an integer-labeled forest in contour order, where the labels
differ by at most $1$ in absolute value between adjacent vertices, and
the labels of the forest's roots is $0$.  We let $\bS_r$ be the set of
discrete snakes on $r$ tree components.

If $(\bm,[\bl])\in \BLM_{g,k}$ we set, whenever $0\leq i\leq \tau^\eg$
is such that $\eg(i)\in (e_{m-1},e_m]$ (recall the notations of
Sect.\ \ref{sec:map-reductions-1}),
$$z^\eg(i)=\bl(\eg(i))-\bl(e_m)\, ,$$ and we interpolate linearly
between integer times.  The pair $(c^\eg,z^\eg)$ is then the discrete
snake associated with the forest branching on the edge $\eg$, where
the labels in each tree component have been shifted so that the roots
have labels $0$. With the data $(w^\eg,c^\eg,z^\eg)$ defined so far,
the label function $l^\eg$ introduced at the beginning of the
subsection is recovered as
\begin{equation}\label{eq:13}
l^\eg(i)=z^\eg(i)+w^\eg(\un{c}^\eg(i))\, .
\end{equation}
We also use (\ref{eq:13}) to define the function $l^\eg$ out of any
$(\mg,\lgo,t^*)\in \bC_{g,k}$, where $\un{c}^\eg(t)=-\inf_{0\leq s\leq
  t}c^\eg(s)$ for $0\leq s\leq \tau^\eg$, and define the function $l$
as in (\ref{eq:21}).

Summing up the study of this section, we obtain

\begin{prp}\label{sec:labeling-maps-1}
{\rm (i)} Let $(g,k)\notin\{(0,1),(0,2)\}$. Then $\BLM_{g,k}$ is in
bijection with the subset $\bC_{g,k}^{\rm map}\subset \bC_{g,k}$ of
elements $(\mg,\lgo,t^*)$ such that 
  \begin{enumerate}
\item $(w^\eg,\eg\in E(\mg))\in \BWN_\mg$,
\item $(c^\eg,z^\eg)\in \bS_{r_\be}$ for every $\eg\in E(\mg)$, and
\item $t^*$ is n integer in $[0,\tau(c^{\eg_*}))$.
\end{enumerate}

\noindent{\rm (ii)} $\BLM_{0,1}$ is in bijection with the set
$\bC_{0,1}^{\rm map}=\{(c,z)\in \bS_1:\tau(c)>1\}$. 

\noindent{\rm (iii)} $\BLM_{0,2}$ is in bijection with the subset
$\bC_{0,2}^{\rm map}\subset\bC_{0,2}$ of elements $(\mg,\lgo,t^*)$
with
\begin{enumerate}
\item $w^{\eg_*}\in \bW_{0\to 0}$ and $w^{\ov{\eg}_*}=\ov{w}^{\eg_*}$,
\item $(c,z)\in \bS_1$ and $t^*\in[0,\tau(c))$ is an integer,
\item $(c^{\eg_*},z^{\eg_*})\in \bS_{r_{\eg_*}-1}$ and
  $(c^{\ov{\eg}_*},z^{\ov{\eg}_*})\in \bS_{r_{\eg_*}}$.
\end{enumerate}
\end{prp}

\section{Discrete and continuous Boltzmann
  measures}\label{sec:discr-cont-path}

We now explore natural measures and probability distributions on
tessellated genus-$g$ maps, and their continuous counterparts.

We define the measure $\QQ_{g,k}$ on the set $\bQ_{g,k}$ as the
$\sigma$-finite measure assigning mass $12^{-|F(\bq)|}$ to every
element $(\bq,\bx,D)$. Otherwise said, if we recall the notations of
Section \ref{sec:mains-results}, one has
\begin{equation}\label{eq:19}
\QQ_{g,k}(\d(\bq,\bx,D))=\QQ_g(\d\bq)\,
(V_\bq\mu_\bq)^{\otimes k}(\d \bx)\, \#_{{\cal D}(\bq,\bx)}(\d D)\, ,
\end{equation}
 where $\#_{{\cal D}(\bq,\bx)}(\d D)$ is the counting measure on the
 (possibly empty) set of delays in $\bq$ between sources $\bx$.  Up to
 Sect.\ \ref{sec:exceptional-cases} we assume
 $(g,k)\notin\{(0,1),(0,2)\}$.

According to Theorem \ref{BIJ} and Corollary \ref{sec:rooting-1}, the
push-forward of $\QQ_{g,k}$ by the mapping $\Psi_{g,k}$ is the measure
on the set $\BLM_{g,k}$ that assigns mass $2\cdot 12^{-|\bE(\bm)|}$ to
any element $(\bm,[\bl])$.  Using the considerations of Section
\ref{sec:struct-label-maps}, we see that this measure can in turn be
obtained out of measures on the more elementary pieces that compose a
labeled maps (walk networks and snakes). More precisely, consider the
bijective mapping $\Xi_{g,k}:\BLM_{g,k}\to \bC_{g,k}^{\rm map}$ that
is described in Proposition \ref{sec:labeling-maps-1}. 
Define
$$\wt{\Psi}_{g,k}:=\Xi_{g,k}\circ\Psi_{g,k}\quad :\quad\bQ_{g,k}\to
\bC_{g,k}^{\rm map}\, .$$ Then the image measure $\LM_{g,k}$ of
$\QQ_{g,k}$ by $\wt{\Psi}_{g,k}$ assigns weight $2\cdot
12^{-\tau_\lgo/2}$ to $(\mg,\lgo,t^*)\in \bC_{g,k}^{\rm map}$, as one
checks that $\tau_\lgo/2=|\bE(\bm)|=|F(\bq)|=V_\bq-\chi(g)$ whenever
$(\bq,\bx,D),(\bm,[\bl]),(\mg,\lgo,t^*)$ are associated by our
bijections.  We rewrite this measure in the form
\begin{equation}\label{eq:18}
\LM_{g,k}=2\sum_{\mg\in \bM_{g,k}^{\geq 3}}\LM_{\mg}
\end{equation}
where
$$\LM_{\mg}=\delta_{\{\mg\}}\WN_{\mg}(\d
(w^\eg)_{\eg\in E(\mg)})\, \bigotimes_{\eg\in
  E(\mg)}\DSN_{\tau(w^\eg)}(\d (c^\eg,z^\eg)) \,
\#_{\tau(c^{\eg_*})}(\d t^*)
$$
and
\begin{itemize}
\item $\WN_{\mg}$ is the measure assigning weight $3^{-\sum_{\be\in
    \bE(\mg)}r_\be}$ to a walk network $(w^\eg)\in \BWN_{\mg}$,
\item $\DSN_r$ is the measure assigning weight
  $2^{-\tau(c)}3^{-(\tau(c)-r)/2}$ to an element $(c,z)\in \bS_r$,
\item lastly, $\#_{\tau(c^{\eg_*})}$ is the counting measure on
  $\{0,\ldots,\tau(c^{\eg_*})-1\}$.
\end{itemize}

Individually, these measures are $\sigma$-finite measures, which are
finite except for $\WN_{\mg}$. The measure $\DSN_r$ is a probability
measure, under which the $r$ trees involved in the forest are
independent Galton-Watson trees with geometric $(1/2)$ offspring
distribution, and with uniformly distributed labels among the allowed
ones, conditionally on the forest. To see this, recall that $(c,z)\in
\bS_r$ encodes a forest with $r$ tree components, and $\tau(c)$
oriented edges (counting the $r$ floor edges between tree
components). A component of this forest, which is a labeled tree with
$m'$ edges say, contributes a probability factor $2^{-2m'-1}3^{-m'}$
to the weight $2^{-\tau(c)}3^{-(\tau(c)-r)/2}$. Also, the measure
$\DSN_r(\d(c,z))\#_{\tau(c)}(\d t^*)$ is an infinite measure.

\subsection{Walk networks}\label{sec:free-paths-measures}

For simplicity we let in this section $V=V(\mg),E=E(\mg),\bE=\bE(\mg)$
for some $\mg\in \bM_{g,k}^{\geq 3}$. We let $\wh{\WN}_{\mg}$ be the
image measure of $\WN_{\mg}$ under $(w^\eg,\eg\in E)\mapsto
((r_\be,\be\in \bE),(\wh{w}^\eg,\eg\in E))$.

\subsubsection{Gaussian laws on
  potentials}\label{sec:gauss-laws-potent}

Fix $\mg\in\bM_{g,k}$, and let $v_*=\eg_*^-$ be the root
vertex. Recall the notations of Section \ref{sec:labeling-maps}.  We
endow $\ell_0(V)$ with the scalar product $(\bx,\bx')=\sum_{v\neq
  v_*}x_vx_v'$, and endow the vector space $\nabla_V:=\nabla\ell_0(V)\subset
\ell_-(E)$ of potentials on $\mg$ with the scalar product
$$\langle\nabla \bx,\nabla\bx'\rangle=(\bx,\bx')\, .$$ This space is
isometric to $\ell_0(V)$, and of dimension $|V|-1=\dim \ell_0(V)$.
We let $\lambda_\nabla$ be the measure on $\nabla_V$ which is
the image of $\prod_{v\in V\setminus\{v_*\}}\d y_v$ by
$\nabla:\ell_0(V)\to \ell_-(E)$, and call it the Lebesgue measure
on $\nabla_V$: it assigns mass $1$ to the unit cube delimited
by the vectors of the orthonormal basis $(\nabla\ind_v,v\neq
v_*)$. For $(r_\be,\be\in \bE)$ a family of positive weights, we also
let
$${\cal E}_\br(\bx)=\sum_{\be\in
  \bE}x_\be^2/r_\be\, ,$$ which defines a positive definite quadratic
form on $\nabla_V$. The Gaussian distribution on
$\nabla_V$ associated with ${\cal E}_\br$ is the measure
$${\cal N}^\mg_\br(\d\bx)=
\frac{1}{(2\pi)^{\frac{|V|-1}{2}}}\exp\left(-\frac{{\cal
    E}_\br(\bx)}{2}\right)\sqrt{\det {\cal E}_\br}\,
\lambda_\nabla(\d\bx)\, ,$$ which is a probability distribution. One
should be careful that $\det {\cal E}_\br$ is the determinant of the
quadratic form defined on $(\nabla_V,\langle
\cdot,\cdot\rangle)$ and not the whole space $\ell_-(E)$.

\begin{lmm}\label{sec:gauss-laws-potent-1}
Let ${\rm ST}(\mg)$ denote the set of spanning trees of $\mg$. Then,
one has
$$\det {\cal E}_\br=\sum_{\mathfrak{a}\in{\rm ST}(\mg)}\prod_{\be\in
  \bE(\mathfrak{a})}(r_\be)^{-1}\, .$$ Otherwise said, this
determinant is a partition function for spanning trees on $\mg$ with
weight $1/r_\be$ for the edge $\be$.
\end{lmm}

\proof In the orthonormal basis $(\nabla\ind_{\{v\}},v\neq v_*)$ of
$(\nabla_V,\langle\cdot,\cdot\rangle)$, it is easy to check
that the matrix $(m_{uv},u,v\in V\setminus\{v_*\})$ of the quadratic
form ${\cal E}_\br$ is given by $m_{vv}=\sum_{\eg:\eg^-=v}1/r_\eg$ and
$m_{uv}=-\sum_{\eg:\eg^-=u,\eg^+=v}1/r_\eg$. Otherwise said, this is a
minor (excluding row and column corresponding to $v_*$) of the
weighted Laplacian matrix on $V$ with weight $1/r_\eg$ on edge $\eg$,
and one concludes by the matrix-tree theorem (and the result does not
depend on the root of $\mg$ but only on the map). \cq

\bigskip

The following result will also be useful when estimating integrals
with respect to $\lambda_\nabla$. For a spanning tree $\mathfrak{a}\in
{\rm ST}(\mg)$, we privilege the orientation of its edges
$E_{1/2}(\mathfrak{a})$ by letting them point away from the root
vertex $v_*$. If $\by\in \R^{E_{1/2}(\mathfrak{a})}$, we define a
function $\theta^\mathfrak{a}_\by\in \ell_-(E)$ extending $\by$ by
letting $\theta^\mathfrak{a}_\by(\ov{\eg})=-y_\eg$ for $\eg\in
E_{1/2}( \mathfrak{a})$ and then
$$\theta^\mathfrak{a}_\by(\eg)=\sum_{\eg'}y_{\eg'}\, ,$$ where the sum
is over the half-edges $\eg'\in E(\mathfrak{a})$ of the unique
oriented simple chain going from $\eg^-$ to $\eg^+$. It is easy to see
that $\theta_\by^{\mathfrak{a}}\in \nabla_V$ is a potential, and any
potential is recovered from its restriction to $E_{1/2}(\mathfrak{a})$
by applying $\theta^{\mathfrak{a}}$.

\begin{lmm}\label{sec:gauss-laws-potent-2}
Let $H:\R^{E}\to \R_+$ be measurable. Then for $\mathfrak{a}\in
{\rm ST}(\mg)$,
$$\int_{\nabla_V}\lambda_\nabla(\d\bx)H(\bx)=\int_{\R^{E_{1/2}(\mathfrak{a})}}\d
\by H(\theta^\mathfrak{a}_\by)\, .$$
\end{lmm}

\proof
By definition
$$\int_{\nabla_V}\lambda_\nabla(\d\bx)H(\bx)=\int_{\R^{V\setminus\{v_*\}}}\d\bz
H(\nabla\bz)\, .$$ The result follows from the change of variables
$\by(\eg)=\nabla\bz(\eg),\eg\in E_{1/2}(\mathfrak{a})$, provided we
can justify that its Jacobian is $1$. This results from the fact that
this change of variables is linear, and that its representative matrix
in an appropriate basis is triangular with only $1$'s on the
diagonal. To see this, perform the search-depth exploration of
$\mathfrak{a}$ starting from any half-edge $\eg_1\in E(\mathfrak{a})$
with $\eg_1^-=v_*$, and letting, for $i\geq 2$,
$\eg_i=\varphi_{\mathfrak{a}}^{m(i)}(\eg_{i-1})$ where $m(i)$ is the
first integer $j\geq 0$ such that
$(\varphi_{\mathfrak{a}}^j(\eg_{i-1}))^+$ is a vertex that does not
belong to the set $\{v_*,\eg_1^+,\ldots,\eg_{i-1}^+\}$ of already
explored vertices: this is possible as long as $i\leq |V|-1$.  This
gives a labeling $\eg_1,\ldots,\eg_{|V|-1}$ of the half-edges of
$\mathfrak{a}$ pointing out from $v_*$, such that the source $\eg_i^-$
is always a vertex belonging to the set $\{\eg_j^+:j<i\}$, while the
target $\eg_i^+$ is a newly explored vertex.  \cq

\subsubsection{Key lemmas}\label{sec:key-lemmas}

Define the quantity
$${\cal Z}^\mg(\br)=\sum_{\mathfrak{a}\in {\rm
    ST}(\mg)}\prod_{\be\notin \bE(\mathfrak{a})}2\pi r_\be\, , \qquad
\br\in \R_+^E\, ,$$ which is a partition function for the complement
of a spanning tree, with weight $r_\be$ on edge $\be$, and defines a
homogeneous function of degree $|\bE|-|V|+1$. We define a measure
$\wh{\PN}_\mg$ on $\R_+^\bE\times \nabla_V$ by
$$\wh{\PN}_\mg\left(\d\br,\d\bx\right)=\frac{\d \br}{\sqrt{{\cal
      Z}^\mg(\br)}} {\cal N}^\mg_\br(\d \bx)\, .$$

Let $\Delta_\bE=\{\bu=(u_\be,\be\in \bE)\in\R_+^\bE:\sum_{\be\in
  \bE}u_\be=1\}$ be the $|\bE|-1$-dimensional simplex indexed by
$\bE$, and $\lambda_\Delta$ denote the (scaled) Lebesgue measure on
$\Delta_\bE$, satisfying
$$\int_{\R_+^\bE}H(\br)\d \br =\int_{\R_+}\rho^{|\bE|-1}\d
\rho\int_{\Delta_\bE}\lambda_\Delta(\d \bu)H(\rho \bu)$$ for every
integrable $H$. In particular, taking $H(\br)=\exp(-\sum_{\be\in \bE}
r_\be)$, we see that its total mass equals $1/(|\bE|-1)!$.

We define the measure $\wh{\PN}_\mg^{(\rho)}$ on $\rho\Delta_\bE\times
\nabla_V$ as the image measure by $(\bu,\bx)\mapsto(\rho\bu,\bx)$ of
$$\rho^{|\bE|-1} \frac{\lambda_\Delta(\d\bu)}{\sqrt{{\cal
      Z}^\mg(\rho\bu)}}{\cal N}_{\rho\bu}^\mg(\d\bx)\, ,$$ hence
giving a disintegration of $\wh{\PN}_\mg$ with respect to
$\rho(\br):=\sum_{\be\in \bE}r_\be$.

\begin{lmm}\label{sec:key-lemmas-2}
For $\bu\in n\Delta_{\bE}$, we let $\wt{\bu}$ be the point of
$n\Delta_\bE\cap \Z^\bE$ that is closest to $\bu$ for the Euclidean
metric (taking some arbitrary convention for ties). Then for any
bounded function $H:\R_+^\bE\to\R$,
$$\left|\frac{1}{n^{|\bE|-1}}\sum_{\bv\in n\Delta_\bE\cap
  \Z^\bE}H(\bv)-\int_{\Delta_{\bE}}\lambda_\Delta(\d\bu)H(\wt{n\bu})\right|\leq
(\sup |H|)\eps_n\, ,$$ where the sequence $(\eps_n,n\geq 1)$ is
independent of $H$ and goes to $0$ as $n\to\infty$.
\end{lmm}

\proof To see this, note that 
\begin{equation}\label{eq:28}
\int_{\Delta_\bE}\lambda_\Delta(\d\bu)H(\wt{n\bu})=\sum_{\bv\in
  n\Delta_\bE\cap \Z^\bE}H(\bv)\lambda_\Delta(\{\bu\in
\Delta_\bE:\wt{n\bu}=\bv\})\, .
\end{equation} 
Now, as long as $\bv\in n\Delta_\bE\cap \Z^\bE$ has all its
coordinates positive, i.e.\ does not belong to the boundary
$n\partial\Delta_\bE$ and therefore stands at distance at least $1$
from it as it has integer coordinates, the sets $\{\bu\in
\Delta_\bE:\wt{n\bu}=\bv\}$ do not intersect $\partial \Delta_\bE$,
and have equal $\lambda_\Delta$-measure by an obvious property of
translation-invariance. There are $\binom{n-1}{|\bE|-1}$ elements in
$n\Delta_\bE\cap \N^\bE$, which is the number of integer compositions
of $n$ into $|\bE|$ positive parts. We conclude that for any $\bv\in
n\Delta_\bE\cap \N^\bE$,
$$0\leq
\lambda_\Delta(\Delta_\bE)-\binom{n-1}{|\bE|-1}\lambda_\Delta(\{\bu\in
\Delta_\bE:\wt{n\bu}=\bv\})\leq
\lambda_\Delta(V_{1/n}(\partial\Delta_\bE)) \, ,$$ where
$V_{1/n}(\partial \Delta_\bE)$ is the $1/n$-neighborhood of
$\partial\Delta_\bE$, and
$\lambda_\Delta(V_{1/n}(\partial\Delta_\bE))\to 0$ as $n\to\infty$.
Since 
$$\binom{n-1}{|\bE|-1}\sim
n^{|\bE|-1}/(|\bE|-1)=n^{|\bE|-1}\lambda_\Delta(\Delta_\bE)\, ,$$ we
conclude that $\lambda_\Delta(\{\bu\in \Delta_\bE:\wt{n\bu}=\bv\})=
(1+o(1))/n^{|\bE|-1}$ as $n\to\infty$ independently of $\bv\in
n\Delta_\bE\cap\N^{\bE}$, and since $|n\partial \Delta_\bE\cap
\Z^\bE|=o(n^{|\bE|-1})$, the conclusion follows easily from
(\ref{eq:28}).  \cq

\begin{lmm}\label{sec:key-lemmas-1}
Set $\varsigma=(8/9)^{1/4}$. Let $G_a,a>0$ be a uniformly bounded
family of functions on $\N^\bE\times \Z^E$ such that
$$G_a(\br_a,\bx_a)\to G(\br,\bx)$$ whenever $\br_a/(2a)^{1/2}\to \br$
and $\bx_a/(\varsigma a^{1/4})\to \bx\in \R^E$. Let also
$\rho_a/(2a)^{1/2}\to \rho$ with $\rho_a\in \N$. Then
$$\frac{\varsigma^{k+2g-1}}{2^{\frac{|\bE|-1}{2}}a^{\frac{|V|+|\bE|-3}{4}}}
\wh{\WN}_\mg(G_a\ind_{\{\rho(\br)=\rho_a\}})\to
\wh{\PN}^{(\rho)}_\mg(G)\, .$$
\end{lmm}

\proof Let $\bW_{0\to x}^r$ be the set of walks with duration $r$ that
end at $x$, and $p_n(N)=3^{-r}|\bW_{0\to N}^r|$ be the probability
that a uniform walk of duration $r$ ends at $N$.  Note that
$k+2g-1=|\bE|-|V|+1$ by (\ref{eq:17}), and write the left-hand side in
the statement as
\begin{eqnarray*}
\lefteqn{
  \frac{\varsigma^{|\bE|-|V|+1}}{2^{\frac{|\bE|-1}{2}}a^{\frac{|V|+|\bE|-3}{4}}}
  \sum_{\br\in \N^\bE}3^{-\sum_{\be\in \bE}r_\be}\sum_{\by\in
    \Z^{V\setminus\{v_*\}}}\sum_{(w^\eg)\in \BWN_\mg(\nabla \by)}
  G_a\left(\br,\nabla \by\right)\ind_{\{\rho(\br)=\rho_a\}}}\\ &=&
\frac{\varsigma^{|\bE|-|V|+1}}{2^{\frac{|\bE|-1}{2}}a^{\frac{|V|+|\bE|-3}{4}}}
\sum_{\br\in \rho_a\Delta_\bE\cap \N^\bE}\sum_{\by\in
  \Z^{V\setminus\{v_*\}}} G_a\left(\br,\nabla \by\right)\prod_{\eg\in
  E_{1/2}} p_{r_\eg}(\nabla\by(\eg)) \\ &\build\sim_{a\to\infty}^{}&
\rho_a^{|\bE|-1}\frac{\varsigma^{|\bE|-|V|+1}}{2^{\frac{|\bE|-1}{2}}
  a^{\frac{|V|+|\bE|-3}{4}}}
\int_{\Delta_\bE}\lambda_\Delta(\d\bu)\int_{ \R^{V\setminus\{v_*\}}}\!\!\!\!\!\d
\by\, G_a\left(\wt{\rho_a\bu},\nabla \lfloor\by\rfloor\right)
\prod_{\eg\in E_{1/2}}p_{\wt{\rho_a\bu}(\eg)}(\nabla\by(\eg)) \\ &
\build\sim_{a\to\infty}^{}& \rho^{|\bE|-1}
\int_{\Delta_\bE}\lambda_\Delta(\d \bu)\int_{
  \R^{V\setminus\{v_*\}}}\!\!\!\!\!\d \by\, G_a\left(\wt{\rho_a\bu},\nabla
\lfloor a^{1/4}\varsigma \by\rfloor\right) \prod_{\eg\in
  E_{1/2}}\varsigma a^{1/4}p_{\wt{\rho_a\bu}(\eg)}(\nabla\lfloor
a^{1/4}\varsigma\by\rfloor(\eg))\, ,
\end{eqnarray*}
where we made use of Lemma \ref{sec:key-lemmas-2} for $n=\rho_a$ at
the third step of the computation. Here, $\lfloor \bx\rfloor$ is the
vector whose components are the integer parts of the components of
$\bx$, and we write $\wt{\bu}(\eg)$ for the coordinate of $\wt{\bu}$
corresponding to the edge $\be=\{\eg,\ov{\eg}\}$.

By the local limit theorem \cite[Theorem VII.3.16]{petrov75}, there
exists\footnote{Until the end of the proof, the same notation $C$ will
  be used for different such bounding constants.} some $C\in(0,\infty)$
such that for every $n\geq 1$, $N\in \Z$ and letting $x=N/
(n^{1/2}\sqrt{2/3})$,
$$\left|\sqrt{\frac{2}{3}}
n^{1/2}p_n(N)-\frac{1}{\sqrt{2\pi}}e^{-\frac{x^2}{2}}\right|\leq
n^{-1/2}\frac{C}{1+|x|^2}\, .$$ From this, and the easily checked fact
that $n\partial\Delta_\bE$ is at distance $C_\bE>0$ (independent on
$n$) of the set $\{n\bu:\bu\in \Delta_\bE:\wt{n\bu}\notin
n\partial\Delta_\bE\}$ (so that the components $nu_\eg$ are
bounded away from $0$ for $n\bu$ in this set), it follows that
$$\varsigma a^{1/4}p_{\wt{\rho_a\bu}(\eg)}(\nabla\lfloor
a^{1/4}\varsigma \by\rfloor(\eg))\leq
C\min\left(\frac{\sqrt{\rho}}{x_\eg^2},\frac{1}{\sqrt{\rho
    u_\eg}}\right)\, ,$$ where $x_\eg=|\nabla\by(\eg)|$.  Combined
with the fact that the $G_a$ are uniformly bounded allows us to bound
the last displayed double integral by
$$C\max(\rho^{3|\bE|/2-1},\rho^{|\bE|/2-1})\int_{\nabla_V}\lambda_\nabla(\d
\bx) \int_{\Delta_\bE}\lambda_\Delta(\d\bu) \prod_{\eg\in
  E_{1/2}}\left(\frac{1}{x_\eg^2}\ind_{\{x_\eg>1\}}+ \frac{1}{
  \sqrt{u_\eg}}\ind_{\{x_\eg\leq 1\}}\right)\, .$$ Developing the
product and integrating in $\bu$, we finally obtain the bound
$$C\max(\rho^{3|\bE|/2-1},\rho^{|\bE|/2-1})
\int_{\nabla_V}\lambda_\nabla(\d\bx)\prod_{\eg\in
  E_{1/2}}\min(x_\eg^{-2},1)\, .$$ Taking a spanning tree
$\mathfrak{a}\in{\rm ST}(\mg)$ and using Lemma
\ref{sec:gauss-laws-potent-2}, this is bounded by
\begin{equation}\label{eq:26}
\max(\rho^{3|\bE|/2-1},\rho^{|\bE|/2-1})
\int_{\R^{E_{1/2}(\mathfrak{a})}}\d\bx \left(\sum_{\eg\in
  E_{1/2}(\mathfrak{a})}\min(x_\eg^{-2},1)\right)^{|\bE|}\, ,
\end{equation}
and this is finite. Therefore, we can apply dominated convergence in
the integral expression above and pass to the limit as $a\to \infty$
to get
\begin{eqnarray*}
\lefteqn{\frac{\varsigma^{|\bE|-|V|+1}}{2^{\frac{|\bE|-1}{2}}
    a^{\frac{|V|+|\bE|-3}{4}}}\wh{\WN}_\mg\left(
  G_a\ind_{\rho(\br)=\rho_a} \right)}\\ &\build\longrightarrow_{a\to
  \infty}^{}&
\rho^{|\bE|-1}\int_{\Delta_\bE}\lambda_\Delta(\d\bu)\int_{\nabla_V}\lambda_\nabla(\d\bx)
G(\rho\bu,\bx)\prod_{\eg\in E_{1/2}}\frac{e^{-x_\eg^2/(2\rho
    u_\eg)}}{\sqrt{2\pi \rho u_\eg}} \\ &=&\rho^{|\bE|-1}
\int_{\Delta_\bE}\lambda_\Delta(\d\bu)\int_{\nabla_V}\lambda_\nabla(\d\bx)
\frac{\exp\left(-\frac{{\cal E}_{\rho
      \bu}(\bx)}{2}\right)}{\prod_{\be\in \bE}\sqrt{2\pi \rho u_\be}}
G(\rho \bu,\bx)\, ,
\end{eqnarray*}
and the result follows from the definition of ${\cal
  N}^{\mg}_\br$ and ${\cal Z}^\mathfrak{\mg}(\br)$
and from Lemma \ref{sec:gauss-laws-potent-1}. \cq

\begin{lmm}\label{sec:path-network-measure}
Let $G_a,a>0$ be functions $\N^\bE\times
\Z^{E}\to \R$ such that
$$|G_a(\br,\bx)|\leq K\exp\left(-c\, a^{-1/2}\sum_{\be\in
  \bE}r_\be\right)$$ for some $K,c>0$, and satisfying the same
convergence assumptions as in the previous lemma.  Then
$$\frac{\varsigma^{k+2g-1}}{2^{\frac{|\bE|}{2}}
a^{\frac{|V|+|\bE|-1}{4}}}\wh{\WN}_\mg\left(
G_a\right) \build\longrightarrow_{a\to \infty}^{}\wh{\PN}_\mg(G)\, .$$
\end{lmm}

\proof This is an immediate consequence of the previous lemma, using
easy dominated convergence arguments along with the $\rho$-dependent
bound (\ref{eq:26}). Details are left to the reader. \cq

\subsubsection{Path network measure}\label{sec:free-path-network}

We are now able to describe the scaling limit of the measure
$\WN_\mg$.  
Let $P^r_{0\to x}(\d w)$ be the law of the standard Brownian bridge
with duration $t$ from $0$ to $x$, and $P^r_{0\to x}(\d w\d w')$ be
the image measure of $P^r_{0\to x}(\d w)$ under $w\mapsto
(w,\ov{w})\in \W^2$.

Take an orientation $E_{1/2}$ of $E$.  We define a measure by
$$\PN_\mg\left(\d\left(w^\eg,\eg\in E\right)\right)=\int
\wh{\PN}_{\mg}(\d\br,\d\bx) \prod_{\eg\in E_{1/2}}P^{r_e}_{0\to
  x_e}(\d w^\eg\d w^{\ov{\eg}})\, ,$$ which is a measure on paths
networks, i.e.\ collections $(w^\eg,\eg\in E)$ of paths with
durations $r_\eg:=\tau(w^\eg)$ such that $w^{\ov{\eg}}=\ov{w}^\eg$ and
$(\wh{w}_\eg,\eg\in E)\in \nabla_V$ is a potential. This
measure does not depend on the choice of $E_{1/2}$ because of the
well-known invariance property \cite{revyor} of $P^t_{0\to x}$ under
$w\mapsto \ov{w}$.

We define a scaling operator $\phi^a$ on $\W$ by
$$\phi^a (w)= \left( \frac{9}{8a} \right)^{1/4}w((2a)^{1/2}s)\,
,\qquad 0\leq s\leq \tau(w)/(2a)^{1/2}\, .$$ We extend the scaling
function $\phi^a$ to a function on powers of $\W$ in the natural way,
i.e.\ it acts on every component by applying the scaling $\phi^a$.

\begin{prp}\label{sec:path-network-measure-1}
Let $H$ be a continuous function $\W^{E}\to \R$ such that
$$|H((w^\eg,\eg\in E))|\leq K\exp\left(-c \sum_{\eg\in
  E}\tau(w^\eg)\right)\, ,$$ for some $K,c>0$.  Then
$$\frac{\varsigma^{k+2g-1}}{2^{\frac{|\bE|}{2}}a^{\frac{|V|+|\bE|-1}{4}}}\WN_\mg\left(
H\circ\phi^a\right) \build\longrightarrow_{a\to\infty}^{}\PN_\mg(H)\,
.$$
\end{prp}

\proof We rely on Lemma \ref{sec:path-network-measure} and the
following statement. 

\begin{lmm}\label{sec:path-network-measure-2}
For $n\in\N$ and $x\in\Z$, let $\bP^{n}_{0\to x}$ be the uniform
distribution on $\bW^n_{0\to x}$.  Let $r>0,x\in \R$ be fixed and
assume that $r_a,x_a$ are integers such that $r_a/(2a)^{1/2}\to r$ and
$(9/8a)^{1/4}x_a\to x$. Then the probability measures
$\phi^a_*\bP^{r_a}_{0\to x_a}$ converge in distribution to $P^r_{0\to
  x}$ in $\W$.
\end{lmm}

\proof This is a classical elaboration on Donsker's invariance
principle. See for instance \cite[Lemmas 4.5--6]{duq02} for the scheme
of the proof (in this reference, the author considers more general
stable processes than Brownian motion but is also concerned only in
the case where $x=0$, however the proof in the general $x$ case is the
same). \cq

\medskip

To end the proof of Proposition \ref{sec:path-network-measure-1}, 
we apply Lemma \ref{sec:path-network-measure} to the function
$$G_a(\br,\bx)= \int_{\W^{E}}\bigotimes_{\eg\in
  E_{1/2}}\bP^{r_\eg}_{0\to x_\eg}(\d w^\eg\d w^{\ov{\eg}})
H(\phi^a(w^\eg,\eg\in E))\, ,
$$
where as above the measure $\bP^{r}_{0\to x}(\d w\d
w')$ is the law of $(w,\ov{w})$ under $\bP^{r}_{0\to x}$.  
\cq

\bigskip

\rem The path network measure $\PN_\mg$ actually has the following
interpretation, which we state without proof as this is not needed in
the present article. Assume that the numbers $r_\be>0, \be\in \bE$ are
given and consider a graph with same structure as $\mg$ but where the
edge $\be$ is isometric to a real segment of length $r_\be$. This
defines a {\em metric graph} \cite[Chapter 3.2.2]{burago01}
$\mg_\br$. Let $\Theta_{\mg_\br}$ be the set of continuous functions
on $\mg_\br$ that vanish at $v_*$. With any path network
$(w^\eg,\eg\in E)$ we may associate an element of $\Theta_{\mg_\br}$
in the following way: for $x\in \mg_\br$ on the half-edge $e$ and at
distance $0\leq s\leq r_\be$ from its origin, pick a chain
$e_1,\ldots,e_m$ of half edges starting at $v_*$ and ending at
$e^-$. Then let $H(x)=\sum_{i=1}^m\wh{w}^{e_i}+w^\eg(s)$. The reader
will convince himself that this definition does not depend on the
choice of the orientation of $e$ nor on that of the chain
$e_1,\ldots,e_m$, because of the properties of path networks. Thus,
the measure
$${\cal N}^\mg_\br(\d\bx)\prod_{\eg\in E_{1/2}}P^{r_\eg}_{0\to
  x_\eg}(\d w^\eg\d w^{\ov{\eg}})$$ induces a law
$\mathsf{GFF}_{\mg_\br}(\d H)$ on $\Theta_{\mg_\br}$. This law is what
is referred to as the Gaussian free field on $\mg_\br$. We refer to
\cite{shefGFF} for a nice introduction to this object (here the metric
graphs are really $1$-dimensional objects, so the associated free
fields are well-defined random functions that are much less elaborate
than higher-dimensional free fields on which \cite{shefGFF}
focuses). To explain what $\mathsf{GFF}$ is, imagine that we run a
``Brownian motion'' on $\mg_\br$. It behaves as a linear Brownian
motion when moving on an edge of $\mg_\br$, and when it encounters a
vertex of the graph, its excursions out of the vertex choose
independently at random in which incident edge they will diffuse. The
diffusion is killed when it first encounters the vertex $v_*$. We let
$G(x,y)$ be the density\footnote{The function $y\mapsto G(x,y)$ can
  also be defined as the solution of a Poisson equation
  $\Delta_{\mg_\br}H=\delta_x$ vanishing at $v_*$, for the appropriate
  definition of the Laplacian $\Delta_{\mg_\br}$ on regular functions
  of $\Theta_{\mg_\br}$.} at vertex $y$ of the occupation measure of
the diffusion starting at $x$. Then $\mathsf{GFF}_{\mg_\br}$ is the
law of the centered Gaussian field $(H(x),x\in \mg_\br)$ with
covariance
$$\cov(H(x),H(y))=G(x,y)\, ,\qquad x,y\in \mg_\br\, .$$

\subsection{Discrete and continuous snakes}\label{sec:discr-cont-snak}

We now discuss the scaling limits of the measures $\DSN_r$ on discrete
snakes. This is more classical from a probabilistic point of view as
we are back into the realm of probability measures. We define the
continuous snake measures $\SN_r,r>0$, as the law of the pair
$((B_s)_{0\leq s\leq I^{-1}_r},(S_s)_{0\leq s\leq I^{-1}_r})$, where
\begin{itemize}
\item $(B_s,s\geq 0)$ is a standard Brownian motion in $\R$, 
\item $(I_s=\inf_{0\leq s'\leq s}B_{s'},s\geq 0)$, is the infimum
  process of $B$ and $I^{-1}$ is its right-continuous inverse,
\item conditionally on $B$, the process $S$ is centered Gaussian in
  $\R$, with covariance $$\cov(S_s,S_{s'})=\inf_{s\wedge s'\leq s''\leq
    s\vee s'}(B_{s''}-I_{s''})\, .$$
\end{itemize}
As an application of Kolmogorov's continuity criterion and the
standard fact that trajectories of a Brownian motion are a.s.\ locally
H\"older continuous, a continuous version of $S$ exists, and $\SN_t$
is a distribution on $\W^2$.

We rely on a version of invariance principles such as
\cite{CSise,MMsnake,janmarck05}, but holding for forests. It is easily
derived from Theorem 3 of \cite{miergwmulti} by using the methods of
\cite[Section 2.4]{duqleg02} to derive the results for the contour
process of the tree rather than the height process. There is a
difference between the label processes used in this reference and the
present paper (due to a different way of exploring the trees), but the
proof are adapted in a straightforward way.

We define the scaling functions $\vartheta^a,\psi^a$ on
  $\W$ by
$$\vartheta^a(c)=\frac{1}{(2a)^{1/2}} c(2as)\, ,\qquad 0\leq s\leq
\tau(c)/2a\, ,$$ and
$$\psi^a(z)=\left(\frac{9}{8a}\right)^{1/4} z(2as)\, ,\qquad 0\leq s\leq
  \tau(z)/2a\, .$$

\begin{lmm}\label{sec:discr-cont-snak-1}
  Fix $t>0$ and assume $r_a\in \N$ is such that $r_a/(2a)^{1/2}\to
  r$. Then, the measure $(\vartheta^a,\psi^a){}_*\DSN_{r_a}$ converges
  weakly to $\SN_r$ as $a\to\infty$. If $U_n$ is the uniform
  distribution over $\{0,1,\ldots,n-1\}$, we also have that the
  measure defined by
$$\int\DSN_{r_a}(\d(c,z))U_{\tau(c)}(\d
  t)H(\vartheta^a(c),\psi^a(z),t/2a)\, ,$$ for positive measurable
  $H$, converges weakly towards the measure defined by
$$\int\SN_r(\d(c,z))\int_{\R_+}\frac{\d t\ind_{\{0\leq t\leq
    \tau(c)\}}}{\tau(c)} H(c,z,t)\, .$$
\end{lmm}

\subsection{Scaling limits}\label{sec:cont-meas-label}

In this section, we assume $(g,k)\notin\{(0,1),(0,2)\}$.  We define a
measure on $\bC_{g,k}$ by
\begin{equation} 
\CLM_{g,k}= 2\sum_{\mg\in\bM^3_{g,k}}\CLM_{\mg}
\end{equation}
and $\CLM_{\mg}$ is the appropriate image measure obtained by
rearranging terms of
$$ \delta_{ \mg}\, \PN_{\mg}(\d
(w^\eg)_{\eg\in E(\mg)})\, \bigotimes_{\eg\in
  E(\mg)}\SN_{r_\eg}(\d (c^\eg,z^\eg))\, \d t^*\ind_{\{0\leq
  t^*\leq \tau^{\eg_*}\}}\, .$$ 
Note that the exponent $\geq 3$ of the
first term defining $\LM_{g,k}$ is replaced by a $3$ here.  We define
the scaling operator $\xi^a$ by
$$\xi^a\left(\mg,\lgo,t^*\right)
=\left(\mg,((\phi^a(w^\eg),\vartheta^a(c^\eg),\psi^a(z^\eg)),\eg\in
  E(\mg)),t^*/2a\right).$$

\begin{thm}\label{sec:cont-meas-label-1}
Assume $(g,k)\notin \{(0,1),(0,2)\}$. Let $H$ be a continuous function on
$\bC_{g,k}$ such that $$|H|\leq
K\tau_\lgo^{-1}\exp\left(-\beta\tau_\lgo\right)$$ for some
$K,\beta>0$. 
Then it holds that
$$\frac{\varsigma^{k+2g-1}}{2^{3g-2+\frac{3k}{2}}
  a^{\frac{5g}{2}+\frac{5k-7}{4}}}\LM_{g,k}(H\circ \xi^a)
\build\longrightarrow_{a\to\infty}^{}\CLM_{g,k}(H)\, .$$
\end{thm}


\medskip

\proof First of all, we check that for every $\br\in \N^\bE$, 
\begin{equation}\label{eq:16}
\bigotimes_{\eg\in
  E}\DSN_{r_\eg}\left(\exp\left(-\frac{\beta}{a}\tau_{\lgo}\right)\right)
\leq K'\exp\left(-c_\beta a^{-1/2}\sum_{\be\in \bE}r_\be\right)\, .
\end{equation} Indeed, under
$\DSN_r(\d(c,z))$, the terminal time $\tau$ of $c$ has same law as the
first hitting time of $-r$ by a simple random walk on the
integers. Thus $\tau_\lgo$ has same law under $\bigotimes_{e\in
  E}\DSN_{r_\eg}$ as the first hitting time of $-2\sum_{\be\in
  \bE}r_\be$ by the random walk. Standard computations show that
\begin{equation}\label{eq:27}
\DSN_1(s^{\tau(c)})=\frac{1-\sqrt{1-s^2}}{s}\, ,\qquad 0< s\leq 1\, ,
\end{equation}
so that  
$$\bigotimes_{\eg\in
  E(\mg)}\DSN_{r_\eg}\left(\exp\left(-\frac{\beta}{a}\tau_\lgo
\right)\right)=\exp\left(-2\phi\left(\frac{\beta}{a}\right)\sum_{\be\in
  \bE}r_\be\right)\, ,$$ where
$\phi(\beta)=-\log(e^\beta(1-\sqrt{1-e^{-2\beta}}))$ is the Laplace
exponent of the first hitting time of $-1$ by the walk, so
$\phi(\beta/a)\sim (2\beta/a)^{1/2}$ as $a\to \infty$, hence
(\ref{eq:16}). We can thus write $\LM_{\mg}(H)$ as
$$\int \WN_{\mg}(\d(w^\eg)) \int\bigotimes_{\eg\in
  E}\DSN_{r_\eg}(\d(c^\eg,z^\eg))\int U_{\tau^{\eg_*}}(\d t^*) \tau^{\eg_*}
H(\mg,\lgo,t^*)\, ,$$ with the notation $U_n$ of Lemma
\ref{sec:discr-cont-snak-1}. Set
$$G_a(\br,\bx)=\bigotimes_{\eg\in E_{1/2}}{\bf P}^{r_\eg}_{0\to
  x_\eg}\bigotimes_{\eg\in E}\DSN_{r_\eg}\left[ U_{\tau^{\eg_*}}\left(
  \frac{\tau^{\eg_*}}{2a}H(\xi^a(\mg,\lgo,t^*))\right)\right]\, .$$
Because of the domination hypothesis on $H$, the fact
$\tau^{\eg_*}\leq \tau_\lgo$, and inequality (\ref{eq:16}), the
function $G_a$ satisfies the domination hypothesis of Lemma
\ref{sec:path-network-measure}. Applying the latter and Lemmas
\ref{sec:path-network-measure-2} and \ref{sec:discr-cont-snak-1}
allows us to conclude that for every $\mg\in\bM^{\geq 3}_{g,k}$,
$$(2a)^{-1}\frac{\varsigma^{k+2g-1}}{
  2^{\frac{|\bE|}{2}}a^{\frac{|\bE|+|V|-1}{4}}}
\LM_{\mg}(H\circ\xi^a)\build\to_{a\to\infty}^{} \CLM_{\mg}(H)
$$ Note that for $\mg$ ranging in $\bM^{\geq 3}_{g,k}$, the
quantity $|V|+|\bE|$ appearing in the scaling factor is maximal if and
only if $\mg$ belongs to the set $\bM^3_{g,k}$ of trivalent
maps of genus $g$ with $k$ faces. Indeed, $|\bE|-|V|=k+2g-2$ is a fixed
quantity, so that $|V|+|\bE|$ will be maximal if and only if $|\bE|$ is
maximal subject to the constraints. On the other hand, we already
noticed that $|\bE|\leq 3k+6g-6$, with equality if and only if
$\mg\in \bM^3_{g,k}$, in which case
$(|V|+|\bE|+3)/4=5g/2+(5k-7)/4$. Since $\bM^{\geq 3}_{g,k}$ is a finite
set, we deduce that the only terms in the sum (\ref{eq:18}) remaining
in the scaling limit will be those belonging to $\bM^3_{g,k}$.  Hence
the result. \cq

\bigskip

\noindent{\bf Remark: the Hurwitz measure. }It should be noted that a
non-labeled analog of $\CLM_{g,k}$ has already been considered in the
literature \cite{MR1802530,OkPan}, in problems dealing with the
enumeration of maps, random matrices, and related algebraic geometry
problems. Specifically, the first marginal $\mathsf{P}_t$ of the
measure $\SN_t$ can be considered as the law of the contour process of
a random forest. The measure considered in the above references is
$$\mathsf{CM}^\circ_{g,k}:=\sum_{\mg\in
  \bM^{3,\circ}_{g,k}}\frac{1}{|{\rm Aut}(\mg)|}\int_{\R_+^\bE}\d \br
\bigotimes_{\eg\in E}\mathsf{P}_{r_\eg}(\d c^\eg) ,$$ where the first
sum is over unrooted maps in $\bM^3_{g,k}$ (but still with indexed
faces) and the factor involving $\aut(\mg)$, the automorphisms of
$\mg$, allows one to get rid of symmetries of the graph. Note that we
can regroup the factors in the above measure as $\bigotimes_{\eg\in
  E_{1/2}}\d r_\eg\mathsf{P}_{r_\eg}(\d c^\eg)\mathsf{P}_{r_\eg}(\d
c^{\ov{\eg}})$. In this form, and by using an interpretation in terms
of trees of a well-known decomposition of the Brownian excursion due
to Bismut (see Le Gall \cite{legall93}), we see that we have replaced
each edge $\mg$ by a doubly-marked Brownian continuum random tree
taken under the Itô measure, using the two marks as endpoints of the
edge. The measure $\mathsf{CM}^\circ_{g,k}$ appears as a scaling limit
of a random element of $\bM_{g,k}$ where $k$ remains fixed but the
number of edges tends to infinity. Of particular importance is the
{\em Hurwitz measure}, which is the image of $\mathsf{CM}^\circ_{g,k}$
by the map
$$(c^\eg,\eg\in E)\longmapsto\left(\sum_{e\in f_i(\bm)}\tau(c^\eg),
1\leq i\leq k\right)\, \, ,$$ as its Laplace transform has an
interpretation in terms of generating series related to intersection
theory on the moduli space of curves \cite{OkPan}.

\subsection{Exceptional cases}\label{sec:exceptional-cases}

This section is devoted to the study of the planar cases
$(g,k)\in\{(0,1),(0,2)\}$. An element $(\bm,[\bl])$ of $\BLM_{0,1}$ is
just a rooted labeled tree and is described
by a snake $\Xi_{0,1}(\bm,[\bl])=(c,z)\in \bS_1$. The image measure of
$\QQ_{0,1}$ by $\wt{\Psi}_{0,1}:=\Xi_{0,1}\circ\Psi_{0,1}$ will assign
mass $2\cdot 12^{-n}=4\cdot 2^{-2n-1}3^{-n}$ to each discrete snake
$(c,z)$ with duration $\tau(c)=2n+1$, coding a labeled tree with $n$
edges. This is precisely the measure
$$\LM_{0,1}:=4\, \DSN_1\, .$$ We define its scaling limit $\CLM_{0,1}$
as follows. Let $N(\d c)$ be the Itô measure on $\W$, governing the
intensity of the excursions of standard reflected Brownian motion
(with normalization $N(\tau(c)\in \d t)=(2\pi t^3)^{-1/2}\d t$). Let
$\bN(\d(c,z))$ be the measure on $\W^2$ such that under $\bN$, $c$ has
``law'' $N$, and conditionally on $c$, $z$ is a centered Gaussian
process with covariance
$$\cov(z(s),z(s'))=\inf_{s\wedge s'\leq s''\leq s\vee s'}c(s'')\, .$$
Finally, set $\CLM_{0,1}=4\bN$.

\begin{prp}\label{sec:exceptional-cases-1}
Let $H$ be a continuous functional on $\W^2$ such that 
$$|H(c,z)|\leq K\tau\exp(-\beta \tau)\, ,$$ where $\tau:=\tau(c)$,
for some $K,\beta>0$. Then 
$$(2a)^{1/2}\LM_{0,1}(H\circ
(\vartheta^a,\psi^a))\build\longrightarrow_{a\to\infty}^{}
\CLM_{0,1}(H)\, .$$
\end{prp}

\proof First integrate out the label part, i.e.\ assume $H$ is only a
function of the first component and write
\begin{eqnarray*}
(2a)^{1/2}\LM_{0,1}(H(\vartheta^a(c)))&=&2(2 a)^{1/2}\int_{\R_+}\d t
  \frac{1}{2^{2\lceil t\rceil }}{\rm Cat}_{\lceil t\rceil }P^{\lceil
    t\rceil}_{{\rm exc}}(H(\vartheta^a(c)))\\ &=&2\sqrt{2}
  a^{3/2}\int_{\R_+}\d t \frac{1}{2^{2\lceil at\rceil}}{\rm
    Cat}_{\lceil at\rceil }P^{\lceil at\rceil}_{{\rm
      exc}}(H(\vartheta^a(c)))
\end{eqnarray*} 
where $P^{n}_{{\rm exc}}$ is the uniform distribution over Dyck paths
with duration $2n+1$, the latter set having cardinality ${\rm Cat}_n$,
the $n$-th Catalan number. It is well-known (and can be obtained from
Lemma \ref{sec:path-network-measure-2} with $x=-1$ from a cyclic shift
and Vervaat's theorem \cite{pitmancsp02}) that the law of
$\vartheta^a(c)$ under $P^{\lceil at\rceil}_{{\rm exc}}$ converges
weakly to $N^{(t)}=N(\cdot\, |\, \tau=t)$. On the other hand, it can
be checked that $4^{-n}(n+1)\sqrt{n}{\rm Cat}_n$ increases in $n$
towards its limit $\pi^{-1/2}$, so that $a^{3/2}2^{-2\lceil
  at\rceil}{\rm Cat}_{\lceil at\rceil}P^{\lceil at\rceil}_{\rm
  exc}(H(\vartheta^a(c)))$ is dominated by $K'e^{-2\beta t}/\sqrt{t}$
for some $K'>0$. By dominated convergence, we get
$$(2a)^{1/2}\LM_{0,1}(H\circ
\vartheta^a)\build\longrightarrow_{a\to\infty}^{}4\int_{\R_+}\frac{\d
  t}{\sqrt{2\pi t^3}}N^{(t)}(H)=4N(H)\, ,$$ by the fact that
$N(\tau\in \d t)=\d t/\sqrt{2\pi t^3}$. Incorporating the spatial
displacements is now easy from the convergence of the contour
function, see for instance \cite{janmarck05}. \cq

\bigskip

Finally, the measure associated with $\QQ_{0,2}$ through the mapping
$\wt{\Psi}_{0,2}$ is the measure $\LM_{0,2}=2\sum_{\mg\in\bM^{\geq
    3}_{0,2}}\LM_{\mg}$ where
$$\LM_\mg=\delta_\mg\WN_{\mg}(\d w^{\eg_*}\d w^{\ov{\eg}_*})
\DSN_1(\d(c,z))\#_\tau(\d t^*)\DSN_{r_{\eg^*}-1}(\d
(c^{\eg_*},z^{\eg_*}))\DSN_{r_{\eg^*}}(\d
(c^{\ov{\eg}_*},z^{\ov{\eg}_*}))\, .$$ Its continuous counterpart is
of course $\CLM_{0,2}=2\sum_{\mg\in \bM^{\geq 3}_{0,2}}\CLM_{\mg}$,
where
$$\CLM_{\mg}= \delta_\mg\PN_{\mg}(\d w^{\eg_*}\d w^{\ov{\eg}_*})
\bN(\d(c,z))\d t^*\ind_{\{0\leq t^*\leq \tau\}}\SN_{r_{\eg_*}}(\d
(c^{\eg_*},z^{\eg_*}))\SN_{r_{\eg^*}}(\d (c^{\ov{\eg}_*},z^{\ov{\eg}_*}))\, ,$$
  and $\PN_{\mg}$ takes the particularly simple form
$$\PN_{\mg}(\d w^{\eg_*} \d w^{\ov{\eg}_*}) =\int_{\R_+}\frac{\d
  r}{\sqrt{2\pi r}}P^r_{0\to 0}(\d w^{\eg_*} \d w^{\ov{\eg}_*})$$ We define
the scaling operation $\xi^a$ by
$$\xi^a(\mg,\lgo,t^*)=(\mg,(\vartheta^a(c),\psi^a(z)),(\phi^a(w^{\eg_*}),
  \vartheta^a(c^{\eg_*}),\psi^a(z^{\eg_*})),
  (\phi^a(w^{\ov{\eg}_*}),\vartheta^a(c^{\ov{\eg}_*}),\psi^a(z^{\ov{\eg}_*})),
  t^*/2a)\, .$$

\begin{prp}\label{sec:exceptional-cases-2}
For every continuous functional $H$ satisfying 
$$|H|\leq K\exp(-\beta\tau_\lgo)\, ,$$ for some
$K,\beta>0$. Then 
$$\frac{\varsigma}{2a^{3/4}}\LM_{0,2}(H\circ
\xi^a)\build\longrightarrow_{a\to\infty}^{}\CLM_{0,2}(H)\, .$$
\end{prp}

The proof uses both the ingredients of the proofs of Theorem
\ref{sec:cont-meas-label-1} and Proposition
\ref{sec:exceptional-cases-1}, but presents no new difficulty and is
left as an exercise to the reader.

\section{Asymptotic enumeration in the non-planar case}\label{sec:appl-asympt-enum}

We want to derive Theorem \ref{sec:main-results} by a study of the
generating function
$$G(s):=\sum_{n\geq 1}\frac{|\bQ_g^n|}{12^n}s^n\, ,$$ which will be
done by elementary singularity analysis. We will restrict our
attention to the non-planar cases $g\geq 1$, as the planar case is
already known.  First note that we have the following alternative
expressions
$$G(s)=\sum_{\bq\in\bQ_g}\left(\frac{s}{12}\right)^{|F(\bq)|}=
\sum_{(\bq,x)\in\bQ_{g,1}}V_\bq^{-1}\left(\frac{s}{12}\right)^{|F(\bq)|}
=\QQ_{g,1}(V_\bq^{-1}s^{|F(\bq)|})\, ,$$ which is finally equal to
$\LM_{g,1}((\tau_\lgo/2+\chi(g))^{-1}s^{\tau_\lgo/2})$. The
coefficients' asymptotic behavior is the same as that of
$G_0(s):=\LM_{g,1}(2\tau_\lgo^{-1}s^{\tau_\lgo/2})$, so we may and will
work with the latter.  Now write $s= e^{-1/a}$,
$H=\tau_\lgo^{-1}e^{-\tau_\lgo}$, and apply Theorem
\ref{sec:cont-meas-label-1}. We readily obtain that
$$G_0(e^{-1/a})\build\sim_{a\to\infty}^{} C_g' a^{\frac{5g-3}{2}}\,
,$$ where $C_g'\in (0,\infty)$ is defined by
$C_g'=2^{\frac{3g-1}{2}}3^g\CLM_{g,1}(\tau_\lgo^{-1}e^{-\tau_\lgo})$.
Otherwise said,
\begin{equation}\label{eq:24}
G_0(s)\sim C_g'(1-s)^{-(5g-3)/2}
\end{equation} 
as $s\to 1$ along $[0,1]$. From this, Tauberian theorems allow to
obtain
$$\sum_{k=1}^n12^{-n}|\bQ_g^n|\sim
\frac{C'_g}{\Gamma\left(\frac{5g-1}{2}\right)}n^{\frac{5g-3}{2}}\, ,$$
which is weaker than (but consistent with) what we are aiming at.  In
order to apply standard transfer theorems of singularity analysis, we
need to show that the asymptotics (\ref{eq:24}) hold for $s\to 1$ in a
domain of the form
$$D_{R,\phi}:=\{s\in \C:|s|<R,s\neq 1, |\arg(s-1)|>\phi\}
$$ for some $R>1,\phi\in(0,\pi/2)$ and on which $G_0$ has to be
extendible analytically (here we take the determination of the
argument in $(-\pi,\pi]$).

\begin{lmm}\label{sec:asympt-enum-non}
Let $\alpha> 0$ and $f(z)=\sum_{n\geq 0}h_n z^n$ be a power series
with non-negative coefficients satisfying $h_n\sim
n^{\alpha-1}/\Gamma(\alpha)$ as $n\to\infty$. Then $f(z)\sim
(1-z)^{-\alpha}$ for $z\to 1$ inside a domain of the form
$D_{1,\phi}$ with $\phi\in(\pi/2,\pi)$.
\end{lmm}

\proof Let $(1-z)^{-\alpha}=\sum_{n\geq 0}h^0_n z^n$, so that
$h^0_n\sim h_n$. Fix $\eps>0$, take $N$ large enough so that
$|1-h_n/h^0_n|\leq \eps$. From the fact that $h_n,h^0_n\geq 0$, we
easily get
$$\left|(1-z)^\alpha f(z)-1\right|\leq
|1-z|^\alpha\sum_{n=0}^N(h_n+h^0_n)|z|^n
+\eps\frac{|1-z|^\alpha}{(1-|z|)^\alpha}\, ,$$ and the first term goes
to $0$ as $z\to 1$. As for the second, write $1-z=\rho e^{\ii \theta}$
for $|\theta|<\pi-\phi$. For $\rho$ sufficiently close to $0$, we get
that $|1-z|/(1-|z|)\leq
\rho/(1-(1+\rho^2-2\rho\cos(\pi-\phi))^{1/2})$, which admits a finite
positive limit as $\rho\to 0$, allowing us to bound the right-hand
side of the last displayed expression by a constant ($\phi$-dependent)
multiple of $\eps$. Hence the result. \cq

\medskip

Now rewrite $G_0(s)$ as
$$G_0(s)=2\sum_{\mg\in \bM^{\geq
    3}_{g,1}}\int\wh{\WN}_\mg(\d\br,\d\bx)\bigotimes_{\eg\in
  E(\bm)}\DSN_{r_\eg}
\left(\frac{2\tau^{\eg_*}}{\tau_\lgo}s^{\tau_\lgo/2}\right)\, .$$ We use
the following trick: for a function $H$ on $\bM_{g,1}^{\geq 3}$ we
have
$$\sum_{\mg\in \bM_{g,1}^{\geq 3}}H(\mg)=\sum_{\mg^\circ\in
  \bM_{g,1}^{\geq 3,\circ}}|\aut(\mg^\circ)|^{-1}\sum_{\eg_*\in
  E(\mg^\circ)} H((\mg^\circ,\eg_*))\, ,$$ where $\bM_{g,1}^{\geq
  3,\circ}$ denotes unrooted maps and $(\mg^\circ,\eg_*)$ is the map
$\mg^\circ$ rooted at $\eg_*$.  Indeed, it holds that
$(\mg^\circ,\eg)=(\mg^\circ,\eg')$ as a rooted map for edges
$\eg,\eg'\in E(\mg)$ if and only if there exists an automorphism of
$\mg$ sending $\eg$ to $\eg'$. Thus, the orbits of the natural action
of $\aut(\mg^\circ)$ on $E(\mg)$ correspond to the distinct possible
rootings of $\mg^\circ$. Moreover, if $p\cdot \eg=\eg$, then $p$ is
the identity, which can be seen by reasoning in a step-by-step fashion
from the edge $\eg$ and going to neighbors.
Therefore, all the orbits have same cardinality $|\aut(\mg^\circ)|$,
and the claim follows. Applying it twice, once in each direction,
allows us to get (with a straightforward definition for
$\wh{\WN}_{\mg^\circ}$)
\begin{eqnarray*}
G_0(s)&=&4\sum_{\mg^\circ\in \bM^{\geq
    3,\circ}_{g,1}}|\aut(\mg^\circ)|^{-1}
\int\wh{\WN}_{\mg^\circ}(\d\br,\d\bx)\bigotimes_{\eg\in
  E(\mg^\circ)}\DSN_{r_\eg} \left(s^{\tau_\lgo/2}\right)\\ &=& 4\sum_{\mg\in
  \bM^{\geq
    3}_{g,1}}|E(\mg)|^{-1}\int\wh{\WN}_{\mg}(\d\br,\d\bx)\bigotimes_{\eg\in
  E(\mg)}\DSN_{r_\eg} \left(s^{\tau_\lgo/2}\right) \, .
\end{eqnarray*}
By (\ref{eq:27}) we have $\DSN_r(s^{\tau(c)})=((1-\sqrt{1-s^2})/s)^r$,
so that we get $G_0(s)=G_1\circ G_2(s)$, where
$$G_1(z)=4\sum_{\mg\in \bM^{\geq
    3}_{g,1}}|E(\mg)|^{-1}\int\wh{\WN}_{\mg}(\d\br,\d\bx)\,
z^{\sum_{\be\in \bE(\mg)}r_\be}\, ,$$ and
$$G_2(s)=\DSN_2(s^{\tau(c)/2})=\frac{2-s-2\sqrt{1-s}}{s}\, ,$$ which
can be extended analytically to a domain $D_{R,\phi}$ with $R>1$ and
$\phi\in(0,\pi/2)$. Moreover, elementary computations show that
$D_{R,\phi}$ is mapped to a domain $D_{1,\theta}$ for some
$\theta\in(\pi/2,\pi)$.

On the other hand, applying Lemma \ref{sec:key-lemmas-1} to $G_a=1$
and $\rho_a=n=\lceil (2a)^{1/2}\rceil$ gives that the $n$-th
coefficient in $G_1$ satisfies
$$[z^n]G_1(z)\sim 4\left(\frac{3}{2}\right)^g n^{5g-4}\sum_{\mg\in
  \bM^3_{g,1}}|E(\mg)|^{-1}\int_{\Delta_E}
\frac{\lambda_\Delta(\d\bu)}{\sqrt{{\cal Z}^\mg(\bu)}}=C_g''n^{5g-4}\,
.$$ Note that $|E(\mg)|=2(6g-3)$ for every term in the sum. From Lemma
\ref{sec:asympt-enum-non}, we deduce that $G_1(z)\sim
\Gamma(5g-3)C_g''(1-z)^{-(5g-3)}$ for $z\to 1$ inside
$D_{1,\theta}$. Combined with the fact that $G_2(s)=
1+2\sqrt{1-s}+O(1-s)$ as $s\to 1$ in $D_{R,\phi}$, we finally obtain
that $G_0$ has an analytic continuation in $D_{R,\phi}$, and is
equivalent to $2^{-5g+3}\Gamma(5g-3)C_g''(1-s)^{-(5g-3)/2}$ as $s\to
1$ inside this domain. An elementary transfer theorem \cite[Corollary
  VI.1]{FlSe} allows us to conclude that $[s^n]G_0(s)\sim
C_gn^{5(g-1)/2}$, where the constant $C_g$ is computed as
\begin{equation}\label{eq:25}
C_g=16\left(\frac{3}{64}\right)^g\frac{\Gamma(5g-3)}{(6g-3)
  \Gamma\left(\frac{5g-3}{2}\right)}\sum_{\mg\in
  \bM^3_{g,1}}\int_{\Delta_{\bE}}\frac{\lambda_\Delta(\d\bu)}{\sqrt{{\cal
      Z}^\mg(\bu)}}\, .
\end{equation}
Elementary (but tedious) computations give $C_1=1/24$, consistently
with \cite{BenCan86}. \cq

\medskip

\rems

\noindent{\bf 1.} It is expected, and checked on numerous examples, see Gao
\cite{gao93}, that ``reasonable'' families of random maps with size
$n$ will have the universal property that their cardinality is
equivalent to $C\alpha^n n^{-5\chi(g)/4}$ for some $C,\alpha$
depending on the model. A nice feature of the above approach is that,
let aside the prefactor $16(3/64)^g$, it is reasonable to think that
the gamma terms and summation over $\bM^3_{g,1}$ are not a particular
feature of quadrangulations, and will appear for more general families
of maps. However, it is extremely hard to determine the exact value of
this sum as it involves a number of terms growing very quickly with
$g$, and less and less tractable integrals.

\medskip

\noindent{\bf 2.} In \cite{ChMaSc}, the authors give an alternative
expression for $C_g$, still in the form of a prefactor times a sum
over trivalent maps with one face (and with an extra labeling of
vertices). The expression of \cite{ChMaSc} is somewhat simpler than
ours, as it does not involve integrals but only finite products of
quantities that depend in an elementary way on the labeled map
involved in the sum. It would be interesting to show directly that the
two expressions are the same.

\section{Metrics on metric spaces}
\label{sec:grom-hausd-weight}

In a metric space $(X,d)$, we let $B_X(x,r)$ (or simply $B(x,r)$ if
there is no ambiguity) be the open ball centered at $x$ with radius
$r$. We let $d_H$ be the usual Hausdorff distance between closed
subsets of $X$
$$d_H(C,C')=\inf\{\eps>0:C\subseteq (C')^\eps\mbox{ and }C'\subseteq
C^\eps\}\, ,$$ where $A^\eps=\{x\in X:\inf_{y\in A}d(x,y)<\eps\}$
is the $\eps$-neighborhood of $A$. We also
let $d_P$ be the Prokhorov metric between Borel probability measures
on $(X,d)$ defined by
$$d_P(\mu,\mu')=\inf\{\eps>0:\mu(C)\leq \mu'(C^\eps)+\eps\mbox{ for
  all closed }C\}\, .$$ The topology associated with $d_P$ is that of
weak convergence of probability measures on $X$. See \cite[Chapter
  11]{dudley02} for references.

\subsection{Gromov-Hausdorff distance}\label{sec:grom-hausd-dist}

We first recall the definition and first properties of the
Gromov-Hausdorff distance, taken from
\cite{gromov99,burago01,evpiwin}. If $\X=[X,d],\X'=[X',d']\in \M$, we
write
$$\dgh(\X,\X')=\inf_{\phi,\phi'}\delta_H(\phi(X),\phi'(X'))\, ,$$
where the infimum is taken over all isometries $\phi,\phi'$ from
$X,X'$ into a common metric space $(Z,\delta)$.

An alternative useful description is given as follows.  A
correspondence between two sets $X,X'$ is a subset $\RR\subseteq
X\times X'$ such that $\pi(\RR)=X$ and $\pi'(\RR)=X'$, where
$\pi:X\times X'\to X,\pi':X\times X'\to X'$ are the canonical
projections. We use $x\RR x'$ as an alternative notation for
$(x,x')\in \RR$, and we let $\CC(X,X')$ be the set of correspondences
between $X$ and $X'$. If $(X,d),(X',d')$ are metric spaces, and $\RR$
is a correspondence between $X$ and $X'$, the distortion of $\RR$ is
defined as
$$\dis\RR=\sup\{|d(x,y)-d'(x',y')|:x\RR x',y\RR y'\}\, .$$ 
Then
\begin{equation}\label{eq:2}
\dgh(\X,\X')=\frac{1}{2}\inf_{\RR\in \CC(X,X')}\dis \RR\, ,
\end{equation}
which of course does not depend on the particular representatives
$X,X'$ of $\X,\X'$.

It is an easy exercise to check that $\dis\ov{\RR}=\dis\RR$ where
$\ov{\RR}$ is the closure of $\RR$ in $X\times X'$ for the product
topology, so that the infimum in (\ref{eq:2}) could be taken over the
set $\CC_c(X,X')$ of compact correspondences (this set depends on the
distances $d,d'$ and not only on the sets $X,X'$). The following
properties \cite{burago01,evpiwin}, combined with the forthcoming
Proposition \ref{sec:relat-with-prev-1}, will be useful in proving
Theorem \ref{sec:mains-results-1}, in
Sect.\ \ref{sec:proof-theor-refs-1}.

\begin{prp}\label{sec:grom-hausd-dist-1}
{\rm (i)} The function $\dgh$ is a metric on $\M$. The space
$(\M,\dgh)$ is separable and complete.


{\rm (ii)} A subset $\bA\subseteq \M$ is relatively compact if and
only if the set $\{\diam(\X):\X\in \bA\}$ is bounded, and for every
$\eps>0$, there exists $N>0$ such that every $\X\in \bA$ can be
covered with $N$ open balls of radius $\eps$.
\end{prp}

\subsection{Gromov-Hausdorff-Prokhorov metric}\label{sec:weight-grom-hausd-3}

Let $\X=[X,d,\mu],\X'=[X',d',\mu']\in \mw$ We let
$$\dw(\X,\X')=\inf_{\phi,\phi'}\big(\delta_H(\phi(X),\phi'(X'))
\vee\delta_P(\phi_*\mu,\phi'_*\mu')\big)\, ,$$ where the infimum is taken
over all isometries $\phi,\phi'$ from $X,X'$ into a metric space
$(Z,\delta)$. 

\begin{thm}\label{sec:weight-grom-hausd}
The function $\dw$ defines a distance on $\mw$.
\end{thm}

This result is well-known, see for instance \cite[p.762]{villani09},
with the slight difference that the supremum between $\delta_H$ and
$\delta_P$ in the definition of $\dw$ is replaced by a sum. In order
to make this article self-contained, we will provide a full proof of
Theorem \ref{sec:weight-grom-hausd}, to which the remaining part of
this subsection is devoted.  In order to handle the distance $\dw$, we
develop an alternative definition with the help of correspondences and
couplings. This is inspired in part from \cite{gromov99,GPW06}.

If $\mu,\mu'$ are Borel probability measures on $(X,d)$ and $(X',d')$,
we say that a Borel measure $\nu$ on the product space $X\times X'$ is
a coupling between $\mu$ and $\mu'$ if $\pi_*\nu=\mu$ and
$\pi'_*\nu=\mu'$. We let $\MM(\mu,\mu')$ be the set of couplings
between $\mu$ and $\mu'$. It is never empty as it contains the product
measure $\mu\otimes\mu'$, and is closed (and even compact since
$X\times X'$ is compact) for the weak topology.

A useful feature of couplings is the following. Since all spaces
that are considered are Polish, it is known that we can write a
coupling $\nu\in \MM(\mu,\mu')$ in the form $\mu(\d x)Q_\nu(x,\d y)$,
where $Q_\nu$ is a Markov kernel from $X$ to $X'$. In particular, if
$\mu_1,\ldots,\mu_k$ are probability measures on the compact spaces
$X_1,\ldots,X_k$ and $\nu_i,1\leq i\leq k-1$ are couplings between
$\mu_i$ and $\mu_{i+1}$ respectively, we can define a Markov chain
$(A_1,\ldots,A_k)$ such that $A_i$ is a random variable with
distribution $\mu_i$ on $X_i$, with initial distribution $\mu_1$ and
transition kernel $Q_{\nu_i}$ at step $i$. The joint law of
$(A_1,A_k)$ is then denoted by $\nu_1\ldots\nu_k$ and is a coupling
between $\mu_1$ and $\mu_k$.

\begin{prp}\label{sec:weight-grom-hausd-2}
Let $\X=[X,d,\mu],\X'=[X',d',\mu']\in \mw$. Then 
\begin{equation}\label{eq:1}
\dw(\X,\X')=\inf\left\{\eps>0:\begin{array}{l}\exists
  \nu\in\MM(\mu,\mu'),\ \exists \RR\in\CC_c(X,X'),\\ 
\nu(\RR)\geq 1-\eps\quad \mbox{{\rm and}}\quad \frac{1}{2}\dis\RR\leq \eps
  \end{array}\right\}\, .
\end{equation}
\end{prp}

\proof Assume that $\dw(\X,\X')<\eps$ and let $\phi,\phi'$ be
isometric embeddings of $X,X'$ into $(Z,\delta)$ such that
$\delta_H(\phi(X),\phi'(X'))<\eps$ and
$\delta_P(\phi_*\mu,\phi'_*\mu')<\eps$. Arguing as in \cite[Remark
  7.3.12]{burago01}, we may and will assume that $(Z,\delta)$ is
separable. Then the set $\RR= \{(x,x')\in X\times
X':\delta(\phi(x),\phi'(x'))\leq \eps\}$ is an element of
$\CC_c(X,X')$. To evaluate its distortion, let $x\RR x'$ and $y\RR
y'$, and write
$$d(x,y)-d'(x',y')=\delta(\phi(x),\phi(y))-\delta(\phi'(x'),\phi'(y'))
\leq \delta(\phi(x),\phi'(x'))+\delta(\phi'(y'),\phi(y))\leq 2\eps\,
,$$ and the lower bound is obtained in a similar way. Thus
$\dis\RR\leq 2\eps$.

Now, by a characterization of the Prokhorov distance (see
\cite[Corollary 11.6.4]{dudley02}) on separable metric spaces, we can
find a coupling $\tilde{\nu}\in\MM(\phi_*\mu,\phi'_*\mu')$, which is a
measure on $Z\times Z$, such that
$$\tilde{\nu}(\{(x,y):\delta(x,y)\geq \eps\})\leq \eps\, .$$ We let
$\nu_1\in\MM(\mu,\phi_*\mu)$ be the image measure of $\mu$ by
$x\mapsto (x,\phi(x))$ and $\nu_2\in\MM(\phi'_*\mu',\mu')$ be the
image measure of $\mu'$ by $x'\mapsto (\phi'(x'),x')$.  Then the
coupling $\nu=\nu_1\tilde{\nu}\nu_2$ obviously satisfies
$\nu(\RR^c)=\tilde{\nu}(\{(x,y):\delta(x,y)\geq \eps\})\leq \eps$
(here $\RR^c$ is the complementary set of $\RR$).

Conversely, assume given $\RR\in\CC_c(X,X')$ and $\nu\in\MM(\mu,\mu')$
as in the right-hand side of (\ref{eq:1}) for the value $\eps$. As
in \cite{burago01} we endow $X\sqcup X'$ with a metric $\delta$ that is
compatible with $d,d'$, i.e.\ its restrictions to $X$ and $X'$
coincide with $d,d'$, and such that $\delta(x,x')=\eps$ whenever
$x\RR x'$. To do this, we set for every $x\in X,x'\in X'$
$$\delta(x,x')=\inf\{d(x,y)+\eps+d'(y',x'): y\RR y'\}\, .$$ The
only non-trivial point to show is the triangular inequality, and we
prove only the hardest case, the others being symmetric or left as an
exercise. Let $x,y\in X$ and $x'\in X'$, and consider
$z\RR z',t\RR t'$. Since by definition $|d(z,t)-d'(z',t')|\leq
2\eps$ we estimate
\begin{eqnarray*}
\lefteqn{(d(x,z)+\eps+d'(z',x'))+(d(y,t)+\eps+d'(t',x'))}\\
&\geq& (d(x,z)+d(z,t)+d(y,t))+(d'(z',x')+d'(t',x')-d'(z',t'))\\
&\geq & d(x,y) + 0\, ,
\end{eqnarray*}
which is obtained by two applications of the triangular inequality in
$X$ and $X'$. Taking the infimum over $(z,z'),(t,t')$ allows us to
conclude that $d(x,y)=\delta(x,y)\leq \delta(x,x')+\delta(x',y)$.

Now $\nu$ induces a coupling $\tilde{\nu}$ on the disjoint union
$X\sqcup X'$ between $i_*\mu$ and $i'_*\mu'$, where $i,i'$ are the
canonical injections $X\to X\sqcup X'$ and $X'\to X\sqcup X'$. Now
obviously for every $a>1$
$$\tilde{\nu}(\{(x,x'):\delta(x,x')\geq a\eps\})\leq \nu(\RR^c)\leq
\eps\, ,$$ so that $\delta_P(i_*\mu,i'_*\mu')\leq \eps$. Hence the
result. \cq

\medskip

We are now ready to prove that $\dw$ is indeed a distance. Symmetry is
obvious, as well as $\dw(\X,\X)=0$.

\noindent{\em Separation. } Let $\X=[X,d,\mu],\X'=[X',d',\mu']\in\mw$,
with $\dw(\X,\X')=0$. Then for every $n\geq 1$ we can find $\RR_n\in
\CC_c(X,X')$ and $\nu_n\in\MM(\mu,\mu')$ with $\dis\RR_n\leq 2^{-n}$
and $\nu_n(\RR_n^c)\leq 2^{-n}$. 

Endow the product space $X\times X'$ with the distance
$\delta((x,x'),(y,y'))=\max(d(x,y),d'(x',y'))$. Then $(X\times
X',\delta)$ is a compact metric space, so that the set of compact
subsets of $X\times X'$, endowed with the Hausdorff distance
$\delta_H$, is compact \cite[Theorem 7.3.8]{burago01}. Therefore,
$(\RR_n,n\geq 1)$ converges up to extraction to some compact $\RR$. It
is immediate to check that the latter is an element of $\CC_c(X,X')$,
with $\dis \RR=0$, and this implies that $\RR$ is the graph of a
bijective isometry $h:X\to X'$.

On the other hand, since $\MM(\mu,\mu')$ is compact for the weak
topology, the couplings $(\nu_n,n\geq 1)$ converge weakly, up to
re-extracting, to some $\nu\in \MM(\mu,\mu')$. Since $\RR_n$ converges
to $\RR$ for the distance $\delta_H$, it holds that for every
$\eps>0$, $\RR_n\subseteq \ov{\RR}^\eps=\{(x,x'):\inf_{(y,y')\in
  \RR}\delta((x,x'),(y,y'))\leq \eps\}$ for every $n$ large
enough. For such $n$, we have $\nu_n(\ov{\RR}^\eps)\geq
\nu_n(\RR_n)\geq 1-2^{-n}$, and by well-known properties of weak
convergence of probability measures, since $\ov{\RR}^\eps$ is closed,
this entails that $\nu(\ov{\RR}^\eps)\geq \limsup
\nu_n(\ov{\RR}^\eps)=1$. Letting $\eps\to 0$ shows that $\nu(\RR)=1$,
and we conclude that $\nu$ is supported on the graph of the isometry
$h$. Thus, for every measurable $g:X'\to \R_+$,
$$\int_Xg(h(x))\mu(\d x)=\int_{X\times X'}g(h(x))\nu(\d x,\d
x')=\int_{X\times X'}g(x')\nu(\d x,\d x')=\int_{X'}g(x')\mu'(\d x')\,
,$$ so $\mu'=h_*\mu$. Since $X'=h(X)$, we obtain $\X=\X'$.

\noindent{\em Triangular inequality. }
Let $\X=[X,d,\mu],\X'=[X',d',\mu'],\X''=[X'',d'',\mu'']\in \mw$, and let us assume $\dw(\X,\X')<\eps_1$ and
$\dw(\X',\X'')<\eps_2$. We wish to show that $\dw(\X,\X'')\leq
\eps_1+\eps_2$. By Proposition \ref{sec:weight-grom-hausd-2}, we can
find $\RR_1\in\CC_c(X,X'),\RR_2\in\CC_c(X',X'')$ and $\nu_1\in
\MM(\mu,\mu'),\nu_2\in\MM(\mu',\mu'')$ such that $\dis\RR_i\leq
2\eps_i$ and $\nu_i(\RR_i^c)\leq \eps_i,i\in\{1,2\}$.

Let $\nu_3=\nu_1\nu_2$ and define a correspondence $\RR_3\in
\CC(X,X'')$ by $x\RR_3 x''$ if and only if there exists $x'\in X'$
such that $x\RR_1x'$ and $x'\RR_2 x''$, that is $\RR_3$ is the image
by the (continuous) canonical projection $X\times X'\times X''\to
X\times X'$ of the compact set $(\RR_1\times X'')\cap(X\times \RR_2)$,
and hence is compact.
To evaluate
the distortion, let $x\RR_1 x',x'\RR_2 x''$ and $y\RR_1 y', y'\RR_2
y''$ and write
$$|d(x,y)-d''(x'',y'')|\leq
|d(x,y)-d'(x',y')|+|d'(x',y')-d''(x'',y'')| \leq 2(\eps_1+\eps_2)\,
,$$ so that $\dis\RR_3\leq 2(\eps_1+\eps_2)$.

Also note that $\ind_{(\RR_3)^c}(x,x'')\leq
\ind_{\RR_1^c}(x,x')+\ind_{\RR_2^c}(x',x'')$ for every $(x,x',x'')\in
X\times X'\times X''$. Indeed, if $(x,x'')\notin \RR_3$, then
for all $x'\in X'$, either $(x,x')\notin \RR_1$ or $(x',x'')\notin
\RR_2$. We deduce 
\begin{eqnarray*}\nu_3(\RR_3^c)&\leq & 
\int_{X\times X'\times X''}\mu(\d x)Q_1(x,\d x')Q_2(x',\d x'')
\ind_{\RR_1^c}(x,x')\\ & &+\int_{X\times X'\times X''}\mu(\d
x)Q_1(x,\d x')Q_2(x',\d x'') \ind_{\RR_2^c}(x',x'')\\ &=&
\nu_1(\RR_1^c)+\nu_2(\RR_2^c)\leq \eps_1+\eps_2\, .
\end{eqnarray*}
This is the wanted result, concluding the proof of Theorem
\ref{sec:weight-grom-hausd}. \cq

\subsection{Relation with other metrics}\label{sec:relat-with-prev}

It turns out that the topology determined by $\dw$ has already been
considered by \cite{fukaya,evanswinter}, where a different distance
was introduced. Note that in \cite{evanswinter}, the authors
restricted their attention to compact spaces which are $\R$-trees, but
all of their results up to Theorem 2.5 in this reference can be
rewritten for the space $\mw$ without any change.

An $\eps$-isometry from $X$ to $X'$
is a mapping $f:X\to X'$ such that $(f(X))^{\eps}=X'$ (we say that
$f(X)$ is an $\eps$-net in $X'$) and for every $x,y\in X$,
$$|d(x,y)-d'(f(x),f(y))|\leq \eps\, .$$ We say that
$\X=[X,d,\mu],\X'=[X',d',\mu']\in \mw$ are $\eps$-close if there exist two
measurable $\eps$-isometries $f:X\to X'$ and $g:X'\to X$ such that
$d_P'(f_*\mu,\mu')\vee d_P(g_*\mu',\mu)\leq \eps$.  Obviously, being
$\eps$-close is indeed an isometry class property and does not depend
on the choice of representatives.

We define $\Delta_{{\rm GHP}}(\X,\X')$ as the infimal $\eps$ such that
$\X,\X'$ are $\eps$-close. This does {\em not} define a true distance,
as the triangle inequality is satisfied only up to a factor of
$2$. This quasi-distance is however separated and thus sufficient to
define a separated topology on $\mw$. It is also fully discussed in
\cite{evanswinter}, building on metrization theorems for uniform
spaces \cite{bourbaki71}, how to build a true distance using
$\Delta_{{\rm GHP}}$, which is yet different from $\dw$.

\begin{prp}\label{sec:weight-grom-hausd-1}
It holds that $3^{-1}\Delta_{{\rm GHP}}\leq \dw\leq 2\Delta_{{\rm
    GHP}}$.
\end{prp}

\proof Assume that $\dw(\X,\X')<\eps$, and let $\RR,\nu$ be a
correspondence and a coupling between $\X$ and $\X'$ such that $\dis
\RR\leq 2\eps$ and $\nu(\RR)\geq 1-\eps$. We build a measurable
$3\eps$-isometry $f:X\to X'$ out of $\RR$, by taking
$\{x_1,\ldots,x_N\}$ a minimal $\eps$-net in $X$ , choosing for every
$i$ some $x'_i$ such that $x_i\RR x'_i$, and finally defining
$f(B_i)=\{x'_i\}$ where $B_i=B(x_i,\eps)\setminus \bigcup_{1\leq j\leq
  i-1}B(x_j,\eps)$. Then the image measure $\tilde{\nu}$ of $\nu$ by
$(x,x')\mapsto (f(x),x')$ is a coupling between $f_*\mu$ and
$\mu'$. Moreover, from the way $f$ is constructed, if
$d'(f(x),x')>3\eps$ then $(x,x')\notin \RR$, so that
$\tilde{\nu}\{(x',y'):d'(x',y')>3\eps\}\leq \nu(\RR^c)\leq \eps$. Thus
$d'_P(f_*\mu,\mu')\leq 3\eps$, and by a symmetrical argument we
conclude that $\X$ and $\X'$ are $3\eps$-close.

Conversely assume given $f:X\to X'$ an $\eps$-isometry with
$d'_P(f_*\mu,\mu')\leq \eps$. Define a correspondence between $X$ and
$X'$ by $x\RR x'$ if $d'(x',f(x))\leq \eps$. Its
distortion is estimated as follows:
\begin{eqnarray*}\lefteqn{|d(x,y)-d'(x',y')|}\\
&\leq&
|d(x,y)-d(f(x),f(y))|+|d(f(x),f(y))-d(f(x),y')|
+|d(f(x),y')-d(x',y')|\, ,
\end{eqnarray*} 
which is less than $4\eps$ whenever $x\RR x',y\RR y'$. We take the
closure of $\RR$, which has same distortion, and still call it
$\RR$. Take a coupling $\tilde{\nu}$ between $f_*\mu$ and $\mu'$ with
$\tilde{\nu}\{(x',y'):d'(x',y')\geq \eps\}\leq \eps$ and let $\nu_1$
be the image measure of $\mu$ by $x\mapsto (x,f(x))$. Then
$\nu=\nu_1\tilde{\nu}$ satisfies $\nu(\RR)\geq 1-\eps$. We conclude
that $\dw(\X,\X')\leq 2\eps$. \cq

As a consequence of this and \cite[Lemma 2.3]{evanswinter}, we obtain
that the metric properties of Evans and Winter's distance are the same
as those of $\dw$. In particular: 

\begin{prp}\label{sec:relat-with-prev-1}
{\rm (i)} The metric space $(\mw,\dw)$ is separable and complete. 

{\rm (ii)} A subset $\bA$ of $\mw$ is relatively compact if and only
if the subset $\{[X,d]:[X,d,\mu]\in \bA\}$ is relatively compact in
$(\M,\dgh)$.
\end{prp}

\subsection{Marked topologies}\label{sec:marked-topologies}

In this section we introduce a marked variant of the Gromov-Hausdorff
topology. 

A $k$-marked metric space is a triple $(X,d,(C_1,\ldots,C_k))$ where
$(X,d)$ is a compact metric space and $C_1,\ldots,C_k$ are compact
subsets of $X$. The isometry classes of marked spaces
$[X,d,(C_1,\ldots,C_k)]$ are defined in the obvious way, and their set
is denoted by $\M^*_k$. We let, with obvious notations,
$$\dgh^k(\X,\X')=\inf_{\phi,\phi'}
\left(\delta_H(\phi(X),\phi'(X'))\vee \max_{1\leq i\leq
  k}\delta_H(\phi(C_i),\phi'(C_i'))\right)\, ,$$ the infimum being
taken over all isometric embeddings $\phi,\phi'$ of $\X,\X'$ into a
common metric space $(Z,\delta)$.

The following statement is a straightforward generalization of the
classical results of Proposition \ref{sec:grom-hausd-dist-1}, and its
proof is left as an exercise to the reader.

\begin{prp}\label{sec:marked-topologies-1}
{\rm (i)} The function $\dgh^k$ is a distance on $\M^*_k$, which is
alternatively described by
$$\dgh^k(\X,\X')=\frac{1}{2}\inf_{\RR}\dis \RR\, ,$$ the infimum being
taken over correspondences $\RR\in \CC_c(X,X')$ such that $\RR\cap
(C_i\times C_i')\in \CC(C_i,C_i')$.

{\rm (ii)} The metric space $(\M^*_k,\dgh^k)$ is separable and
complete. 

{\rm (iii)} A subset $\bA\subseteq \M^*_k$ is relatively compact if and
only if $\{[X,d]:[X,d,(C_1,\ldots,C_k)]\in A\}$ is relatively compact
in $(\M,\dgh)$.
\end{prp}

When $C_1,\ldots,C_k$ are singletons, we
write for notational convenience $(X,d,\bx)$ instead of
$(X,d,(\{x_1\},\ldots,\{x_k\}))$. The reader will easily be convinced
that the subspace of such marked spaces is closed in $\M^*_k$.

\subsection{Randomly marked spaces}\label{sec:rand-mark-spac}

We now introduce random marks on metric spaces. If $\X=[X,d,\mu]\in
\mw$, we can consider an i.i.d.\ sequence $(U_1,\ldots,U_k)$ which is
$\mu$-distributed in $X$, and we want to consider the isometry class
of $[X,d,(U_1,\ldots,U_k)]$ as a random variable in $\M^*_k$. We also
wish to consider randomly marked versions of spaces that are
themselves random. To do this implies that we check a couple of
measurability issues. We introduce for each $k$ a kernel from $\mw$ to
$\M^*_k$, defined as follows:
\begin{equation}\label{eq:3}\mm_k(\X,\bA)=\int_{X^k}
\mu^{\otimes k}(\d\bx)\ind_\bA([X,d,\bx])\, ,\qquad \X=[X,d,\mu]\in
\mw, \bA\in{\cal B}(\M^*_k)\, .
\end{equation} 
Note that the integral makes sense as the function $X^k\to \M^*_k$ which
maps $\bx$ to $[X,d,\bx]$ is obviously continuous, so that
$\bx\mapsto \ind_{\bA}([X,d,\bx])$ is measurable by composition. Also
notice that it is independent of the particular choice of the
representative $(X,d,\mu)\in \X$, and hence is unambiguous.

\begin{lmm}\label{sec:rand-mark-spac-1}
The formula (\ref{eq:3}) defines a Markov kernel.
\end{lmm}

\proof It is obvious from the definition that for every $\X$, the
function $\bA\mapsto \mm_k(\X,\bA)$ defines a measure on $\M^*_k$.

Now let $\bA$ be closed in $\M^*_k$. We show $\X\mapsto \mm_k(\X,\bA)$
is upper-semicontinuous. Taking $\eps>0$, for any $\X'\in\mw$ with
$\dw(\X,\X')<\eps$, and with obvious notations, we can find
$\RR\in\CC_c(X,X')$ and $\nu\in\MM(\mu,\mu')$ such that $\nu(\RR)>
1-\eps$ and $\dis \RR< 2\eps$. We write
\begin{eqnarray*}
\lefteqn{\int_{X^k}\mu^{\otimes
    k}(\d\bx)\left(1-\ind_{\bA^{\eps}}([X,d,\bx]) \right)
  =\int_{(X\times X')^k}\nu^{\otimes
    k}(\d\bx,\d\bx')\left(1-\ind_{\bA^{\eps}}([X,d,\bx])\right)}
\\ &\leq & 1-(1-\eps)^k+\int_{(X\times X')^k}\nu^{\otimes
  k}(\d\bx,\d\bx')\ind_{\{x_i\RR x'_i,1\leq i\leq
  k\}}\left(1-\ind_{\bA^{\eps}}([X,d,\bx])\right) \\ &\leq &
1-(1-\eps)^k+\int_{(X')^k}\d(\mu')^{\otimes
  k}(\bx')\left(1-\ind_\bA([X',d',\bx'])\right)\, ,
\end{eqnarray*}
where we used the fact that $\nu^{\otimes k}(\{(\bx,\bx'):x_i\RR
x'_i,1\leq i\leq k\})\geq (1-\eps)^k$ and that whenever the latter
event occurs, it holds that $\dgh^k([X,d,\bx],[X',d',\bx'])<
\eps$. Since $\bA^{\eps}$ decreases to $\bA$ as $\eps\downarrow 0$ we
see that the left-hand side can be made arbitrarily close to
$1-\mm_k(\X,\bA)$ provided $\eps$ is sufficiently small.  Therefore
$$\limsup_{\X'\to\X}\mm_k(\X',\bA)\leq \mm_k(\X,\bA)\, ,$$ as claimed.
We conclude by a monotone class argument that $\mm_k(\cdot,\bA)$ is
measurable for all Borel $\bA\subseteq \M^*_k$. \cq

A random variable with law $\mm_k(\X,\cdot)$ can be considered as the
rigorous definition of the space $\X$ randomly marked with $k$
independent $\mu$-chosen points. If ${\bf M}(\d\X)$ is a distribution
on $\mw$, a random variable with law ${\bf M}(\d\X)\mm_k(\X,\cdot)$ is
interpreted as a weighted space chosen according to ${\bf M}$, and
then marked randomly using the weight. 

A key property of randomly marked spaces is the following. If $\X\in
\mw$, we let $\X^k$ be the canonical random variable with law
$\mm_k(\X,\cdot)$. A similar notation holds if $\X$ is a random
variable with values in $\mw$.

\begin{prp}\label{sec:rand-mark-spac-2}
Let $(\X_n,n\geq 1)$ be random variables in $\mw$ converging to $\X$
in distribution.  Then the randomly marked spaces $\X_n^k$ converge in
distribution towards $\X^k$, for the topology of $\M^*_k$.
\end{prp}

\proof Using the Skorokhod representation theorem, we may assume that the
probability space $(\Omega,{\cal F},P)$ on which the random variables
are defined is such that the convergence $\X_n\to \X$ holds
a.s. Now, we have from the proof of Lemma \ref{sec:rand-mark-spac-1}
that $\mm_k(\cdot,\bA)$ is upper-semicontinuous for every closed
$\bA\subseteq\M^*_k$. By the reversed Fatou lemma this gives
$$\limsup_{n\to\infty} E[\mm_k(\X_n,\bA)]\leq E[\mm_k(\X,\bA)]\, ,$$
which says that $\limsup P(\X_n^k\in\bA)\leq P(\X^k\in \bA)$, and the
conclusion follows from a well-known characterization of weak
convergence. \cq

\subsection{Three lemmas}\label{sec:three-lemmas}

We end this section with key lemmas to the proof of Theorem
\ref{sec:mains-results-2}.

\subsubsection{Diffuse measures}\label{sec:diffuse-measures}

\begin{lmm}\label{sec:diffuse-measures-1}
Let $[X_n,d_n,\mu_n],n\geq 1$ converge in $\mw$ towards
$[X,d,\mu]$. Assume that
$$\lim_{\eps\downarrow 0}\liminf_{n\to\infty}\int_{X_n}\mu_n(\d x)
\mu_n(B(x,\eps))=0\, .$$ Then $\mu$ is a diffuse measure.
\end{lmm}

\proof From the proof of Lemma \ref{sec:rand-mark-spac-1}, we know
that if $\bA$ is an open set in $\M^*_2$, then $\X\mapsto
\mm_2(\X,\bA)$ is lower-semicontinuous from $\mw$ to $\M^*_2$. Now
write
$$\int_{X_n}\mu_n(\d x) \mu_n(B(x,\eps))=\int_{X_n^2}\mu_n^{\otimes
  2}(\d x,\d y)\ind_{\{d_n(x,y)<\eps\}} =\mm_2(\X_n,\bA)\, ,$$ where
$\bA=\{\X'=[X',d',(x',y')]\in \M^*_2:d'(x',y')<\eps\}$ is
open. Therefore, it holds that
$$\liminf_{n\to\infty}\int_{X_n}\mu_n(\d x) \mu_n(B(x,\eps)) \geq
\int_X\mu(\d x) \mu(B(x,\eps))\, .$$ Letting $\eps\to 0$ and using the
hypothesis gives the result. \cq

\subsubsection{Support}\label{sec:support}

\begin{lmm}\label{sec:marked-topologies-2} Let $[X_n,d_n,\mu_n],n\geq
  1$ be a sequence converging to $[X,d,\mu]$ in $\mw$. Let
  $Y_n\subseteq X_n$ be $\eps_n$-dense in $X_n$, for some $\eps_n\to
  0$ (i.e.\ every point of $X_n$ is at distance at most $\eps_n$ of a
  point of $Y_n$). Assume that for every $\eps>0$,
$$\limsup_{n\to\infty}\inf_{x\in Y_n}\mu_n(B(x,\eps))>0\, .$$ Then
  $\supp \mu=X$. 
\end{lmm}

\proof Take $x\in X$ and $\eps>0$. We must show that
$\mu(B(x,\eps))>0$. To this end take compact correspondences $\RR_n$
between $X_n,X$ with distortion going to $0$, and couplings $\nu_n$
between $\mu_n$ and $\mu$ with $\nu_n(\RR_n)\to 1$. Let $x_n$ be such
that $x_n\RR_n x$ and $x'_n\in Y_n$ be at distance $\leq \eps_n$ from
$x_n$. Then
$$\mu(B(x,\eps))=\int_X\mu(\d y)\ind_{\{d(x,y)<\eps\}}\geq
\int_{X_n\times X}\nu_n(\d y_n,\d y)\ind_{\{d(x,y)<\eps\}}\, .$$ Now
take $n$ large enough so that $\dis \RR_n<\eps/2-\eps_n$. Then
$d_n(x_n',y_n)<\eps/2$ implies $d(x,y)<\eps$ as soon as $y_n\RR_n
y$. Therefore the last integral is greater than or equal to
\begin{eqnarray*}
\int_{\RR_n}\nu_n(\d y_n,\d y)\ind_{\{d_n(x'_n,y_n)<\eps/2\}} &\geq&
\int_{X_n\times X}\nu_n(\d y_n,\d y)\ind_{\{d_n(x'_n,y_n)<\eps/2\}}
-\nu_n(\RR_n^c)\\ &=&\mu_n(B(x'_n,\eps/2))-\nu_n(\RR_n^c)\, .
\end{eqnarray*} 
This is bounded below by $\inf_{x\in
  Y_n}\mu_n(B(x,\eps/2))-\nu_n(\RR_n^c)$, which has positive $\limsup$
as $n\to\infty$. Hence the result. \cq

\subsubsection{$D$-median points in geodesic metric
  spaces}\label{sec:d-interm-points}

Since compact geodesic metric spaces form a closed subset $\pmc$ of
$(\M,\dgh)$ (\cite[Thorem 7.5.1]{burago01}), it is immediate that the
set $\pmw$ of spaces $[X,d,\mu]\in\mw$ for which $[X,d]$ is a geodesic
metric space is closed as well, and a similar statement holds for
$k$-pointed geodesic spaces, whose set is denoted by $\pmc^*_k$.

Let $(X,d)$ be a compact geodesic metric space, and let $x,y\in X$ and
$D\in(-d(x,y),d(x,y))$. A point $z\in X$ is called a $D$-median point
between $x$ and $y$ if $d(x,z)-d(y,z)=D$. Note that the set ${\rm
  med}^D_{xy}(X)$ of such points is compact and separates the space
$X$ into two open subsets, one containing $x$ and the other $y$, since
every path $c:[0,1]\to X$ from $x$ to $y$ is such that
$d(c(t),x)-d(c(t),y)$ ranges from $-d(x,y)$ to $d(x,y)$
continuously. 
We also say that $z$ is on a $\delta$-quasi geodesic from $x$ to $y$
if $d(x,z)+d(y,z)\leq d(x,y)+\delta$. We let ${\rm qg}^\delta_{xy}(X)$
be the set of such points, and ${\rm medqg}^{D,\delta}_{xy}(X)={\rm
  med}^D_{xy}(X)\cap {\rm qg}^\delta_{xy}(X)$. 

\begin{prp}\label{sec:d-interm-points-1}
Let $[X_n,d_n,(x_n,y_n)]\in\pmc^*_2$ be converging to
$[X,d,(x_1,x_2)]$, and let also $D_n\in(-d_n(x_n,y_n),d_n(x_n,y_n))$
be converging to $D\in(-d(x,y),d(x,y))$. We assume that for every
$\eps>0$, we can find $\delta>0$ such that,
$$\limsup_{n\to\infty}\, \diam \left({\rm
  medqg}^{D_n,\delta}_{x_ny_n}(X_n)\right) < \eps\, .$$ Then there
exists a unique element of ${\rm med}^D_{xy}(X)$ which is on a
geodesic path from $x$ to $y$, i.e.\ ${\rm medqg}^{D,0}_{xy}(X)$ is a
singleton.
\end{prp}

\proof Assume that there are two possibly different points $z, z'\in
{\rm med}^D_{xy}(X)$ that belong to some geodesic path between $x$ and $y$.
Choose arbitrarily such a path $\gamma=\gamma_{xy}$ going through $z$
and let $\eps>0$. Let $\eps>\delta>0$ be as in the hypothesis for this
value of $\eps$, and let $0=t_0<t_1<\ldots<t_K=d(x,y)$ be such that
$d(u^k,u^{k+1})\leq \delta/7$ where $u^k=\gamma(t_k)$. Let $\RR_n$ be
a correspondence between $X_n$ and $X$ such that $x_n\RR_n x$ and
$y_n\RR_n y$, and $\dis\RR_n\to 0$ as $n\to\infty$. We may assume that
$n$ is large enough so that $\dis\RR_n\leq \delta/7$. Choose
arbitrarily $u^k_n\RR_n u^k$ and let $\gamma_n$ be the concatenation
of arbitrarily chosen geodesic paths between $u^k_n$ and $u^{k+1}_n$
for $0\leq k\leq K-1$.

Take a point $v$ on $\gamma_n$. Then for some $k$, $v$ lies on a
geodesic from $u^k_n$ to $u^{k+1}_n$. Then
\begin{eqnarray*}
d_n(v,x_n)&\leq& d_n(v,u^k_n)+d_n(u^k_n,x_n)\\ &\leq &
d_n(u^{k+1}_n,u^k_n)+d_n(u^k_n,x_n)\\ &\leq &
d(u^{k+1},u^k)+d(u^k,x)+2\dis\RR_n\\ &\leq &
d(u^k,x)+2\dis\RR_n+\delta/7\leq d(u^k,x)+3\delta/7\, ,
\end{eqnarray*}
so that by a symmetric inequality involving $y_n$ and by summing,
$$d_n(v,x_n)+d_n(v,y_n)\leq d(x,y)+6\delta/7 \leq
d_n(x_n,y_n)+\delta\, .$$ Therefore the image of $\gamma_n$ is
included in ${\rm qg}_{x_ny_n}^\delta(X_n)$. Let $z_n$ be any element
of ${\rm med}^{D_n}_{x_ny_n}(X_n)$ that belongs to $\gamma_n$, and let
$\tilde{z}_n\in X$ such that $z_n\RR_n \tilde{z}_n$.

Now doing the same with a geodesic path $\gamma'$ from $x$ to $y$
passing through $z'$, we can construct elements $z'_n\RR_n
\tilde{z}'_n$ such that $z'_n\in {\rm med}^{D_n}_{x_ny_n}(X_n)\cap{\rm
  qg}^\delta_{x_ny_n}(X_n)$.  Hence $d_n(z_n,z'_n)\leq \eps$ for large
$n$ by hypothesis, and thus $d(\tilde{z}_n,\tilde{z}'_n)\leq
\eps+\dis\RR_n$.

Letting $n$ go to infinity, possibly through a subsequence, we obtain
that $\tilde{z}^n\to \tilde{z}(\eps)$ and
$\tilde{z}'_n\to\tilde{z}'(\eps)$, where necessarily
$\tilde{z}(\eps),\tilde{z}'(\eps)\in {\rm med}^D_{xy}(X)$ by passing
to the limit in
$$|d(\tilde{z}_n,x)-d(\tilde{z}_n,y)-D_n|\leq 2\dis \RR_n\, ,$$ and
similarly for $\tilde{z}'_n$. By a similar argument, since $z_n\in
\gamma_n$ we have $\tilde{z}(\eps)\in \gamma^\delta$, the
$\delta$-neighborhood of the image of $\gamma$.  Since
$\gamma^{\delta}$ decreases to $\gamma$ as $\delta\downarrow 0$ and
${\rm med}^D_{xy}(X)$ is compact, letting $\eps$ (and $\delta$)
decrease to $0$, $\tilde{z}(\eps)$ converges to the unique element of
$\gamma\cap {\rm med}^D_{xy}(X)$, which is $z$. Similarly,
$\tilde{z}'(\eps)$ can be made as close to $z'$ as wanted.  But since
$d(\tilde{z}_n,\tilde{z}_n')\leq \eps+\dis\RR_n$ we get that
$d(\tilde{z}(\eps),\tilde{z}'(\eps))\leq \eps$, and this entails
$z=z'$. \cq

\section{Metric aspects of random quadrangulations}\label{sec:proof-main-results}

We now embark in the proof of Theorems \ref{sec:mains-results-1} and
\ref{sec:mains-results-2}. We start with a preliminary remark. For
technical reasons, the set $\pmc$ and its weighted and marked analogs
are somewhat easier to work with for our purposes. Of course, the main
space of interest, namely $\X_\bq$ for a quadrangulation $\bq$, is not
a geodesic metric space. However, its $\dgh$-distance to $\pmc$ is at most
$1$. Indeed, take copies of the unit segment (with the usual metric)
$\{x^e:0\leq x\leq 1\}$ indexed by an orientation $E_{1/2}(\bq)$ of
the edges of $\bq$, and identify the points $0^e$ for which $e^-$
coincide, and the points $1^e$ for which $e^+$ coincide. The quotient
metric graph \cite[Chapter 3.2.2]{burago01} is a geodesic metric space
whose restriction to the points $0^e,e\in E_{1/2}(\bq)$ is isometric
to $(V(\bq),d_\bq)$, and obviously its Gromov-Hausdorff distance to
the latter is less than $1$. Thus any limit in distribution of
$a^{-1/4}\X_\bq$ must be at $0$ distance from $\pmc$ with probability
$1$, hence belong to $\pmc$. In passing this proves property 1.\ in
Theorem \ref{sec:mains-results-2}. For these reasons, we work with the
metric graph constructed above, which we still call $\X_\bq$, and
restrict our attention to geodesic metric spaces from now on.

\subsection{Proof of Theorem \ref{sec:mains-results-1}}
\label{sec:proof-theor-refs-1}

For $\X\in \M$ and $\eps>0$, we let $\n(\X,\eps)$ be the minimal
number of open balls of radius $\leq \eps$ that are needed to cover
$\X$. For any sequence $N=(N_k,k\geq 1)$ of positive integers, we
introduce the set
$$\K_N=\{\X^{\rm w}=[X,d,\mu]\in \pmw:\n([X,d],2^{-k})\leq
N_k\ \forall k\geq 1\}\, .$$ By Propositions
\ref{sec:relat-with-prev-1} and \ref{sec:grom-hausd-dist-1} and
\cite[Exercise 7.4.14]{burago01}, the sets $\K_N$ are (relatively)
compact. By Prokhorov's theorem, it is thus sufficient to show that
for every $\eps>0$, there is a choice of $N,\Delta$ for which
$$\QQ^{(\beta/a)}_g\left(\left\{a^{-1/4}\X_\bq^{\rm w}\in
\K_N,\Delta^{-1}\leq a^{-1}V_\bq\leq \Delta\right\}\right)\geq
1-\eps\, .$$ To do so, we use the results of Sections
\ref{sec:coding-delay-voron}--\ref{sec:discr-cont-path} in the case
$k=1$. We also separate the case $g\geq 1$ and $g=0$, starting with
the former. Consider $(\bq,x)\in\bQ_{g,1}$ and let
$(\bm,[\bl])=\Psi_{g,1}(\bq,x)$ and
$(\mg,\lgo,t^*)=\wt{\Psi}_{g,1}(\bq,x)$. Recall also the definition
(\ref{eq:21}) of $l$, and write
$$\omega(h,\eps)=\sup_{x,y\in J,|x-y|\leq \eps}|h(x)-h(y)|$$ for the
modulus of continuity of a function $h$ defined on an interval $J$.

\begin{lmm}\label{sec:proof-theor-refs-4}
For $R>0$,
\begin{equation}\label{eq:15}
\left\{ {\cal N}(\X_\bq,R)> N\right\}\subset
\left\{2\omega\left(l,
\left\lceil\frac{\tau_\lgo}{N}\right\rceil\right)\geq R\right\}
\end{equation}
\end{lmm}

\proof Since $\bm$ has only one face of degree $\tau_\lgo$, all the
half-edges of $\bm$ are of the form $e(j)=\varphi_\bm^j(\eg_*(0))$ for
some integer $j\in[0,\tau_\lgo)$ with the notations of
  Sect.\ \ref{sec:map-reductions-1}. Split the lifetime interval
  $[0,\tau_\lgo]$ of $l$ into $N$ disjoint parts $I_1,\ldots,I_N$ of
  lengths at most $\lceil \tau_\lgo/N\rceil$.  By Lemma
  \ref{sec:bounds-distances-1} and the definition (\ref{eq:20}) of
  $l$, if $\omega(l,\lceil \tau_\lgo/N\rceil)<R/2$, then all the
  vertices of $\bm$ of the form $e(j)^-$ for $j\in I_i$, and
  considered as vertices of $\bq$, must be within $d_\bq$-distance at
  most $R$, for each $1\leq i\leq N$. Since all vertices of $\bq$ are
  contained in this way except the distinguished vertex $x$, adding yet
  another ball to cover $x$ yields (\ref{eq:15}). \cq

\bigskip

We now consider a marked version of $\QQ^{(\beta)}_g$, as the
probability on $\bQ_{g,1}$ defined by
$$\QQ^{(\beta)}_g(\d\bq)\mu_\bq(\d x)
=\frac{\QQ_{g,1}(V_\bq^{-1}e^{-\beta V_\bq}\,
  \d(\bq,x))}{\QQ_{g,1}(V_\bq^{-1}e^{-\beta V_\bq})}\, .$$ 
Using Section \ref{sec:discr-cont-path}, its image under
$\wt{\Psi}_{g,1}$ can be written
$$\LM_{g,1}^{(\beta)}:=\frac{\LM_{g,1}\left( \left(\tau_\lgo/2+
  \chi(g)\right)^{-1}\exp\left(-\beta\tau_\lgo/2\right)\,
  \d(\mg,\lgo,t^*)\right)}{\LM_{g,1}\left( \left(\tau_\lgo/2+\chi(g)
  \right)^{-1}\exp\left(-\beta\tau_\lgo/2\right)\right)}\, .$$ Now,
Theorem \ref{sec:cont-meas-label-1} (and elementary bounds to get rid
of the term $\chi(g)$) entail that the law of $\xi^a(\mg,\lgo,t^*)$
under $\LM_{g,1}^{(\beta/a)}$ converges weakly towards
\begin{equation}\label{eq:9}\CLM_{g,1}^{(\beta)}:=\frac{\CLM_{g,1}\left(
  \tau_\lgo^{-1}\exp\left(-\beta \tau_\lgo \right)\,
  \d(\mg,\lgo,t^*)\right)}{\CLM_{g,1}\left(
    \tau_\lgo^{-1}\exp\left(-\beta \tau_\lgo \right)\right)}\, .
\end{equation} 
By the continuous mapping theorem, the distribution of
$(\psi^a(l^\eg),\eg\in E(\mg))$ under $\LM^{(\beta/a)}_{g,1}$
jointly converge towards the law of the processes $(l^\eg,\eg\in
E(\mg))$ (\ref{eq:13}) under $\CLM^{(\beta)}_{g,1}$, and thus the
law of $\psi^a(l)$ under $\LM^{(\beta/a)}_{g,1}$ converges towards
the law of $l$ (\ref{eq:21}) under $\CLM^{(\beta)}_{g,1}$.  Combined
with the fact from Lemma \ref{sec:proof-theor-refs-4} that
$$\left\{ {\cal N}(a^{-1/4}\X_\bq,\delta)> N\right\}\subset
\left\{2\omega\left(\psi^a(l), (2a)^{-1}\left\lceil
\frac{\tau_\lgo}{N}\right\rceil \right)\geq \delta/\varsigma\right\}\, ,$$ and
by usual relative compactness criteria derived from the Ascoli-Arzela
theorem \cite[Theorem 8.2]{billingsley99}, we see that for every
$k\geq 1$, there exists $N_k$ such that for every $a$,
$\QQ^{(\beta/a)}_g ({\cal N}(a^{-1/4}\X_\bq,2^{-k})> N_k)\leq \eps
2^{-k}$. This entails $\QQ^{(\beta/a)}_g(a^{-1/4}\X_\bq\in
\mathbb{K}_N)\geq 1-\eps$ for this choice of $N=(N_k,k\geq 1)$. We see
as well that the random variables $V_\bq/a$ are tight in $(0,\infty)$
as they converge in distribution towards the law of $\tau_\lgo$ under
$\CLM_{g,1}^{(\beta)}$. \cq

\medskip

The case $g=0$ is similar but simpler, the role of $\CLM_{0,1}$ being
performed by the Itô measure of the Brownian snake, and taking into
account the factor $V_\bq^2$ appearing in the definition of
$\QQ^{(\beta)}_0$, we must replace the limit in law of
$\psi^a(l)$ by the law of $z$ under
$$\frac{\CLM_{0,1}(\tau\exp(-\beta\tau)\,
  \d(c,z))}{\CLM_{0,1}(\tau\exp(-\beta\tau))}\, .
$$ One will notice that the denominator would diverge if it were not
for the corrective term $\tau$ in front of the exponential, since
$\CLM_{0,1}(\tau\in \d t)=4(2\pi t^3)^{-1/2}\d t$. The statement of
Lemma \ref{sec:proof-theor-refs-4} remains true without change, and we
conclude in the same way.

\subsection{Proof of Theorem \ref{sec:mains-results-2}}
\label{sec:proof-theor-refs-2}

Property 1.\ of this theorem has been discussed at the beginning of
this section. Let us set things up for Property 2. We give the proofs
for $g\geq 1$, the case $g=0$ being similar and easier.

\subsubsection{Diffuseness of $\mu$}\label{sec:support-diffuseness}

By Theorem \ref{sec:mains-results-1}, there exists an increasing
sequence $(a_n,n\geq 1)$ along which the distribution under
$\QQ_g^{(\beta/a)}$ of $(a^{-1/4}\X_\bq^{\rm w},a^{-1}V_\bq)$
converges to a limit $\mathscr{S}_g^{(\beta)}$. From now on, we will
always suppose that $a$ is taken along this sequence. From this and
the discussion above (\ref{eq:9}), we obtain that the laws of
$(a^{-1/4}\X_\bq^{\rm w},a^{-1}V_\bq,\xi^a(\wt{\Psi}_{g,1}(\bq,x)))$
under $\QQ_g^{(\beta/a)}(\d\bq)\mu_\bq(\d x)$ form a relatively
compact family of probability distributions on $\pmw\times
(0,\infty)\times \bC_{g,1}$. Up to a further extraction, and by a use
of Skorokhod's embedding theorem, we may and will assume that we are
working on a probability space $(\Omega,\FF,P)$ under which are
defined random variables $(\bq_a,x_a)$ with respective laws
$\QQ^{(\beta/a)}_g(\d\bq)\mu_\bq(\d x)$ such that
\begin{enumerate}
\item
letting $\X_a^{\rm w}=[X_a,d_a,\mu_a]:=\X_{\bq_a}^{\rm w}$, we have
$a^{-1/4}\X_a^{\rm w}\to \X^{\rm w}=[X,d,\mu]$, a random variable with
values in $\pmw$ (and we let $\X_a=[X_a,d_a]$),
\item
letting $(\bm_a,[\bl_a])=\Psi_{g,1}(\bq_a,x_a)$ and
$(\mg_a,\lgo_a,t_a^*)=\wt{\Psi}_{g,1}(\bq_a,x_a)$, then
$\xi^a(\mg_a,\lgo_a,t_a^*)\to (\mg,\lgo,t^*)$, a random variable with
law $\CLM_{g,1}^{(\beta)}$,
\end{enumerate}
these convergences holding almost-surely as $a\to\infty$ (along the
proper sequence).  We will adopt the notation
$\lgo_a=(w^\eg_a,c^\eg_a,z^\eg_a,\eg\in E(\mg_a))$,
$l^\eg_a=z^\eg_a+w^\eg_a\circ\un{c}^\eg_a$, let $\tau^\eg_a$ be the
duration of $c^\eg_a$, and let $l_a$ be the process defined as in
(\ref{eq:21}) or (\ref{eq:20}). Note that $(2a)^{-1}\tau_{\lgo_a}$ and
$a^{-1}V_{\bq_a}$ both converge to $\tau_\lgo$ as $a\to\infty$, so the
random variable $(\X^{\rm w},\tau_\lgo)$ has law
$\mathscr{S}^{(\beta)}_g$. By the previous convergences, the processes
$\psi^a(l^\eg_a)$ converge to $l^\eg$, so that $\psi^a(l_a)$ converges
a.s.\ as $a\to\infty$ towards the concatenated process $l$ associated
with $(\mg,\lgo,t^*)$.

\begin{lmm}\label{sec:support-diffuseness-2} 
We have
\begin{eqnarray*}
\lefteqn{\int_{a^{-1/4}X_a}\mu_a(\d x)\mu_a(B(x,\eps))}\\ &\leq&
\left(\frac{4a}{\tau_{\lgo_a}}\right)^2\int_{[0,\tau_{\lgo_a}/(2a))^2}\d
  u\d v\ind_{\{|\psi^a(l_a)(\lfloor 2au\rfloor/(2a))-\psi^a(l_a)(\lfloor
    2av\rfloor/(2a))|<\eps/\varsigma\}}+4\tau_{\lgo_a}^{-1}\, .
\end{eqnarray*}
\end{lmm}

\proof With each integer time $i$ corresponds the exploration of a
corner $e(i)=\varphi_{\bm_a}^i(\eg_*(0))$ of $\bm_a$, and each corner
of $\bm_a$ is explored in this way because $\bm_a$ has only one
face. Considering the vertices $e(i)^-$ as vertices of $\bq_a$, we
have $d_a(e(i)^-,e(j)^-)\geq |l_a(i)-l_a(j)|$ by Lemma
\ref{sec:bounds-distances-1}. From this, it follows that
$$\sum_{x,y\in V(\bq_a)}\ind_{\{d_a(x,y)<a^{1/4}\eps\}}\leq
\int_{[0,\tau_{\lgo_a})^2}\d u\d v\ind_{\{|l_a(\lfloor
    u\rfloor)-l_a(\lfloor
    v\rfloor)|<a^{1/4}\eps\}}+2(\tau_{\lgo_a}/2+\chi(g))\, ,$$ where
  the last term amounts from the fact that the selected vertex $x_a$
  is the only vertex of $\bq_a$ not corresponding to a vertex of
  $\bm_a$, and a rough upper bound. We conclude by dividing by
  $(\tau_{\lgo_a}/2+\chi(g))^2=V_{\bq_a}^2$ and a linear change of
  variables.  \cq

\begin{lmm}\label{sec:support-diffuseness-1}
Almost-surely, the level sets $\{ t\in[0, \tau_\lgo):l(t)=x\}$ have
zero Lebesgue measure for every $x\in \R$.
\end{lmm}

\proof From the construction of $l$ as the concatenation of the paths
$l^\eg$, and given $\mg$, it suffices to show the same property for
the processes $l^\eg=z^\eg+w^\eg\circ \un{c}^\eg$ of
(\ref{eq:13}). Under $\CLM_{g,1}$, recall that conditionally on
$(w^\eg,c^\eg,\eg\in E(\mg))$, $(z^\eg(t),0\leq t\leq \tau^\eg)$ is a
Gaussian process with covariance
$\cov(z^\eg(s),z^\eg(t))=\inf_{s\wedge t\leq u\leq s\vee t}c^\eg(u)$,
by definitions of snakes. Therefore, under this conditioned measure,
the laws of pairs $(l^\eg(s),l^\eg(t))$ have densities, and
consequently, for fixed $s,t$,
$$\CLM_{g,1}(\{l^\eg(s)=l^\eg(t),\tau^\eg\geq s\vee t\})=0\, ,$$ and
the same is true for $\CLM_{g,1}^{(\beta)}$ instead of
$\CLM_{g,1}$. By Fubini's theorem, and with a slight abuse of
notation, this implies that for any $T>0$, if $U,V$ are uniform,
independent random variables on $[0,T]$ independent of
$(\mg,\lgo,t_*)$ then it holds that
$$\CLM_{g,1}^{(\beta)}(\{l^\eg(U)=l^\eg(V),\tau^\eg\geq T\})=0\, .$$
But obviously, if some level set $\{l^\eg=x\}$ had positive Lebesgue
measure with positive probability, then the latter probability would
be positive for some $T>0$. \cq

\medskip

We now conclude that $\mu$ is diffuse. Indeed, by the
convergence of $\psi^a(l_a)$ to $l$ and of $\tau_{\lgo_a}/(2a)$ to
$\tau_\lgo\in(0,\infty)$, we obtain from Lemma
\ref{sec:support-diffuseness-2} and dominated convergence that
$$\liminf_{a\to\infty}\int_{a^{-1/4}X_a}\mu_a(\d
x)\mu_a(B(x,\eps))\leq \frac{16}{\tau_\lgo^2}\int_{[0,\tau_\lgo)^2}\d
  u\d v\ind_{\{|l(u)-l(v)|\leq \eps\}}\, .$$ The last integral
  decreases to $\int_{[0,\tau_\lgo)^2}\d u\d v\ind_{\{l(u)=l(v)\}}$
    as $\eps\to 0$, and this is $0$ by Lemma
    \ref{sec:support-diffuseness-1}. Lemma
    \ref{sec:diffuse-measures-1} allows us to conclude.

\subsubsection{Support of $\mu$}\label{sec:support-1}

The setting is the same as in the previous section. For $e\in
E(\bm_a)$ we let $x(e)=e^+$ if $d_{\bm_a}(e^+,V(\bm_a^{\geq
  2}))>d_{\bm_a}(e^-,V(\bm_a^{\geq 2}))$ and $x(e)=e^-$ otherwise
(intuitively, it selects the vertex incident to $e$ that is the
furthest away from $\bm^{\geq 2}$). For $u\in[0,\tau_{\lgo_a})$ we
  let
$$\langle u\rangle_a=\left\{\begin{array}{cl} \lfloor u\rfloor &
  \mbox{if }x(e(\lfloor u\rfloor))=e(\lfloor u\rfloor)^-\\ \lceil
  u\rceil &\mbox{if }x(e(\lfloor u\rfloor))=e(\lfloor
  u\rfloor)^+\end{array}\right.\, .$$

\begin{lmm}\label{sec:support-3}
The image measure of $2^{-1}\d u\ind_{[0,\tau_{\lgo_a})}(u)$ by
  $u\mapsto e(\langle u\rangle_a)^-$ is the measure on
  $V(\bm_a)=V(\bq_a)\setminus\{x_a\}$ assigning mass $1$ on every
  vertex of $V(\bm_a)\setminus V(\mg_a)$, where $\mg_a$ is seen as a
  subgraph of $\bm_a$, and assigning mass $\deg_{\mg_a}(x)/2$ to every
  vertex $x\in V(\mg_a)$, the latter denoting the degree in the graph
  $\mg_a$.
\end{lmm}

\proof Assume $x$ is an element of $V(\bm_a)\setminus V(\bm_a^{\geq
  2})$. Then $e(\langle u\rangle_a)^-=x$ if and only if $i=\lfloor
u\rfloor$ is either the first integer such that $e(i)^+=x$ or the last
integer such that $e(i)^-=x$, corresponding to the two visits of the
parent edge of $x$ in the tree rooted on $V(\bm^{\geq 2})$ containing
$x$. The total mass of such $u$'s is thus $1/2+1/2=1$. 

Next, assume that $x\in V(\bm_a^{\geq 2})\setminus V(\mg_a)$. This
time $e(\langle u\rangle_a)^-=x$ if $\lfloor u\rfloor$ is one of the two
integers $i$ such that $e(i)\in E(\bm^{\geq 2}_a)$ and $e(i)^-=x$.

Finally, if $x\in V(\mg_a)$ then $e(\langle u\rangle_a)^-=x$ if $\lfloor
u\rfloor$ is an integer $i$ such that $e(i)\in E(\bm_a^{\geq 2})$ and
$e(i)^-=x$. There are $\deg_{\mg_a}(x)$ such integers. \cq

\medskip The measure appearing in Lemma \ref{sec:support-3} is of
total mass $\tau_{\lgo_a}/2$, and we let $\tilde{\mu}_a$ be the
probability measure obtained by renormalizing it by this
quantity. Obviously, the Prokhorov distance in $a^{-1/4}\X_a$ between
$\mu_a$ and $\tilde{\mu}_a$ vanishes as $a\to\infty$, so that
$[X_a,a^{-1/4}d_a,\tilde{\mu}_a]$ converges to the same limit
$[X,d,\mu]$ as $[X_a,a^{-1/4}d_a,\mu_a]$.

\begin{lmm}\label{sec:support-2}
It holds that
$$ \min_{x\in
  V(\bq_a)\setminus\{x_a\}}\tilde{\mu}_a(B_{a^{-1/4}X_a}(x,\eps)) \geq
\min_{s\in[0,\tau_{\lgo_a}/(2a))}\frac{2a}{\tau_{\lgo_a}}\int_0^{\tau_{\lgo_a}/(2a)}\d
  u\, \ind_{\{\omega(\psi^a(l_a),|\langle 2au\rangle_a-2as|/(2a))\leq
    \eps/(2\varsigma)\}}\, .
$$
\end{lmm}

\proof Let $x$ be a vertex in $V(\bq_a)\setminus\{x_a\}$, and let
$i_x$ be an integer such that $e(i_x)^-=x$. If $j$ is another integer
time, then $d_a(x,e(j)^-)\leq l_a(j)+l_a(i_x)-2\min_{j\wedge i_x\leq
  u\leq j\vee i_x} l_a(u)\leq 2\omega(l_a,|j-i_x|)$ by Lemma
\ref{sec:bounds-distances-1}. Thus
\begin{eqnarray*}
(\tau_{\lgo_a}/2)\tilde{\mu}_a\left(\{y\in
V(\bq_a):d_a(x,y)<a^{1/4}\eps\}\right) &=&
\frac{1}{2}\int_{0}^{2\tau_{\lgo_a}}\d u\, \ind\{d_a(x,e(\langle
u\rangle_a)^-)<a^{1/4}\eps\} \\
&\geq &\frac{1}{2}\int_0^{2\tau_{\lgo_a}}\d u\, \ind\{\omega(l_a,
|\langle u\rangle_a-i_x|)<a^{1/4}\eps/2\}\, ,
\end{eqnarray*} 
where we used Lemma \ref{sec:support-3} in the first step. The result
is obtained by taking the infimum over $x$ and by simple
rearrangements. \cq

\medskip

Since $\psi^a(l_a)$ converges uniformly to the continuous limit $l$,
there exists $\eta$ such that for every (large enough) $a$,
$\omega(\psi^a(l_a),\eta)\leq \eps/(4\varsigma)$. From this, we
readily obtain from the last lemma that
$$\limsup_{a\to\infty}\inf_{x\in
  V(\bq_a)\setminus\{x_a\}}\tilde{\mu}_a(B_{a^{-1/4}X_a}(x,\eps))\geq
\eta/2>0\, .$$ Since $V(\bq_a)\setminus\{x_a\}$ is $a^{-1/4}$-dense in
$\X_a$, one concludes by Lemma \ref{sec:marked-topologies-2} that
$\mu$ is of full support.

\subsubsection{Uniqueness of geodesics}\label{sec:uniqueness-geodesics}

We now specialize the results of Sect.\ \ref{sec:discr-cont-path} to
the case $k=2$, starting with some general remarks. A delay
$D=[d_1,d_2]\in {\cal D}(\bq,(x_1,x_2))$ is identified with the
quantity $d_1-d_2$, which is an integer in the interval
$(-d_\bq(x_1,x_2),d_\bq(x_1,x_2))$ such that $D+d_\bq(x_1,x_2)$ is
even. 

Let $(\bq,(x_1,x_2),D)\in \bQ_{g,2}$ and
$(\bm,[\bl])=\Psi_{g,2}(\bq,(x_1,x_2),D)$. Consider any chain in $\bq$
from $x_1$ to $x_2$ with length $d_\bq(x_1,x_2)$. There must be a
(unique) vertex $x$ in this chain that belongs to $V_1\cap V_2$ (where
$V_i=V_i(\bq,\bx,D)$), i.e.\ for which
$d_\bq(x_1,x)-d_\bq(x_2,x)=D$. We claim that $\bl(x)=\min_{V_1\cap
  V_2} \bl(y)$. Indeed, we know from Theorem \ref{BIJ} that
$d_\bq(y,x_1)=\bl(y)-\bl(x_1)$ and $d_\bq(y,x_2)=\bl(y)-\bl(x_2)$ for
$y\in V_1\cap V_2$, and since
$d_\bq(x,x_1)+d_\bq(x,x_2)=d_\bq(x_1,x_2)$ since $x$ is on a geodesic,
there must be equality in the triangle inequality
$2\bl(y)-\bl(x_1)-\bl(x_2)\geq d_\bq(x_1,x_2)$, and $\bl(x)$ is indeed
minimal. In passing, we have obtained the formula
\begin{eqnarray}d_\bq(x_1,x_2)&=&
2\min_{e\in f_1:\ov{e}\in
  f_2}\bl(e)-\bl(x_1)-\bl(x_2)\nonumber\\ &=&2\min_{e\in f_1:\ov{e}\in
  f_2}\bl(e)-\min_{e\in f_1}\bl(e)-\min_{e\in
  f_2}\bl(e)-2\nonumber\\ &=&2\min_{e\in f_*:\ov{e}\in
  \ov{f}_*}\bl(e)-\min_{e\in f_*}\bl(e)-\min_{e\in \ov{f}_*}\bl(e)-2\,
,\label{eq:6}
\end{eqnarray}
where $f_*$ is the face of $\bm$ incident to its root $e_*$, and
$\ov{f}_*$ is the other face. Indeed, with a vertex in $V_1\cap V_2$
corresponds a vertex of $\bm$ incident to both $f_1$ and $f_2$, and
such a vertex is the origin of at least one edge incident to $f_1$
with reversal incident to $f_2$. This gives the first equality, and
the second one follows by symmetry in the indices of the second
formula.

For an element $(\mg,\lgo,t^*)\in \bC_{g,2}$, we let 
$\fgo_*$ be the face of $\mg$ incident to $\eg_*$, and $\ov{\fgo}_*$
be the other face.  Let 
$$\eg_i:=\varphi_\mg^{i-1}(\eg_*)\, ,\qquad 1\leq i\leq
\deg_\mg(\fgo_*)$$ be the list of edges of $\mg$ incident to $\fgo_*$
(this notation should not be confused with the notation $\eg(i)$,
which denotes edges of $\bm$).  Besides the process $l$ defined at
(\ref{eq:20}) we will need two other auxiliary processes.  First, let
$i_0$ be the first index $i$ such that $\ov{\eg}_i\in
\ov{\fgo}_*$. Let
$$\eg'_1=\varphi_\mg(\ov{\eg}_{i_0})=
\sigma_\mg^{-1}(\varphi_\mg^{i_0-1}(\eg_*))\, ,$$ which is incident to
$\ov{\fgo}_*$ and has the same origin as $\eg_{i_0}$, and
$$\eg'_i=\varphi_\mg^{i-1}(\eg'_1)\, ,\qquad 1\leq i\leq
\deg_\mg(\ov{\fgo}_*)\, .$$
We define a process $l'$ as the concatenation  
$$l'=l^{\eg'_1}\ldots l^{\eg'_{\deg_\mg(\ov{\fgo}_*)}}\, .$$ 
Next, let $w$ be defined as the concatenation
$$w=w^{\eg_*}w^{\varphi_\mg(\eg_*)}\ldots
w^{\varphi_\mg^{\deg_\mg(\eg_*)-1}(\eg_*)}\, .$$ 
We also consider the function defined for $0\leq t\leq
\tau(w)$ by $\Upsilon(t)=i$ if $\sum_{0\leq
  j<i}r_{\varphi_\mg^j(\eg_*)}\leq t<\sum_{0\leq j\leq
  i}r_{\varphi_\mg^j(\eg_*)}$, and $\Upsilon(\tau(w))=\deg_\mg
(\fgo_*)-1$.

When $(\bm,[\bl])=\Psi_{g,2}(\bq,\bx,D)$ and
$(\mg,\lgo,t^*)=\Xi_{g,2}(\bm,[\bl])$ for some $(\bq,\bx,D)\in
\bQ_{g,2}$, it is straightforward, though a little tedious to check
formally, that $l'$ is alternatively described as the process
interpolating linearly between the values
\begin{equation}\label{eq:8}
l'(i)=\bl(\varphi_\bm^i(e'_1))-\bl(e'_1)\, ,\qquad 0\leq i\leq
\deg_\bm(\ov{f}_*)-1\, ,
\end{equation} and $l'(\deg_\bm(\ov{f}_*))=0$,
where
$$e'_1=\sigma_\bm^{-1}\varphi_\bm^{j_0}(e_*)\, ,$$ and where $j_0$ is
the first index $j$ such that $\ov{\varphi_\bm^{j_0}(e_*)}\in
\ov{f}_*$.  One checks in this case that
\begin{equation}\label{eq:22}
\bl(e'_1)=\sum_{i=1}^{i_0-1}\wh{l}^{\eg_i}\, ,
\end{equation} so that the
function $\sum_{i=1}^{i_0-1}\wh{l}^{\eg_i}+l'$ reproduces faithfully
the labels of half-edges incident to $\ov{f}_*$.  Still in this case,
the function $w$ reproduces the labels $\bl$ along the face of
$\bm^{\geq 2}$ incident to $e_*$, while $\Upsilon(t)$ records the
label of the edge of $\mg$ on which the edge of $\bm^{\geq 2}$
explored at time $\lfloor t\rfloor$ projects. 

\begin{lmm}\label{sec:proof-theor-refs-3}
Let $(\bq,\bx,D)\in \bQ_{g,2}$ and
$(\mg,\lgo,t^*)=\wt{\Psi}_{g,2}(\bq,\bx,D)$.  Assume that $\min
\{w(t):\ov{\eg}_{\Upsilon(t)}\in \ov{\fgo}_*\}$ is attained for a
unique value $i_m$ of $\Upsilon_t$, such that $\ov{\eg}_{i_m}\in
\ov{\fgo}_*$: Also assume that for some $R>0$,
$$\min \{w(t):\ov{\eg}_{\Upsilon(t)}\in \ov{\fgo}_*,\Upsilon(t)\neq
i_m\}>\min \{w(t):\ov{\eg}_{\Upsilon(t)}\in \ov{\fgo}_*\}+R/2\, .$$
Then
\begin{eqnarray*}
\lefteqn{\diam \left({\rm medqg}_{x_1x_2}^{D,R}(X_\bq)\right)}\\ &\leq
& 2\max_{\eg\in E} \omega\left(l^\eg,\diam\{i<
\tau^\eg:w^\eg(\un{c}^\eg(i))\leq \inf w^\eg+R/2\}\right)+R+2\, ,
\end{eqnarray*}
\end{lmm}

\proof Recall that $X_\bq$ is considered as a geodesic metric space, so
that the quantity ${\rm medqg}_{x_1x_2}^{D,R}$ makes sense.  Assume
that $x\in{\rm medqg}_{x_1x_2}^{D,R}(X_\bq)$. We first argue that
there exists an element of ${\rm medqg}_{x_1x_2}^{D,R}(X_\bq)\cap
V(\bq)$ at distance at most $1$ from $x$. Indeed, if $x\notin V(\bq)$,
then it must stand on an edge of the metric graph $\X_\bq$ that links
two vertices $y,z$. Thus any geodesic path from $x$ to $x_1$ or $x_2$
must pass through $y$ or $z$, and by symmetry we may assume
$d_\bq(x,x_1)=d_\bq(x,y)+d_\bq(y,x_1)$. If we also have
$d_\bq(x,x_2)=d_\bq(x,y)+d_\bq(y,x_2)$ then by subtracting
$d_\bq(y,x_1)-d_\bq(y,x_2)=D$ and $y$ is the wanted vertex. Otherwise,
it holds that $d_\bq(x,x_2)=d_\bq(x,z)+d_\bq(z,x_2)$ so that
$$D=d_\bq(x,x_1)-d_\bq(x,x_2)=d_\bq(y,x_1)-d_\bq(z,x_2)+
(d_\bq(x,y)-d_\bq(x,z))\, .$$ As edges have length $1$, we have
$d_\bq(x,y)=1-d_\bq(x,z)>0$, so the quantity in brackets is in
$(-1,1)$, and the others are integers, so we must have
$d_\bq(x,y)=d_\bq(x,z)=1/2$. But since $|d_\bq(z,x_2)-d_\bq(y,x_2)|=
1$, this implies that $d_\bq(y,x_1)-d_\bq(y,x_2)$ has not the same
parity as $D$. This contradicts the fact that $\bl_1-\bl_2$ takes
values in $2\Z$ established in Sect.\ \ref{sec:bijection}. Thus, we
may, and will, assume that $x\in V(\bq)$ up to increasing the diameter
by $2$.

We next argue that there is a vertex $y\in V_1\cap V_2$ at
$d_\bq$-distance at most $R/2$ from $x$. If $x$ itself is not such a
vertex, then it is incident only to edges of $E_i(\bq,\bx,D)$ for some
$i\in\{1,2\}$. Assume $i=1$ by symmetry and note that with the
notations of Sect.\ \ref{sec:bijection}, it holds that
$\bl_1(x)=\bl_2(x)$. Take a chain of edges of $\bq$ starting from $x$
and ending at $x_2$, along which the function $\bl$ decreases, so that
this chain is geodesic. It is easy to check that the first element $e$
of this chain which belongs to $E_2(\bq,\bx,D)$ is such that $y=e^-$
is an element of $V_1\cap V_2$. Therefore, the chain from $x$ to $y$
can be extended as well to a geodesic from $x$ to $x_1$. Thus
$D=d_\bq(x,x_1)-d_\bq(x,x_2)=d_\bq(y,x_1)-d_\bq(y,x_2)$ and $y$ is an
element of of ${\rm med}^D_{x_1x_2}(\X_\bq)\cap V(\bq)$. Since $x\in
{\rm qg}^{R}_{x_1x_2}$ and $y$ is on a geodesic from $x$ to $x_1$ and
from $x$ to $x_2$, it is easily computed that $y$ is another element
of ${\rm qg}^R_{x_1x_2}$ at distance at most $R/2$ from $x$.

Therefore, we may, and will, restrict our attention to $x\in {\rm
  medqg}^{D,R}_{x_1x_2}\cap V_1\cap V_2$, up to increasing the
diameter by $R+2$. Letting $(\bm,[\bl])=\Psi_{g,2}(\bq,\bx,D)$, the
corresponding $x\in V(\bm)$ satisfies
$d_\bq(x,x_1)+d_\bq(x,x_2)=2\bl(x)-\bl(x_1)-\bl(x_2)\leq
d_\bq(x_1,x_2)+R$, and together with (\ref{eq:6}) this yields
$$\bl(x)\leq \min_{V_1\cap V_2}\bl+R/2\, .$$ Now since $x$ is incident
to $f_*$ and $\ov{f}_*$, we can find a time $t$ with $\bl(x)=w(t)$,
and $\ov{\eg}_{\Upsilon(t)}$ is incident to $\ov{\fgo}_*$, and one
sees after a moment's thought that the last inequality can be
rewritten as
$$\bl(x)=w(t)\leq \min\{w(s):\ov{\eg}_{\Upsilon(s)}\in
\ov{\fgo}_*\}+R/2\, .$$ These facts imply by hypothesis that
$\Upsilon(t)=i_m$.  Moreover, by Lemma \ref{sec:bounds-distances-1},
if $y$ is another vertex such as $x$, we have $d_\bq(x,y)\leq
2\omega(l^{e_{i_m}},|p-p'|)$, where $\eg_{i_m}(p)^-=x$ and
$\eg_{i_m}(p')^-=y$ with the notations of Section
\ref{sec:map-reductions-1}. Since $x$ and $y$ are incident to $f_*$
and $\ov{f}_*$, we have $w^\eg(\un{c}^\eg_p)=\bl(x)-\bl(\eg_{i_m}(0))$
and $w^\eg(\un{c}^\eg_{p'})=\bl(y)-\bl(\eg_{i_m}(0))$ which must be
both $\leq \min w^\eg+R/2$. Thus the result. \cq

\medskip

By Theorem \ref{sec:mains-results-1}, along the sequence $(a_n,n\geq
0)$, the spaces $a^{-1/4}\X_\bq^{\rm w}$ under $\QQ_g^{(\beta/a)}$
converge to $[X,d,\mu]$ in distribution, for some limiting $[X,d,\mu]$
with law $\mathscr{S}_g^{(\beta)}$. Together with lemma
\ref{sec:rand-mark-spac-2}, this entails that the marked space
$[X_\bq,a^{-1/4}d_\bq,(x,y)]$ under $\QQ_g^{(\beta/a)}(\d
\bq)\mu_\bq^{\otimes 2}(x,y)$ converges to a space $[X,d,(x,y)]$ with
law $\mathscr{S}_g^{(\beta)}(\d\X)\mu^{\otimes 2}(\d(x,y))$.  On the
other hand, the law of $\wt{\Psi}_{g,2}(\bq,(x,y),D)$ under the law
$\QQ_{g,2}^{(\beta/a)}$ defined by the formula
$$\QQ_g^{(\beta/a)}(\d\bq)\mu_\bq^{\otimes 2}(\d
(x,y))\frac{\#_{{\cal D}(\bq,(x,y))}(\d D)}{|{\cal
    D}(\bq,(x,y))|}=\frac{\QQ_{g,2}\left(V_\bq^{-2}|{\cal
    D}(\bq,(x,y))|^{-1}e^{-\beta V_\bq/a}\,
  \d(\bq,(x,y),D)\right)}{\QQ_{g,1}(V_\bq^{-1}e^{-\beta V_\bq/a})}$$
(recall (\ref{eq:19})) equals
$$\LM^{(\beta/a)}_{g,2}(\d(\mg,\lgo,t^*)):=
\frac{\LM_{g,2}\left((\tau_\lgo/2+\chi(g))^{-2}|{\cal
    D}(\bq,(x,y))|^{-1}e^{-\beta(\tau_\lgo/2+\chi(g))/a}\,
  \d(\mg,\lgo,t^*)\right)}{\LM_{g,1}\left((\tau_\lgo/2+\chi(g))^{-1}
  e^{-\beta(\tau_\lgo/2+\chi(g))/a}\right)}\, .$$ here, we take the
convention $|{\cal D}(\bq,(x,y))|^{-1}=0$ when ${\cal D}(\bq,(x,y))$
is empty. The random variable $|{\cal
  D}(\bq,(x,y))|=(d_\bq(x,y)-1)\vee 0$ is a continuous function of
$(\mg,\lgo,t^*)$, since (\ref{eq:6}), (\ref{eq:8}) and (\ref{eq:22})
and the definition of $w$ imply
\begin{equation}\label{eq:23}
d_\bq(x,y)=2\min_{t:\ov{\eg}_{\Upsilon(t)}\in \ov{\fgo}_*} w(t)
-\min l-\min l'-\sum_{i=1}^{i_0-1}\wh{l}^{\eg_i} -2\, .
\end{equation} 
We let $\mathfrak{d}$ denote the latter quantity with the $-2$
omitted, which can be associated with any element of $\bC_{g,2}$ (not
necessarily in $\bC_{g,2}^{\rm map}$).

\begin{lmm}
The law of $\xi^a(\mg,\lgo,t^*)$ under $\LM^{(\beta/a)}_{g,2}$
converges towards the law $\CLM^{(\beta)}_{g,2}$ defined by 
$$\CLM^{(\beta)}_{g,2}(\d(\mg,\lgo,t^*))=3
\frac{\CLM_{g,2}\left(\tau_\lgo^{-2}\mathfrak{d}^{-1}e^{-\beta
    \tau_\lgo}\,
  \d(\mg,\lgo,t^*)\right)}{\CLM_{g,1}\left(\tau_\lgo^{-1} e^{-\beta
    \tau_\lgo}\right)}\, .$$
\end{lmm}

\proof Let $H$ be bounded continuous on $\bC_{g,2}$. Let $h_\eps$ be
$0$ on $[0,\eps/2]$, equal to $1$ on $[\eps,\infty)$ and interpolate
  linearly in-between.  Let
  $H_\eps(\mg,\lgo,t^*)=H(\mg,\lgo,t^*)h_\eps(\mathfrak{d})h_\eps(\tau_\lgo)$.
  Then a straightforward use of Theorem \ref{sec:cont-meas-label-1}
  entails that $\LM_{g,2}^{(\beta/a)}(H_\eps\circ \xi^a)$ converges
  to $\CLM_{g,2}^{(\beta)}(H_\eps)$. When $H$ is nonnegative, this
  quantity converges as $\eps\to 0$ to $\CLM_{g,2}^{(\beta)}(H)$,
  and in particular $\CLM_{g,2}^{(\beta)}$ is a sub-probability
  measure. Thus, by dominated convergence, this last convergence holds
  for $H$ not necessarily non-negative.

On the other hand, we know that $(\tau_\lgo/2+\chi(g),\mathfrak{d})$
has same law under $\LM_{g,2}^{(\beta/a)}$ as $(V_\bq,d_\bq(x,y)+2)$
under $\QQ_g^{(\beta/a)}(\d \bq)\mu_\bq^{\otimes 2}(\d (x,y))$. So
$(a^{-1}(\tau_\lgo/2+\chi(g)),a^{-1/4}\mathfrak{d})$ under
$\LM_{g,2}^{(\beta/a)}$ converges in distribution to $({\cal
  V},d(x,y))$ under $\mathscr{S}_g^{(\beta)}(\d(\X,{\cal
  V}))\mu^{\otimes 2}(\d(x,y))$ (here we used Proposition
\ref{sec:rand-mark-spac-2} for the convergence of the distance).
Since we know that $\mu$ is a.s.\ diffuse, it holds that $d(x,y)>0$
a.s.\ under the latter measure, entailing the tightness of the
variables $a^{-1/4}\mathfrak{d}$ under $\LM_{g,2}^{(\beta/a)}$ in
$(0,\infty)$. Since $a^{-1}V_\bq=a^{-1}(\tau_\lgo/2+\chi(g))$ under
$\LM_{g,2}^{(\beta/a)}$ are tight in $(0,\infty)$ as well, this
entails that for some positive $\eta(\eps)\to 0$ as $\eps\to 0$,
$$\left|\LM_{g,2}^{(\beta/a)}(H\circ
\xi^a)-\LM_{g,2}^{(\beta/a)}(H_\eps\circ \xi^a)\right|\leq
\|H\|_\infty\eta(\eps)\, .$$ By first taking $\eps$ small enough and
then letting $a\to\infty$, we conclude that
$\LM_{g,2}^{(\beta/a)}(H\circ\xi^a)\to \CLM_{g,2}^{(\beta)}(H)$,
as wanted. \cq

\medskip

As in the beginning of Section \ref{sec:support-diffuseness}, we make
the assumption that we are working with a probability space
$(\Omega,\FF,P)$ under which are defined random variables
$(\bq_a,(x_a,y_a),D_a)$ with respective laws $\QQ_g^{(\beta/a)}(\d
\bq)\mu_\bq^{\otimes 2}(\d (x,y))|{\cal D}(\bq,(x,y))|^{-1}\#_{{\cal
    D}(\bq,(x,y))}(\d D)$, such that, along $a_n,n\geq 1$,
\begin{enumerate}
\item
the marked spaces $[X_a,a^{-1/4}d_a,(x_a,y_a)]$ converge a.s.\ in
$\pmc^*_2$ to $[X,d,(x,y)]$ with law
$\pi^1_*\mathscr{S}_g^{(\beta)}(\d \X)\mm_2(\X,\d (x,y))$, where
$\pi^1:\pmw\times (0,\infty)\to \pmw$ is the first projection,
\item
$a^{-1/4}D_a\to D$ a.s., the latter having uniform law on
  $(-d(x,y),d(x,y))$ given $d(x,y)$, and
\item
letting $\Psi_{g,2}(\bq_a,(x_a,y_a),D_a)=(\bm_a,[\bl_a])$ and
$\wt{\Psi}_{g,2}(\bq_a,(x_a,y_a),D_a)=(\mg_a,\lgo_a, t^*_a)$, we have
a.s.\ $\xi^a(\mg_a,\lgo_a, t^*_a)\to (\mg,\lgo,t^*)$ where the limit
has law $\CLM_{g,2}^{(\beta)}$.
\end{enumerate}

\begin{lmm}\label{sec:proof-theor-refs-5}
The following events are of full measure under $\CLM_{g,2}$:
\begin{itemize}
\item
The quantity $\inf\{w(t):\ov{\eg}_{\Upsilon(t)}\in \ov{\fgo}_*\}$ is
attained for a unique time $t$, which is not a jump time of
$\Upsilon(t)$ and is distinct from $0$ and $\tau(w)$.
\item
The quantities $\diam\{0\leq s\leq \tau^\eg:w^\eg(\un{c}^\eg(s))\leq \inf
w^\eg+\delta\}$ for $\eg\in E$ decrease to $0$ as $\delta\downarrow 0$.
\end{itemize}
\end{lmm}

\proof Under the measure $\CLM_{g,2}$ and conditionally on
$\mg,\br,(\wh{w}^\eg,\eg\in E(\mg))$, the processes $w^\eg$ are
independent Brownian bridges. Since the minimum value attained by a
Brownian bridge has a diffuse law, and is attained a.s.\ once by the
path at a time which is neither the starting nor the ending time of
the bridge, the first point is clear.

Now conditionally on $\mg,\br,(w^\eg,\eg\in E(\mg))$, the paths
$(c^\eg,z^\eg)$ are independent Brownian snakes under $\CLM_{g,2}$. If
$s_m^\eg$ is the (unique) time at which $w^\eg$ attains its minimum,
then it follows from standard properties of Brownian motion that
$\un{c}^\eg$ (the infimum process of $c^\eg$) attains the value $s_m$
at a unique time $t_m$ with probability $1$. This means that
$\un{c}^\eg_t<\un{c}^\eg_{t_m}<\un{c}^\eg_{t'}$ for every $t<t_m<t'$,
and by continuity the two bounding terms converge towards the middle
one as $t,t'\to t_m$. Combining this with the continuity of $w^\eg$
and the fact that the value $\inf w^\eg$ is attained only at
$s^\eg_m$, we obtain the result. \cq

\medskip

We can now finish the proof of Theorem \ref{sec:mains-results-2}. By
the first point in Lemma \ref{sec:proof-theor-refs-5} and the
a.s.\ convergence of $\xi^a(\mg_a,\lgo_a,t^*_a)$ to $(\mg,\lgo,t^*)$,
we conclude that $(\mg_a,\lgo_a,t^*_a)$ satisfies the hypotheses of
Lemma \ref{sec:proof-theor-refs-3} for $R$ in the form $\delta
a^{1/4}$ with $\delta$ small enough, and every large enough $a$,
almost-surely. From the latter lemma and to be able to apply Lemma
\ref{sec:d-interm-points-1} to the spaces
$[X_a,a^{-1/4}d_a,(x_a,y_a)]$ and the delays $D_a$, it suffices to
notice that for every $\eta>0$, there exists $\delta\in(0,\eta)$ such
that for every $a,\eg$,
$$\diam\left(\{i<\tau_a^\eg:w_a^\eg(\un{c}_a^\eg(i))\leq \inf
w_a^\eg+a^{1/4}\delta/2\}\right)\leq 2a\eta\, ,$$ as follows from the
second point of Lemma \ref{sec:proof-theor-refs-5} and the
a.s.\ uniform convergence of the processes
$\phi^a(w^\eg_a),\vartheta^a(c^\eg_a)$ towards $w^\eg,c^\eg$.  Then,
taking $\eps>0$ and choosing $\eta$ in turn so that
$2\varsigma\omega(\psi^a(l_a^\eg),\eta)+\eta< \eps$ for every $a,\eg$,
which is possible by the uniform convergence of $\psi^a(l^\eg_a)$ to
$l^\eg$, Lemma \ref{sec:proof-theor-refs-3} implies that
$$\limsup_{a\to\infty}\diam({\rm
  medqg}^{a^{-1/4}D_a,\delta}_{xy}(a^{-1/4}\X_a))\leq \eps\, .$$
Therefore, Lemma \ref{sec:d-interm-points-1} implies that a.s.\ under
$$\mathscr{S}_g^{(\beta)}(\d(\X,{\cal V}))\mu^{\otimes 2}(\d (x,y))\d
D\ind_{[-d(x,y),d(x,y)]}(D)/2d(x,y)\, ,$$ all geodesic paths from $x$
to $y$ pass through a unique $D$-median point. This shows that
a.s.\ under the measure $\mathscr{S}_g^{(\beta)}(\d(\X,{\cal
  V}))\mu^{\otimes 2}(\d (x,y))$, the geodesics from $x$ to $y$ all
pass through a unique $D$-median point, for Lebesgue-almost every
$D$. Since $D$-median points form a continuous parametrization in
$D\in[-d(x,y),d(x,y)]$ of geodesics from $x$ to $y$, this implies that
a.s.\ the geodesic between $\mu$-almost all $x$ and $y$ is unique,
hence the result.

\def\polhk#1{\setbox0=\hbox{#1}{\ooalign{\hidewidth
  \lower1.5ex\hbox{`}\hidewidth\crcr\unhbox0}}}

\bigskip

\noindent Gr\'egory Miermont\\
DMA, \'Ecole Normale Sup\'erieure\\
45 rue d'Ulm \\
F-75230 Paris Cedex 05 \\ {\tt gregory.miermont@ens.fr}


\begin{thebibliography}{10}

\bibitem{aldouscrt91}
D.~J. Aldous.
\newblock The continuum random tree. {I}.
\newblock {\em Ann. Probab.}, 19(1):1--28, 1991.

\bibitem{ADJ}
J.~Ambj{\o}rn, B.~Durhuus, and T.~Jonsson.
\newblock {\em Quantum geometry. A statistical field theory approach}.
\newblock Cambridge Monographs on Mathematical Physics. Cambridge University
  Press, Cambridge, 1997.

\bibitem{angelpeeling}
O.~Angel.
\newblock Growth and percolation on the uniform infinite planar triangulation.
\newblock {\em Geom. Funct. Anal.}, 13(5):935--974, 2003.

\bibitem{BenCan86}
E.~A. Bender and E.~R. Canfield.
\newblock The asymptotic number of rooted maps on a surface.
\newblock {\em J. Combin. Theory Ser. A}, 43(2):244--257, 1986.

\bibitem{billingsley99}
P.~Billingsley.
\newblock {\em Convergence of probability measures}.
\newblock Wiley Series in Probability and Statistics: Probability and
  Statistics. John Wiley \& Sons Inc., New York, second edition, 1999.
\newblock A Wiley-Interscience Publication.

\bibitem{bourbaki71}
N.~Bourbaki.
\newblock {\em \'{E}l\'ements de math\'ematique. {T}opologie g\'en\'erale.
  {C}hapitres 1 \`a 4}.
\newblock Hermann, Paris, 1971.

\bibitem{BdFGmobiles}
J.~Bouttier, P.~Di~Francesco, and E.~Guitter.
\newblock Planar maps as labeled mobiles.
\newblock {\em Electron. J. Combin.}, 11:Research Paper 69, 27 pp.
  (electronic), 2004.

\bibitem{BoGu08a}
J.~Bouttier and E.~Guitter.
\newblock Statistics in geodesics in large quadrangulations.
\newblock {\em J. Phys. A}, 41(14):145001, 30, 2008.

\bibitem{burago01}
D.~Burago, Y.~Burago, and S.~Ivanov.
\newblock {\em A course in metric geometry}, volume~33 of {\em Graduate Studies
  in Mathematics}.
\newblock American Mathematical Society, Providence, RI, 2001.

\bibitem{ChMaSc}
G.~Chapuy, M.~Marcus, and G.~Schaeffer.
\newblock A bijection for rooted maps on orientable surfaces.
\newblock 2007.
\newblock arXiv:0712.3649.

\bibitem{CSise}
P.~Chassaing and G.~Schaeffer.
\newblock Random planar lattices and integrated super{B}rownian excursion.
\newblock {\em Probab. Theory Related Fields}, 128(2):161--212, 2004.

\bibitem{dudley02}
R.~M. Dudley.
\newblock {\em Real analysis and probability}, volume~74 of {\em Cambridge
  Studies in Advanced Mathematics}.
\newblock Cambridge University Press, Cambridge, 2002.
\newblock Revised reprint of the 1989 original.

\bibitem{duq02}
T.~Duquesne.
\newblock A limit theorem for the contour process of conditioned
  {G}alton-{W}atson trees.
\newblock {\em Ann. Probab.}, 31(2):996--1027, 2003.

\bibitem{duqleg02}
T.~Duquesne and J.-F. Le~Gall.
\newblock Random trees, {L}\'evy processes and spatial branching processes.
\newblock {\em Ast\'erisque}, 281:vi+147, 2002.

\bibitem{duqlegprep}
T.~Duquesne and J.-F. Le~Gall.
\newblock Probabilistic and fractal aspects of {L}\'evy trees.
\newblock {\em Probab. Theory Related Fields}, 131(4):553--603, 2005.

\bibitem{evpiwin}
S.~N. Evans, J.~Pitman, and A.~Winter.
\newblock Rayleigh processes, real trees, and root growth with re-grafting.
\newblock {\em Probab. Theory Related Fields}, 134(1):81--126, 2006.

\bibitem{evanswinter}
S.~N. Evans and A.~Winter.
\newblock Subtree prune and regraft: a reversible real tree-valued {M}arkov
  process.
\newblock {\em Ann. Probab.}, 34(3):918--961, 2006.

\bibitem{FlSe}
P.~Flajolet and R.~Sedgewick.
\newblock {\em Analytic Combinatorics}.
\newblock Cambridge University Press, Cambridge, 2009.

\bibitem{fukaya}
K.~Fukaya.
\newblock Collapsing of {R}iemannian manifolds and eigenvalues of {L}aplace
  operator.
\newblock {\em Invent. Math.}, 87(3):517--547, 1987.

\bibitem{gao93}
Z.~Gao.
\newblock A pattern for the asymptotic number of rooted maps on surfaces.
\newblock {\em J. Combin. Theory Ser. A}, 64(2):246--264, 1993.

\bibitem{GPW06}
A.~Greven, P.~Pfaffelhuber, and A.~Winter.
\newblock Convergence in distribution of random metric measure spaces
  ({$\Lambda$}-coalescent measure trees).
\newblock 2006.

\bibitem{gromov99}
M.~Gromov.
\newblock {\em Metric structures for {R}iemannian and non-{R}iemannian spaces},
  volume 152 of {\em Progress in Mathematics}.
\newblock Birkh\"auser Boston Inc., Boston, MA, 1999.

\bibitem{janmarck05}
S.~Janson and J.-F. Marckert.
\newblock Convergence of discrete snakes.
\newblock {\em J. Theor. Probab.}, 18(3):615--645, 2005.

\bibitem{LaZv04}
S.~K. Lando and A.~K. Zvonkin.
\newblock {\em Graphs on surfaces and their applications}, volume 141 of {\em
  Encyclopaedia of Mathematical Sciences}.
\newblock Springer-Verlag, Berlin, 2004.

\bibitem{legall93}
J.-F. Le~Gall.
\newblock The uniform random tree in a {B}rownian excursion.
\newblock {\em Probab. Theory Relat. Fields}, 96(3):369--383, 1993.

\bibitem{legall06}
J.-F. Le~Gall.
\newblock The topological structure of scaling limits of large planar maps.
\newblock {\em Invent. Math.}, 169(3):621--670, 2007.

\bibitem{legall08}
J.-F. Le~Gall.
\newblock Geodesics in large planar maps and in the {B}rownian map.
\newblock 2008.
\newblock arXiv:0804.3012.

\bibitem{lgp}
J.-F. Le~Gall and F.~Paulin.
\newblock Scaling limits of bipartite planar maps are homeomorphic to the
  2-sphere.
\newblock {\em Geom. Funct. Anal.}, 18(3):893--918, 2008.

\bibitem{therryYM}
T.~L{\'e}vy.
\newblock Yang-{M}ills measure on compact surfaces.
\newblock {\em Mem. Amer. Math. Soc.}, 166(790):xiv+122, 2003.

\bibitem{jfmgm05}
J.-F. Marckert and G.~Miermont.
\newblock Invariance principles for random bipartite planar maps.
\newblock {\em Ann. Probab.}, 35(5):1642--1705, 2007.

\bibitem{MMsnake}
J.-F. Marckert and A.~Mokkadem.
\newblock States spaces of the snake and its tour---convergence of the discrete
  snake.
\newblock {\em J. Theoret. Probab.}, 16(4):1015--1046 (2004), 2003.

\bibitem{MM05}
J.-F. Marckert and A.~Mokkadem.
\newblock Limit of normalized random quadrangulations: the {B}rownian map.
\newblock {\em Ann. Probab.}, 34(6):2144--2202, 2006.

\bibitem{MaSc01}
M.~Marcus and G.~Schaeffer.
\newblock Une bijection simple pour les cartes orientables.
\newblock 2001.
\newblock Preprint available via\\ {\tt
  www.lix.polytechnique.fr/Labo/Gilles.Schaeffer/}.

\bibitem{mierinv}
G.~Miermont.
\newblock An invariance principle for random planar maps.
\newblock In {\em Fourth Colloquium on Mathematics and Computer Sciences
  CMCS'06}, Discrete Math. Theor. Comput. Sci. Proc., AG, pages 39--58
  (electronic). Nancy, 2006.

\bibitem{miergwmulti}
G.~Miermont.
\newblock Invariance principles for spatial multitype {G}alton-{W}atson trees.
\newblock {\em Ann. Inst. H. Poincar{\'e} Probab. Statist.}, 44(6):1128--1161,
  2008.

\bibitem{mierweill}
G.~Miermont and M.~Weill.
\newblock Radius and profile of random planar maps with faces of arbitrary
  degrees.
\newblock {\em Electron. J. Probab.}, 13:no. 4, 79--106, 2008.

\bibitem{MR1802530}
A.~Okounkov.
\newblock Random matrices and random permutations.
\newblock {\em Internat. Math. Res. Notices}, (20):1043--1095, 2000.

\bibitem{OkPan}
A.~Okounkov and R.~Pandharipande.
\newblock Gromov-{W}itten theory, {H}urwitz numbers, and matrix models, {I}.
\newblock 2001.
\newblock arXiv:math.AG/0101147.

\bibitem{petrov75}
V.~V. Petrov.
\newblock {\em Sums of independent random variables}.
\newblock Springer-Verlag, New York, 1975.
\newblock Translated from the Russian by A. A. Brown, Ergebnisse der Mathematik
  und ihrer Grenzgebiete, Band 82.

\bibitem{pitmancsp02}
J.~Pitman.
\newblock {\em Combinatorial stochastic processes}, volume 1875 of {\em Lecture
  Notes in Mathematics}.
\newblock Springer-Verlag, Berlin, 2006.
\newblock Lectures from the 32nd Summer School on Probability Theory held in
  Saint-Flour, July 7--24, 2002, With a foreword by Jean Picard.

\bibitem{revyor}
D.~Revuz and M.~Yor.
\newblock {\em Continuous martingales and {B}rownian motion}, volume 293 of
  {\em Grundlehren der Mathematischen Wissenschaften}.
\newblock Springer-Verlag, Berlin, third edition, 1999.

\bibitem{schaeffer98}
G.~Schaeffer.
\newblock {\em Conjugaison d'arbres et cartes combinatoires al{\'e}atoires}.
\newblock PhD thesis, Universit{\'e} Bordeaux {I}, 1998.

\bibitem{shefGFF}
S.~Sheffield.
\newblock Gaussian free fields for mathematicians.
\newblock {\em Probab. Theory Related Fields}, 139(3-4):521--541, 2007.

\bibitem{villani09}
C.~Villani.
\newblock {\em Optimal transport, old and New}, volume 338 of {\em Grundlehren
  der Mathematischen Wissenschaften [Fundamental Principles of Mathematical
  Sciences]}.
\newblock Springer-Verlag, Berlin, 2009.

\bibitem{weill06}
M.~Weill.
\newblock Asymptotics for rooted planar maps and scaling limits of two-type
  spatial trees.
\newblock {\em Electron. J. Probab.}, 12:Paper no.\ 31, 862--925 (electronic),
  2007.

\end{thebibliography}
\end{document}